\documentclass[UTF-8,reqno]{amsart}
\usepackage{enumerate}
\usepackage{amssymb,url,color, booktabs}
\usepackage{mathrsfs}

\usepackage{color}
\usepackage[colorlinks=true]{hyperref}
\hypersetup{
    linkcolor=blue,          
    citecolor=red,        
    filecolor=blue,      
    urlcolor=cyan
}

\definecolor{MyDarkBlue}{cmyk}{0.8,0.3,0.8,0.4}
\definecolor{yellow}{rgb}{0.99,0.99,0.70}
\definecolor{white}{rgb}{1.0,1.0,1.0}
\definecolor{black}{rgb}{0.00,0.00,0.00}

\numberwithin{equation}{section}

\newcommand{\be}{\begin{eqnarray}}
\newcommand{\ee}{\end{eqnarray}}
\newcommand{\ce}{\begin{eqnarray*}}
\newcommand{\de}{\end{eqnarray*}}
\newtheorem{theorem}{Theorem}[section]
\newtheorem{lemma}[theorem]{Lemma}
\newtheorem{remark}[theorem]{Remark}
\newtheorem{definition}[theorem]{Definition}
\newtheorem{proposition}[theorem]{Proposition}
\newtheorem{Examples}[theorem]{Example}
\newtheorem{corollary}[theorem]{Corollary}

\def\nor{|\mspace{-3mu}|\mspace{-3mu}|}

\def\eps{\varepsilon}

\def\e{\mathrm{e}}

\def\p{\partial}

\def\[{{\Big[}}
\def\]{{\Big]}}
\def\<{{\langle}}
\def\>{{\rangle}}
\def\({{\Big(}}
\def\){{\Big)}}

\def\bx{{\mathbf{x}}}
\def\tr{\mathrm {tr}}

\def\dif{{\mathord{{\rm d}}}}

\def\bbp{{\boldsymbol{p}}}
\def\bbr{{\boldsymbol{r}}}
\def\bbq{{\boldsymbol{q}}}
\def\bba{{\boldsymbol{a}}}
\def\bbw{{\boldsymbol{w}}}
\def\bb2{{\boldsymbol{2}}}

\def\no{\nonumber}
\def\={&\!\!=\!\!&}

\def\bB{{\mathbf B}}

\def\cA{{\mathcal A}}
\def\cB{{\mathcal B}}

\def\cD{{\mathcal D}}

\def\cG{{\mathcal G}}

\def\cI{{\mathcal I}}

\def\cL{{\mathcal L}}
\def\cM{{\mathcal M}}

\def\cP{{\mathcal P}}

\def\cR{{\mathcal R}}
\def\cS{{\mathcal S}}

\def\cW{{\mathcal W}}

\def\mB{{\mathbb B}}
\def\mC{{\mathbb C}}
\def\mD{{\mathbb D}}
\def\mE{{\mathbb E}}

\def\mH{{\mathbb H}}
\def\mI{{\mathbb I}}

\def\mL{{\mathbb L}}
\def\mM{{\mathbb M}}
\def\mN{{\mathbb N}}

\def\mP{{\mathbb P}}

\def\mR{{\mathbb R}}

\def\mW{{\mathbb W}}

\def\bB{{\mathbf B}}
\def\bP{{\mathbf P}}

\def\bE{{\mathbf E}}

\def\1{{\mathbf{1}}}

\def\sA{{\mathscr A}}

\def\sF{{\mathscr F}}

\def\sI{{\mathscr I}}

\def\sM{{\mathscr M}}

\def\sV{{\mathscr V}}

\def\geq{\geqslant}
\def\leq{\leqslant}

\def\div{\mathord{{\rm div}}}

\def\eps{\varepsilon}

\def\e{\mathrm{e}}

\def\p{\partial}

\def\[{{\Big[}}
\def\]{{\Big]}}
\def\<{{\langle}}
\def\>{{\rangle}}
\def\({{\Big(}}
\def\){{\Big)}}

\def\bx{{\mathbf{x}}}
\def\tr{\mathrm {tr}}

\def\dif{{\mathord{{\rm d}}}}

\def\no{\nonumber}
\def\={&\!\!=\!\!&}
\def\bt{\begin{theorem}}
\def\et{\end{theorem}}
\def\bl{\begin{lemma}}
\def\el{\end{lemma}}
\def\br{\begin{remark}}
\def\er{\end{remark}}
\def\bx{\begin{Examples}}
\def\ex{\end{Examples}}
\def\bd{\begin{definition}}
\def\ed{\end{definition}}
\def\bp{\begin{proposition}}
\def\ep{\end{proposition}}
\def\bc{\begin{corollary}}
\def\ec{\end{corollary}}

\def\geq{\geqslant}
\def\leq{\leqslant}

\def\div{\mathord{{\rm div}}}

\def\bP{{\mathbf P}}

\def\<{\langle} \def\>{\rangle}

\allowdisplaybreaks

\begin{document}

\title[Second order SDEs and kinetic FPK equations]
{Second order McKean-Vlasov SDEs and kinetic Fokker-Planck-Kolmogorov equations}

\author{Xicheng Zhang}

\address{Xicheng Zhang:
School of Mathematics and Statistics, Wuhan University,
Wuhan, Hubei 430072, P.R.China\\
Email: XichengZhang@gmail.com
 }

\thanks{
This work is partially supported by NNSFC grants of China (Nos. 12131019, 11731009), and the German Research Foundation (DFG) through the Collaborative Research Centre(CRC) 1283 ``Taming uncertainty and profiting from randomness and low regularity in analysis, stochastics and their applications".
}

\begin{abstract}
In this paper we study second order stochastic differential equations with measurable and density-distribution dependent coefficients.
Through establishing a maximum principle for kinetic Fokker-Planck-Kolmogorov equations with distribution-valued inhomogeneous term, we show the 
existence of weak solutions under mild assumptions. Moreover, by using the H\"older regularity estimate  
obtained recently in \cite{GIMV19}, we also show the well-posedness of generalized martingale problems when diffusion coefficients only depend on the position variable (not necessarily continuous). 
Even in the non density-distribution dependent case, 
it seems that this is the first result about the well-posedness of SDEs with measurable diffusion coefficients.

\bigskip
\noindent 
\textbf{Keywords}: Maximum principle,  De-Giorgi's iteration, Stochastic differential equation, Krylov's estimate, kinetic Fokker-Planck-Kolmogorov equation.\\

\noindent
 {\bf AMS 2010 Mathematics Subject Classification:} 60H10, 35H10
\end{abstract}

\maketitle \rm

\tableofcontents

\section{Introduction}

In this paper we are concerned with the following second order stochastic differential equation
with density-distribution dependent coefficients (also called McKean-Vlasov SDE or simply DDSDE): 
\begin{align}\label{SDE10}
\dif \dot X_t=b_Z(t,Z_t)\dif t+\sigma_Z(t,X_t)\dif W_t,
\end{align}
where $Z_t:=(X_t,\dot X_t)$ and $\dot X_t=:V_t$ 
stands for the velocity, 
$(W_t)_{t\geq 0}$ is a $d$-dimensional standard Brownian motion defined on some stochastic basis $(\Omega,\sF,\bP; (\sF_t)_{t\geq 0})$, and
\begin{align}\label{BB8}
b_Z(t,z):=\int_{\mR^{2d}}b(t,z, \rho_{Z_t}(z), z')\mu_{Z_t}(\dif z'),
\end{align}
and
\begin{align}\label{BB9}
\sigma_Z(t,x):=\sqrt{2a_Z(t,x)},\ \ a_Z(t,x):=\int_{\mR^d}a(t,x, \rho_{X_t}(x), z')\mu_{Z_t}(\dif z').
\end{align}
Here $\mu_{Z_t}(\dif z)=\bP\circ Z^{-1}_t(\dif z)=\rho_{Z_t}(z)\dif z$, $\bP\circ X^{-1}_t(\dif x)=\rho_{X_t}(x)\dif x$, and 
$$
b(t,z,r,z'): \mR_+\times\mR^{2d}\times\mR_+\times\mR^{2d}\to\mR^{d}
$$
and
$$
a(t,x,r,z'): \mR_+\times\mR^{d}\times\mR_+\times\mR^{2d}\to\mM^{d}_{\rm sym}
$$ 
are Borel measurable functions, $\mM^{d}_{\rm sym}$ is the set of all symmetric $d\times d$-matrices. 
Throughout the paper, the density dependence means that $b_Z, a_Z$ point-wisely depend on the density of $Z$, while the distribution dependence means that
$b_Z, a_Z$  non-locally depend on the distribution of $Z$ in the whole space.
Below we suppose that $a$ satisfies that for some $0<\kappa_0<\kappa_1$ and all $t,x,r,z'$,
\begin{align}\label{ELL}
\kappa_0|\xi|^2\leq\xi\cdot a(t,x,r,z')\xi\leq\kappa_1|\xi|^2,\ \ \xi\in\mR^d,
\end{align}
where the dot stands for the inner product of two vectors.
Since the coefficients depend on the density and distribution of the solution itself, equation \eqref{SDE10} is usually 
regarded as a strong nonlinear SDE. Suppose that DDSDE \eqref{SDE10} has a solution, whose meaning is given in Definition \ref{Def2} below. 
By It\^o's formula, it is easy to see that the density $\rho(t,z):=\rho_{Z_t}(z)$ solves the following nonlinear kinetic Fokker-Planck-Kolmogorov equation (abbreviated as FPKE) in the distributional sense:
\begin{align}\label{FPK}
\p_t \rho=\p_{v_i}\p_{v_j}(\bar a_{ij}(\rho) \rho)-v\cdot\nabla_x \rho+\div_v(\bar b(\rho) \rho),
\end{align}
where for $(t,x)\in\mR_+\times\mR^d$,
$$
\bar a(\rho; t,x):=\int_{\mR^{2d}}a(t,x, \<\rho\>(t,x), z')\rho(t,z')(\dif z'),
$$
and  for $(t,z)\in\mR_+\times\mR^{2d}$,
$$
\bar b(\rho; t,z):=\int_{\mR^{2d}}b(t,z, \rho(t,z), z')\rho(t,z')(\dif z').
$$
Here and below, $\<\rho\>(t,x):=\int_{\mR^d}\rho(t,x,v)\dif v$ stands for the mass density, and we use the convention that repeated indices are summed automatically.

\medskip

An important prototype of PDE \eqref{FPK} is the following Landau equation or its variant form (see \cite{Lan36, Li94, AV04, Si16}):
\begin{align}\label{LAN1}
\p_t\rho+v\cdot\nabla_x\rho=\p_{v_i}\left\{\int_{\mR^d}a_{ij}(v-v') \left[\rho(v')\p_{v_j}\rho(v)-\rho(v)\p_{v'_j}\rho(v')\right]\dif v'\right\},
\end{align}
where $\rho=\rho(t,x,v)$ and for $v=(v_1,\cdots,v_d)$,
$$
a_{ij}(v):=c_{d,\gamma}[\delta_{ij}-v_iv_j/|v|^2]/|v|^{\gamma+d},\ \ \gamma\in[-d,+\infty),\ c_{d,\gamma}>0,
$$
which models the behavior of a dilute plasma interacting through binary collisions (see \cite{Sc}).
While nonlinear SDE \eqref{SDE10} with distribution dependent coefficients
was firstly proposed by McKean \cite{Mc} to give a probabilistic explanation for nonlinear Vlasov equations.
Nowadays, McKean-Vlasov SDEs have beed widely used in the study of interacting particle systems and mean-field games 
(cf. \cite{Szn91, CD} and references therein).
Recently, there are many works to study the well-posedness of  first order nondegenerate DDSDEs with rough coefficients.
When $b$ and $a$ does not depend on the density variable $r$ and of at most linear growth, by using the classical Krylov estimates, 
Mishura and Veretennikov \cite{MV} showed the existence of weak solutions for first order DDSDEs. 
The uniqueness is also proved when $a$ does not depend on $z'$ and is Lipschitz continuous in $x$.
Later, in \cite{RZ21} their results were extended to the case that drift $b$ is in $L^p$-spaces. More results and references about McKean-Vlasov SDEs 
are summarized in the paper of \cite{HRW20} (see also the references in \cite{RZ21}). 
For nonlinear SDEs with density dependent coefficients (also called McKean-Vlasov SDEs of Nemytskii-type), 
it was firstly studied in \cite{BR18} to give a probabilistic representation for the solution of nonlinear FPKEs.
In \cite{BR20}, for a large class of time independent coefficients $b$, $a$, Barbu and R\"ockner 
obtained the existence of weak solutions for first order density dependent SDEs. The strategy in \cite{BR18} and \cite{BR20} 
is to solve the associated nonlinear FPKE and then by the well- known superposition principle 
to establish the existence of a weak solution to the first order DDSDE. 
In \cite{HRZ21}, Hao, R\"ockner and the present author give a purely probabilistic 
proof for the existence of weak solutions for the first order density dependent SDE by Euler's scheme.
Our current degenerate kinetic settings are more general and is not studied in the literature yet.

\medskip

As usual, to study DDSDE \eqref{SDE10}, we must have better understanding for the associated kinetic FPKEs. In particular, we need to develop the hypoelliptic regularity estimates for kinetic FPKEs with
rough coefficients. 
In \cite{PP04}, Pascucci and Polidoro obtained the upper bound estimate of weak solutions for a class of ultraparabolic equations of divergence form via Moser's method, 
especially, including the following kinetic equation
\begin{align}\label{FKP0}
\p_t u=\div_v(a\cdot\nabla_v u)+v\cdot\nabla_x u.
\end{align}
In \cite{WZ09, WZ11}, Wang and Zhang showed the H\"older regularity of weak solutions to the above equation.
In \cite{GIMV19}, Golse, Imbert, Mouhot and Vasseur showed the Harnack inequality for \eqref{FKP0} with bounded first order and inhomogeneous terms.
As an application, they also showed the conditional H\"older regularity for the Landau equation \eqref{LAN1}. More recently, weak Harnack inequality and some quantitative estimates 
for kinetic FPKEs were obtained in \cite{GI21, GM21}. More historical backgrounds about Harnack inequalities and kinetic equations are referred to \cite[Section 1.2]{GI21}.

\subsection{Main results}
Let $\cP(\mR^{2d})$ be the space of all probability measures over $\mR^{2d}$.
We first introduce the following notion of weak solutions to DDSDE \eqref{SDE10}.
\bd\label{Def2}
Let $\nu\in\cP(\mR^{2d})$ and ${\frak F}:=(\Omega,\sF,\bP;(\sF_t)_{t\geq 0})$ be a stochastic basis, $Z_t=(X_t,V_t)$ and $W_t$
be $\mR^{2d}$ and $\mR^d$-valued $\sF_t$-adapted processes, respectively.
We call $({\frak F}, Z,W)$ a weak solution of DDSDE \eqref{SDE10} with initial distribution $\nu$ if 
\begin{enumerate}[(i)]
\item $\bP\circ Z^{-1}_0=\nu$, and for Lebesgue almost all $t\geq 0$, 
$$
\bP\circ Z^{-1}_t(\dif z)=\rho_{Z_t}(z)\dif z,\ \ \bP\circ X^{-1}_t(\dif x)=\rho_{X_t}(x)\dif x.
$$
\item $W$ is a standard $d$-dimensional $\sF_t$-Brownian motion.
\item For all $t\geq 0$, it holds that $X_t=X_0+\int^t_0V_s\dif s$ and
$$
V_t=V_0+\int^t_0b_Z(s, Z_s)\dif s+\int^t_0\sigma_Z(s,X_s)\dif W_s,\ \bP-a.s.,
$$
where $b_Z$ and $\sigma_Z$ are defined by \eqref{BB8} and \eqref{BB9}, respectively.
\end{enumerate}
\ed

Before stating our main results, we make the following assumptions about $b$ and $a$:
\begin{enumerate}[{\bf (H$_1$)}]
\item For any $m\in\mN$  and bounded domain $Q\subset\mR_+\times\mR^{2d}\times\mR^{2d}$, it holds that
\begin{align}\label{CON1}
\lim_{h\to 0}\left\|\sup_{r,r'\leq m,|r-r'|\leq h}|b(\cdot,\cdot,r,\cdot)-b(\cdot,\cdot,r',\cdot)|\right\|_{L^1(Q)}=0,
\end{align}
and for all $(t,z,r,z')\in\mR_+\times\mR^{2d}\times\mR_+\times\mR^{2d}$,
\begin{align}\label{DR1}
|b(t,z,r,z')|\leq h(t,z-z')\mbox{ with } \ \nor h\nor_{\widetilde\mL^{q_1}_t(\widetilde\mL^{\bbp_1}_z)}\leq\kappa_2,
\end{align}
where $q_1\in(2,4)$ and $\bbp_1\in(2,\infty)^{2d}$ satisfy $\bba\cdot\frac1{\bbp_1}+\frac2{q_1}<1$, 
the localized norm 
$\nor\cdot\nor_{\widetilde\mL^{q_1}_t(\widetilde\mL^{\bbp_1}_z)}$ is defined by \eqref{LOC1} below.
Here $\bba=(3,\cdots,3,1,\cdots,1)\in\mR^{2d}$ (see \eqref{AA11} below).
\item In addition to ellipticity assumption \eqref{ELL}, we assume that for each $(t,x,z)\in\mR_+\times\mR^{d}\times\mR^{2d}$, $[0,\infty)\ni r\mapsto a(t,x,r,z)\in\mM^d_{\rm sym}$
is continuous.
\end{enumerate}
\br\rm
The continuity of $[0,\infty)\ni r\mapsto a(t,x,r,z)\in\mM^d_{\rm sym}$ together with the boundedness of $a$ implies that
for any bounded $Q\subset \mR_+\times\mR^{d}\times\mR^{2d}$,
\begin{align}\label{CON11}
\lim_{h\to 0}\left\|\sup_{r,r'\leq m,|r-r'|\leq h}|a(\cdot,\cdot,r,\cdot)-a(\cdot,\cdot,r',\cdot)|\right\|_{L^1(Q)}=0.
\end{align}
Indeed, it follows by the dominated convergence theorem and for each $(t,x,z')$,
$$
\lim_{h\to 0}\sup_{r,r'\leq m,|r-r'|\leq h}|a(t,x,r,z')-a(t,x,r',z')|=0.
$$
\er

First of all, we have the following existence result of weak solutions.
 \bt\label{Th55}
 Under {\bf (H$_1$)} and {\bf (H$_2$)}, for any $\nu\in\cP(\mR^{2d})$, 
 there exists at least one weak solution to SDE \eqref{SDE10} with initial distribution $\nu$ in the sense of Definition \ref{Def2}.
 In particular, the probability distributional density $\rho_t$ of $Z_t$ solves \eqref{FPK} with initial value being the probability measure $\nu$.
 Moreover, $\rho$ enjoys the following regularity: for any $\alpha\in(0,1)$ and $(q,\bbp)\in(1,\infty)^{1+2d}$ satisfying
$$
\tfrac2{q}<1+\alpha,\ \  \ \tfrac2{q}+\bba\cdot(\tfrac1{\bbp}-\1)>2\alpha,
$$
it holds that
\begin{align}\label{DL131}
\|\rho\1_{[0,T]}\|_{\mL^{q}_t(\bB^{\alpha}_{\bbp;\bba})}<\infty,\ T>0,
\end{align}
where $\bB^{\alpha}_{\bbp;\bba}$ is the anisotropic Besov space in Definition \ref{bs}.
 \et

\br\label{Re14}\rm
(i) Suppose that $b$ and $a$ do not depend on the density variable $r$, then we can drop the assumption $q_1\in(2,4)$ in  {\bf (H$_1$)}
 but require $\tfrac1{\bbp_1}+\frac{\1}{q_1}<\tfrac1{\bb2}$ (see Remark \ref{Re611} below).  

(ii) We explain the use of 
multi-integrability index $\bbp\in(2,\infty)^{2d}$ (see \eqref{AM9} below). Consider the following second order SDE in $\mR^2$:
$$
\left\{
\begin{aligned}
\dif \dot X^1_t=b_1(X^1_t-X^2_t)\dif t+\dif W^1_t,\\
\dif \dot X^2_t=b_2(X^1_t-X^2_t)\dif t+\dif W^2_t,
\end{aligned}
\right.
$$
where $b_1,b_2:\mR\to\mR$ are two measurable functions.
If we write 
$$
X_t:=(X^1_t,X^2_t),\ \ Z_t:=(X_t, \dot X_t),\ \ W_t:=(W^1_t,W^2_t)
$$
 and
$$
b(x_1,x_2, v_1,v_2):=(b_1,b_2)(x_1-x_2):\mR^4\to\mR^2,
$$ 
then the above SDE can be written in the form of \eqref{SDE10}.
Suppose now that $b_1,b_2\in\widetilde\mL^p$ for some $p>3$. Then by the definition \eqref{LOC101} below,
it is easy to see that \eqref{DR1} holds for $\bbp=(p,p',p',p')$, where $p'$ is chosen large enough so that
$\bba\cdot\frac1{\bbp}=\frac{3}{p}+\frac3{p'}+\frac2{p'}<1$. However, if we do not use the above multi-index $\bbp$, then we have to require $p>8$.
To the author's knowledge, Ling and Xie \cite{LX21} firstly used the mixed-$L^\bbp$ norm for studying singular first order SDEs. We emphasize that
such a feature naturally appears in the study of interacting particle systems with singular interaction force, especially, propagation of chaos (cf. \cite{Szn91}).
We shall study this in a future work. 
\er

Although we have shown the existence of weak solutions under very mild assumptions in the above theorem,
the uniqueness is a more subtle problem.  In fact, in the non-distribution dependent case, that is, 
$a$ and $b$ does not depend on $(r,z')$, when $a$ is uniformly continuous and $b$ is H\"older continuous or in some $L^p$-spaces, 
the uniqueness of weak solutions or strong solutions for second order SDEs was obtained in \cite{Ch17, WZ16, Zh18, CM17} via Zvonkin's transformation. 
Recently, when $a$ is the identity matrix and $b=b_1+b_2$, where $b_1(t,z)$ is a distribution of $z$ and $b_2(t,z,z')$ is bounded measurable, 
the well-posedness of  generalized martingale problem associated with SDE \eqref{SDE10} was obtained in \cite{HZZZ21},
where the generalized martingale problem is taken in the sense of Either and Kurtz \cite[Chapter 4]{EK86}. It should be noted that the notion of generalized martingale solutions
strongly depends on the solvability of the associated backward Kolmogorov equations (see Definition \ref{MP81} below for a precise description).
To show the uniqueness, we now suppose that
\begin{enumerate}[{\bf (H$_3$)}] 
\item $a=a(t,x)$ is independent of $(r,z')$ and satisfies \eqref{ELL}, and $b$ is bounded measurable and satisfies that for 
some $C>0$ and all $(t,z,z')\in\mR_+\times\mR^{2d}\times\mR^{2d}$ and $r,r'\geq 0$,
$$
|b(t,z,r,z')-b(t,z,r',z')|\leq C |r-r'|.
$$
\end{enumerate}

We have the following well-posedness result of generalized martingale problems.
\bt\label{Th15}
Under {\bf (H$_3$)}, for any $s\geq 0$ and $\nu(\dif z)=\rho_0(z)\dif z$, where $\rho_0\in C^1_b(\mR^{2d})$, 
there is a unique generalized martingale solution $\mP\in \widetilde\cM^{a,b}_{s,\nu}$ with initial distribution $\nu$ at time $s$
in the sense of Definition \ref{MP81} below. Moreover, its density $\rho_t(z)$ enjoys the regularity \eqref{DL131} and solves the following nonlinear kinetic FPKE in the distributional sense:
\begin{align}\label{FPK1}
\p_t \rho_t=\p_{v_i}\p_{v_j}(a_{ij}\rho_t)-v\cdot\nabla_x \rho_t+\div_v(b_\rho\rho_t),\ t\geq s,
\end{align}
where $b_\rho(t,z):=\int_{\mR^{2d}}b(t,z,\rho_t(z),z')\rho_t(z')\dif z'$.
\et

\br\rm

(i) If $b$ does not depend on the density variable $r$, then we can drop the regularity assumption on the initial distribution $\nu$ (see Remark \ref{Re71} below).
In particular, suppose that $a=a(x)$ satisfies \eqref{ELL} and $b=b(z)$ is bounded measurable, 
then for each starting point $z\in\mR^{2d}$, the following second order SDE admits a unique generalized martingale solution 
$\mP_z\in \widetilde\cM^{a,b}_{0,\delta_{z}}$ in the sense of Definition \ref{MP81}:
\begin{align}\label{SDE90}
\dif \dot X_t=b(Z_t)\dif t+\sqrt{a}(X_t)\dif W_t,\ \ Z_0=z. 
\end{align}
Moreover, let $b_\eps$ and $a_\eps$ be the smooth approximation of $b$ and $a$. Consider the following SDE:
$$
\dif \dot X^\eps_t=b_\eps(Z^\eps_t)\dif t+\sqrt{a_\eps}(X^\eps_t)\dif W_t,\ \ Z^\eps_0=z.
$$
By the uniqueness, for each $z\in\mR^{2d}$, the law of $Z^\eps$ weakly converges to $\mP_z$ (not necessarily subtracting a subsequence). 
In particular, the family of probability measures $(\mP_z)_{z\in\mR^{2d}}$ 
on the space of continuous functions forms a family of strong Markov processes (see \cite[p.184, Theorem 4.4.2]{EK86}). In fact, from the proof below, one can see that for each $z\in\mR^{2d}$,
$\mP_z$ is also a classical martingale solution in the sense of  Definition \ref{Def3}. However, 
it is not known whether {\it any} classical martingale solution priorily is a generalized martingale solution.
This seems to be an open problem.

(ii) Consider the following first order SDE in $\mR^d$,
$$
\dif X_t=b(X_t)\dif t+\sqrt{a}(X_t)\dif W_t,\ \ X_0=x.
$$ 
When $b$ is bounded measurable and $a$ satisfies the ellipticity condition \eqref{ELL}, 
it is well known that there is at least one weak solution to the above SDE (see \cite{Kry80} or \cite[Ex. 7.3.2]{St-Va}). Moreover, for $d=1, 2$, the weak uniqueness also holds (see \cite[Ex. 7.3.3, 7.3.4]{St-Va}).
However, for $d\geq 3$, to get the weak uniqueness, one usually needs to assume $a$ being continuous (see \cite[Theorem 7.2.1]{St-Va}). 
For the strong uniqueness, the well-known best condition for $a$ seems to be $a\in \mW^{2,p}$, the second order Sobolev space, where $p>d$ (see \cite{Z11}).
It is interesting that for second order SDEs with discontinuous diffusion coefficients, we have the well-posedness of generalized martingale problems.
\er

The crux point in our proof is to establish a maximum principle for the following kinetic FPKE:
$$
\p_t u=\tr(a\cdot\nabla^2_v u)+v\cdot\nabla_x u+b_1\cdot \nabla_v u+\div_v(b_2 u)+f,\ \ u(0)=\varphi,
$$
where an important observation is that since the diffusion-matrix $a$ does not depend on the velocity variable $v$, the principal second order term can be written in divergence form:
\begin{align}\label{DZ8}
\tr(a\cdot\nabla^2_v u)=\div_v(a\cdot\nabla_v u).
\end{align}
Thus one can invoke the classical De-Giorgi method to show the apriori $L^\infty$-estimates. 
Compared with the recent works \cite{GIMV19, GI21, GM21}, the novelty of our boundedness estimates is that $b$ can be in some $\mL^{q_1}_t(\mL^{\bbp_1}_z)$-space, where
$(q_1,\bbp_1)\in(2,\infty)^{1+2d}$ satisfies $\bba\cdot\frac1{\bbp_1}+\frac2{q_1}<1$,
and $f$ is in some  $\mL^{q_0}_t(\mL^{\bbp_0}_z)$-space, where
$(q_0,\bbp_0)\in(1,\infty)^{1+2d}$ satisfies $\bba\cdot\frac1{\bbp_0}+\frac2{q_0}<2$. 
More importantly, $f$  is also allowed to be a distribution in some Besov spaces with negative differentiability index, which plays a crucial role to show 
the existence of weak solutions in Theorem \ref{Th55}. That is, the distribution valued $f$ provides us the regularity of the density so that one can show the strong convergence of the densities of approximating equations (see Lemma \ref{Le57} below).

\subsection{Plan of the paper}
This paper is organized as follows. 
In Section 2, we prepare some preliminary results for later use. More concretely,
in Subsection 2.1, we introduce necessary anisotropic Besov spaces and some basic properties related to the Besov spaces.
It is noted that our $L^p$-spaces have different integrability indices along different components as explained in Remark \ref{Re14} above. 
In Subsection 2.2 we provide an abstract criterion for the local bound of
a function in De-Giorgi's class. Compared with the classical notion of De-Giorgi's class \cite{LSU68}, 
our definition does not depend on any structure of PDEs and is only in the scope of $L^p$-spaces. 

\medskip

In Section 3, we recall some well-known properties about weak (sub)-solutions.
In Subsection 3.2, we show how to improve the regularity in Besov spaces via Duhamel's formula for kinetic equations, which is essentially proved in \cite{HWZ20, ZZ21}.
Moreover, we also show the existence and uniqueness of weak solutions for the Cauchy problem of a linear kinetic FPKE (see Theorem \ref{Th3}), which seems to be new.

\medskip

Section 4 is devoted to showing the local boundedness of weak solutions via De-Giorgi's method. Two cases are considered. In the first case that the nonhomogeneous $f$ is an $L^p$-function,
we use the comparison method as in \cite{GIMV19} to improve the regularity of weak sub-solutions. In this case, we have best integrability conditions on drift $b$ and 
inhomogeneous $f$ as stated above. In the second case that the nonhomogeneous $f$ is allowed to be a distribution, we use the simple fact that for a weak solution $u$,
$(u^+)^2$ is still a weak solution  with $2f u^+-2\<a\cdot\nabla_vu^+,\nabla_vu^+\>$ in place of $f$ (see Lemma \ref{Le213} below). In this case, the $L^\infty$-estimate of weak solutions
in Theorem \ref{Th34} shall provide extra Besov regularity for the density of second order SDEs so that we can treat the density dependent SDEs (see Corollary \ref{Cor1}).

\medskip

In Section 5,  we show the well-posedness and stability of weak solutions for PDE \eqref{PDE0}. 
In particular, we obtain the global bound estimate and stability of weak solutions, where the key point is to use
the localization norm $\nor\cdot\nor_{\widetilde\mL^q_t(\widetilde\mL^\bbp_z)}$ introduced in \cite{ZZ21}. For the stability, we use the H\"older regularity estimate established in
\cite{GIMV19} (see also \cite{WZ09, GI21, GM21}).

\medskip

In Section 6, we prove the existence of weak solutions. As usual, in terms of the equivalence between weak solutions and martingale problems, 
we show the existence of classical martingale solutions associated with DDSDE \eqref{SDE10} via mollifying the coefficients.  Since our diffusion coefficient 
does not depend on velocity variable, the Kolmogorov equation associated SDEs \eqref{SDE10} can be written in the divergence form as explained in \eqref{DZ8}.
Thus one can utilize the maximum principle obtained in Section 5 to show the crucial Krylov estimate for approximating equation.
Then by Aubin-Lions' lemma, one can find a strong convergence subsequence of the densities.

\medskip

In Section 7, we show the well-posedness of generalized martingale problems. As said above,
such a notion strongly depends on the solvability of the associated Kolmogorov equations. When the solution of the associated Kolmogorov equations has enough regularity, saying $C^2$-smooth, 
the classical martingale solution must be a generalized martingale solution. In general, these two notions are not comparable. The advantage of using 
the notion of generalized martingale solutions is that 
the uniqueness is an easy consequence of the unique solvability of the associated Kolmogorov equation. The disadvantage is that it has few flexibility 
and It\^o's formula is not applicable so that it 
is not direct to see that the density of the solution solves the Fokker-Planck-Kolmogorov equation.

\medskip

We conclude this introduction by introducing the following convention: 
Throughout this paper, we use $C$ with or without subscripts to denote an unimportant constant, whose value
may change in different occasions. We also use $:=$ as a way of definition and $a^+:=0\vee a$. By $A\lesssim_C B$ and $A\asymp_C B$
or simply $A\lesssim B$ and $A\asymp B$, we mean that for some unimportant constant $C\geq 1$,
$$
A\leq C B,\ \ C^{-1} B\leq A\leq CB.
$$

\section{Preliminary}

\subsection{Anisotropic Besov spaces}
Fix $N\in\mN$. For multi-index $\bbp=(p_1,\cdots,p_N)\in(0,\infty]^N$, we define
\begin{align}\label{AM9}
\|f\|_{\mL^\bbp}:=\Bigg(\int_\mR\Bigg(\int_\mR\cdots\left(\int_\mR |f(z_1,\cdots,z_N)|^{p_N}\dif z_N\right)^{\frac{p_{N-1}}{p_N}}\cdots
\dif z_2\Bigg)^{\frac{p_1}{p_2}}\dif z_1\Bigg)^{\frac{1}{p_1}}.
\end{align}
When $\bbp=(p,\cdots,p)\in(0,\infty]^N$, we shall simply write 
$$
\mL^p=\mL^\bbp.
$$
Note that if $\bbp\in[1,\infty]^N$, then $\|\cdot\|_{\mL^\bbp_z}$ satisfies the triangle inequality and is a norm,
and for any permutation $\bbp'$ of $\bbp$, 
$$
\|f\|_{\mL^\bbp}\not=\|f\|_{\mL^{\bbp'}}.
$$
For multi-indices $\bbp, \bbq\in(0,\infty]^N$, we shall use the following notations:
$$
\frac1{\bbp}:=\Big(\frac1{p_1},\cdots,\frac1{p_N}\Big),\ \ \ \bbp\cdot\bbq:=\sum_{i=1}^Np_iq_i,
$$
and 
$$
\bbp>\bbq\ (\bbp\geq\bbq;\ \bbp=\bbq)\Leftrightarrow p_i>q_i\ (p_i\geq q_i;\ p_i=q_i),\ \ i=1,\cdots,N.
$$
Moreover,  we use bold number to denote constant vector in $\mR^N$:
$$
\boldsymbol{1}=(1,\cdots,1),\ \ \boldsymbol{2}=(2,\cdots,2).
$$
With a little confusion, for a set $Q\subset\mR^N$, we also use $\1_Q$ to denote the indicator function of $Q$.
For any multi-indices $\bbp,\bbq,\bbr\in(0,\infty]^N$ with $\tfrac1\bbp+\tfrac1\bbr=\tfrac1\bbq$, the following H\"older's inequality holds
$$
\|fg\|_{\mL^\bbq}\leq \|f\|_{\mL^\bbp}\|g\|_{\mL^\bbr}.
$$
For any multi-indices $\bbp,\bbq,\bbr\in[1,\infty]^N$ with $\tfrac1\bbp+\tfrac1\bbr=\boldsymbol{1}+\tfrac1\bbq$, the following Young's inequality holds
$$
\|f*g\|_{\mL^\bbq}\leq \|f\|_{\mL^\bbp}\|g\|_{\mL^\bbr}.
$$
These two inequalities are easy consequences of the well-known H\"older and Young's inequalities for $\bbp=(p,\cdots,p)$.

Let $\bba=(a_1,\cdots,a_N)\in[1,\infty)^N$ be a multi-scaling parameter. 
For $z,z'\in\mR^{N}$, we introduce the following anisotropic distance in $\mR^N$
$$
|z-z'|_\bba:=\sum_{i=1}^N|z_i- z'_i|^{1/a_i}.
$$
For $r>0$ and $z\in\mR^N$, we also introduce the ball with respect to the above distance
$$
B^\bba_r(z):=\{z'\in\mathbb{R}^N:|z'-z|_\bba\leq r\},\ \ B^\bba_r:=B^\bba_r(0).
$$
Let $\phi^\bba_0$ be  a symmetric $C^{\infty}$-function  on $\mathbb{R}^N$ with
$$
\phi^\bba_0(\xi)=1\ \mathrm{for}\ \xi\in B^\bba_1\ \mathrm{and}\ \phi^\bba_0(\xi)=0\ \mathrm{for}\ \xi\notin B^\bba_2.
$$
For $j\in\mathbb{N}$, we define
$$
\phi^\bba_j(\xi):=\phi^\bba_0(2^{-j\bba}\xi)-\phi^\bba_0(2^{-(j-1)\bba}\xi),
$$
where for $s\in\mR$ and $\xi=(\xi_1,\cdots,\xi_N)$,
$$
2^{s\bba }\xi=(2^{sa_1}\xi_1,\cdots, 2^{sa_N}\xi_N).
$$
For an $L^1$-integrable function $f$, let $\hat f$ be the Fourier transform of $f$ defined by
$$
\hat f(\xi):=(2\pi)^{-N/2}\int_{\mR^N} \e^{-{\rm i}\xi\cdot z}f(z)\dif z, \quad\xi\in\mR^N,
$$
and $\check f$ the Fourier  inverse transform of $f$ defined by
$$
\check f(z):=(2\pi)^{-N/2}\int_{\mR^N} \e^{{\rm i}\xi\cdot z}f(\xi)\dif\xi, \quad z\in\mR^N.
$$

Let $\cS$ be the space of all Schwartz functions on $\mR^N$ and $\cS'$ the dual space of $\cS$, called the tempered distribution space.
For given $j\in\mathbb{N}_0$, the block operator  $\mathcal{R}^\bba_j$ is defined on $\cS'$ by
\begin{align}\label{Ph0}
\mathcal{R}^\bba_jf(z):=(\phi^\bba_j\hat{f})\check{\ }(z)=\check{\phi}^\bba_j*f(z),
\end{align}
where the convolution is understood in the distributional sense and by scaling,
\begin{align}\label{SX4}
\check{\phi}^\bba_j(z)=2^{(j-1)|\bba|}\check{\phi}^\bba_1(2^{(j-1)\bba}z).
\end{align}
Here and after, we denote
$$
|\bba|:=a_1+\cdots+a_N=\bba\cdot \1.
$$
For $j\geq 0$, by definition it is easy to see that
\begin{align}\label{KJ2}
\cR^\bba_j=\cR^\bba_j\widetilde\cR^\bba_j,\ \mbox{ where }\ \widetilde\cR^\bba_j:=\cR^\bba_{j-1}+\cR^\bba_{j}+\cR^\bba_{j+1}
\mbox{ with } \cR^\bba_{-1}\equiv 0,
\end{align}
and by the symmetry of $\phi^\bba_j$,
$$
\<\cR^\bba_j f,g\>=\< f,\cR^\bba_jg\>,\ \ f\in\cS', \ g\in\cS.
$$
The cut-off low frequency operator $S_k$ is defined by
\begin{align}\label{EM9}
S_kf:=\sum_{j=0}^{k-1}\cR^\bba_j f\to f,\ \ k\to\infty.
\end{align}
For $f,g\in\cS'$, define
\begin{align}\label{EM39}
f\prec g:=\sum_{k\geq 0} S_{k-1}f\cR^\bba_k g,\ \ f\circ g:=\sum_{|i-j|\leq1}\cR^\bba_i f\cR^\bba_jg.
\end{align}
The Bony decomposition of $fg$ is formally given by (cf. \cite{BCD11, Bon81})
\begin{align}\label{Bony}
fg=f\prec g+ f\circ g+g\prec f.
\end{align}

Now we introduce the following anisotropic Besov spaces (cf. \cite[Chapter 5]{Tri06}).
\begin{definition}(Anisotropic Besov space)\label{bs}
For $s\in\mR$ and $\bbp\in[1,\infty]^N$, we define
$$
\mathbf{B}^s_{\bbp;\bba}:=\left\{f\in \cS': \|f\|_{\mathbf{B}^{s}_{\bbp;\bba}}
:=\sup_{j\geq0}\big(2^{s j}\|\mathcal{R}^\bba_{j}f\|_{\mL^\bbp}\big)<\infty\right\}.
$$
\end{definition}
\br\rm
For $s\in(0,1)$ and $\bbp\in[1,\infty]^N$, it is well known that $\mathbf{B}^s_{\bbp;\bba}$ has the following characterization (see \cite{Tri06} or \cite[Theorem 2.7]{HZZZ21}):
\begin{align}\label{EC1}
\|f\|_{\mathbf{B}^{s}_{\bbp;\bba}}\asymp\|f\|_{\mL^{\bbp}}+\sup_h\|f(\cdot+h)-f(\cdot)\|_{\mL^{\bbp}}/|h|^{s}_\bba.
\end{align}
Although the proof in \cite{HZZZ21} is for $\bbp=(p,\cdots,p)$, it also works for general $\bbp\in[1,\infty]^N$.
\er
The following inequality of Bernstein's type is quite useful (cf. \cite{BCD11}).
\bl\label{Bern}
Let $\bbp,\bbq\in[1,\infty]^N$ with $\bbp\leq\bbq$.
For any $k\in\mN_0$ and $i=1,\cdots,N$, 
there is a constant $C=C(\bbp,\bbq,\bba,k,i)>0$ such that for all $j\geq 0$,
\begin{align}\label{Ber}
\|\p_{z_i}^k\cR_j^\bba f\|_{\mL^\bbq}\lesssim_C2^{j(a_ik+\bba\cdot(\frac{1}{\bbp}-\frac{1}{\bbq}))}\|\cR_j^\bba f\|_{\mL^\bbp},
\end{align}
where $\p_{z_i}^k$ denotes the $k$-order partial derivative with respect to $z_i$. Moreover,
\begin{align}\label{Crapp}
\|\cR_j^\bba f\|_{\mL^\bbp}\lesssim_C\|f\|_{\mL^\bbp},\ j\geq 0.
\end{align}
\el
\begin{proof}
We only prove \eqref{Ber} for $j\geq 1$. Let $\bbr\in[1,\infty]^N$ be defined by
$$
\tfrac1\bbp+\tfrac1\bbr=\1+\tfrac1\bbq.
$$
By \eqref{KJ2}, \eqref{Ph0} and Young's inequality, we have
\begin{align*}
\|\p_{z_i}^k\cR_j^\bba f\|_{\mL^\bbq}=\|\p_{z_i}^k\widetilde\cR_j^\bba\cR_j^\bba f\|_{\mL^\bbq}
=\|(\p_{z_i}^k\check{\widetilde{\phi_{j}^\bba}})*\cR_j^\bba f\|_{\mL^\bbq}
\lesssim \|\p_{z_i}^k\check{\widetilde{\phi_{j}^\bba}}\|_{\mL^\bbr}\|\cR_j^\bba f\|_{\mL^\bbp},
\end{align*}
where 
$$
\widetilde{\phi_{j}^\bba}:=\phi_{j-1}^\bba+\phi_{j}^\bba+\phi_{j+1}^\bba.
$$
By \eqref{SX4} and the change of variables, we have
\begin{align*}
\|\p_{z_i}^k\check{\widetilde{\phi_{j}^\bba}}\|_{\mL^\bbr}= 2^{(j-1)a_ik}
2^{(j-1)\bba\cdot(\1-\frac{1}{\bbr})}\|\p_{z_i}^k\check{\widetilde{\phi_1^\bba}}\|_{\mL^\bbr}.
\end{align*}
Thus we get \eqref{Ber}. 
\end{proof}

We have the following interpolation inequality of  Gagliado-Nirenberge's type.
\bl
Let $\bbp,\bbq,\bbr\in[1,\infty]^N$ and $s,s_0,s_1\in\mR$, $\theta\in[0,1]$. Suppose that
\begin{align}\label{DD01}
\tfrac1{\bbp}\leq\tfrac{1-\theta}{\bbq}+\tfrac{\theta}{\bbr},\ \ s-\bba\cdot\tfrac{1}{\bbp}=(1-\theta)\big(s_0-\bba\cdot\tfrac{1}{\bbq}\big)+\theta\big(s_1-\bba\cdot\tfrac{1}{\bbr}\big).
\end{align}
Then there is a constant $C=C(\bbp,\bbq,\bbr,s,s_0,s_1,\theta)>0$ such that
\begin{align}\label{Sob}
\|f\|_{\bB^s_{\bbp;\bba}}\lesssim_C\|f\|^{1-\theta}_{\bB^{s_0}_{\bbq;\bba}}\|f\|^\theta_{\bB^{s_1}_{\bbr;\bba}}.
\end{align}
\el
\begin{proof}
For $\theta=0,1$, \eqref{Sob} is direct by \eqref{Ber}. Below we assume $\theta\in(0,1)$.
Let $\bbr\leq\bbw\in[1,\infty]^N$ and $s_2\in\mR$ be defined by 
\begin{align}\label{DD02}
\tfrac{1}{\bbp}=\tfrac{1-\theta}{\bbq}+\tfrac{\theta}{\bbw},\ \ s_0(1-\theta)+s_2\theta=s.
\end{align}
By H\"older's inequality, we have
\begin{align*}
\|\cR^\bba_j f\|_{\mL^\bbp}&=\|(\cR^\bba_j f)^{1-\theta}(\cR^\bba_j f)^{\theta}\|_{\mL^\bbp}
\leq\|\cR^\bba_j f\|^{1-\theta}_{\mL^{\bbq}}\|\cR^\bba_j f\|^\theta_{\mL^{\bbw}}.
\end{align*}
Hence,
\begin{align}
2^{js}\|\cR^\bba_j f\|_{\mL^\bbp}&\leq\Big(2^{js_0}\|\cR^\bba_j f\|_{\mL^{\bbq}}\Big)^{1-\theta}
\Big(2^{js_2}\|\cR^\bba_j f\|_{\mL^{\bbw}}\Big)^\theta.\label{DD03}
\end{align}
Note that by \eqref{DD01} and \eqref{DD02},
$$
s_1=s_2+\bba\cdot(\tfrac 1{\bbr}-\tfrac 1{\bbw}),
$$ 
and by \eqref{Ber},
$$
\|\cR^\bba_j  f\|_{\mL^{\bbw}}\lesssim 2^{j\bba\cdot(\frac 1{\bbr}-\frac 1{\bbw})}\|\cR^\bba_j  f\|_{\mL^{\bbr}},\ j\geq 0.
$$
Substituting it into \eqref{DD03}, we obtain
$$
2^{js}\|\cR^\bba_j f\|_{\mL^\bbp}\lesssim\Big(2^{js_0}\|\cR^\bba_j f\|_{\mL^{\bbq}}\Big)^{1-\theta}
\Big(2^{js_1}\|\cR^\bba_j f\|_{\mL^{\bbq}}\Big)^\theta
\leq\|f\|^{1-\theta}_{\bB^{s_0}_{\bbq;\bba}}\|f\|^{\theta}_{\bB^{s_1}_{\bbq;\bba}}.
$$
The proof is complete.
\end{proof}
We also need the following two simple results. Since they are crucial for treating distribution-valued inhomogeneous $f$, 
we provide detailed proofs here for readers' convenience.
\bl\label{Le23}
For any $s'>s>0$ and $\bbp,\bbq\in[1,\infty]^N$ with $\frac{1}{\bbp}+\frac1{\bbq}=\1$, there is a constant $C=C(\bbp,\bbq,s,s',N)>0$
such that for all $f\in\bB^{-s}_{\bbp;\bba}$ and $g\in\bB^{s'}_{\bbq;\bba}$,
$$
|\<\!\!\!\<f,g\>\!\!\!\>|+\|fg\|_{\bB^{-s}_{1;\bba}}\lesssim_C\|f\|_{\bB^{-s}_{\bbp;\bba}}\|g\|_{\bB^{s'}_{\bbq;\bba}},
$$
where $\<\!\!\!\<\cdot,\cdot\>\!\!\!\>$ stands for the dual pair between $\cS'$ and $\cS$.
\el
\begin{proof}
By \eqref{EM9}, \eqref{KJ2} and H\"older's inequality, we have
\begin{align*}
\<\!\!\!\<f,g\>\!\!\!\>=\sum_{j\geq0}\<\!\!\!\<\cR^\bba_jf,\widetilde\cR^\bba_jg\>\!\!\!\>
\leq\sum_{j\geq 0}\|\cR^\bba_jf\|_{\mL^\bbp}\|\widetilde\cR^\bba_jg\|_{\mL^\bbq}
\leq\sum_{j\geq0} 2^{(s-s')j}\|f\|_{\bB^{-s}_{\bbp;\bba}}\|g\|_{\bB^{s'}_{\bbq;\bba}},
\end{align*}
where the last series converges due to $s'>s$. Next we show
$$
\|fg\|_{\bB^{-s}_{1;\bba}}\lesssim_C\|f\|_{\bB^{-s}_{\bbp;\bba}}\|g\|_{\bB^{s'}_{\bbq;\bba}}.
$$
By Bony's decomposition \eqref{Bony}, it suffices to show the following three estimates:
\begin{align}
\|f\prec g\|_{\bB^{s'-s}_{1;\bba}}&\lesssim_C\|f\|_{\bB^{-s}_{\bbp;\bba}}\|g\|_{\bB^{s'}_{\bbq;\bba}},\label{DD1}\\
\|g\prec f\|_{\bB^{-s}_{1;\bba}}&\lesssim_C\|g\|_{\mL^\bbq}\|f\|_{\bB^{-s}_{\bbp;\bba}},\label{DD2}\\
\|f\circ g\|_{\bB^{0}_{1;\bba}}&\lesssim_C\|f\|_{\bB^{-s}_{\bbp;\bba}}\|g\|_{\bB^{s'}_{\bbq;\bba}}.\label{DD3}
\end{align}
Noting that  (see \cite[(2.9)]{HZZZ21}),
$$
\cR^\bba_j(S_{k-1}f\cR^\bba_k g)=0,\ \ |k-j|>2,
$$
by definition \eqref{EM39}, \eqref{Crapp} and H\"older's inequality, we have
\begin{align*}
\|\cR^\bba_j(f\prec g)\|_{\mL^1}&\lesssim \sum_{|k-j|\leq 2} \|S_{k-1}f\cR^\bba_k g\|_{\mL^1}
\leq \sum_{|k-j|\leq 2} \|S_{k-1}f\|_{\mL^\bbp}\|\cR^\bba_k g\|_{\mL^\bbq}\\
&\lesssim \sum_{|k-j|\leq 2} \sum_{i=0}^{k-2}2^{s i}2^{-k s'}\|f\|_{\bB^{-s}_{\bbp;\bba}} \|g\|_{\bB^{s'}_{\bbq;\bba}}\lesssim 
2^{j(s-s')}\|f\|_{\bB^{-s}_{\bbp;\bba}} \|g\|_{\bB^{s'}_{\bbq;\bba}},
\end{align*}
which gives \eqref{DD1} by definition. Similarly,
\begin{align*}
\|\cR^\bba_j(g\prec f)\|_{\mL^1}&\leq \sum_{|k-j|\leq 2} \|S_{k-1}g\cR^\bba_k f\|_{\mL^1}
\lesssim \sum_{|k-j|\leq 2} \|S_{k-1}g\|_{\mL^\bbq}\|\cR^\bba_k f\|_{\mL^\bbp}\\
&\lesssim \sum_{|k-j|\leq 2}\Big(\|g\|_{\mL^\bbq}2^{sk}\|f\|_{\bB^{-s}_{\bbp;\bba}}\Big)\lesssim 2^{sj}\|g\|_{\mL^\bbq}\|f\|_{\bB^{-s}_{\bbp;\bba}},
\end{align*}
which yields \eqref{DD2}.
Moreover,  we also have
\begin{align*}
\|\cR^\bba_j(f\circ g)\|_{\mL^1}&\lesssim\|f\circ g\|_{\mL^1}\leq \sum_{|i-k|\leq 2} \|\cR^\bba_if\cR^\bba_k g\|_{\mL^1}
\lesssim \sum_{|i-k|\leq 2} \|\cR^\bba_if\|_{\mL^\bbp}\|\cR^\bba_k g\|_{\mL^\bbq}\\
&\lesssim \sum_{|i-k|\leq 2} 2^{si}2^{-ks'}\|f\|_{\bB^{-s}_{\bbp;\bba}}\|g\|_{\bB^{s'}_{\bbq;\bba}}\lesssim \|f\|_{\bB^{-s}_{\bbp;\bba}}\|g\|_{\bB^{s'}_{\bbq;\bba}},
\end{align*}
which gives \eqref{DD3}.
The proof is complete.
\end{proof}
\bl\label{Le24}
For any $\bbp\in[1,\infty]^N$ and $s\in(0,1)$, there is a constant $C=C(\bbp,\bba,s)>0$ such that for any nonnegative $u$,
$$
\|u\|^2_{\bB^{s/2}_{2\bbp;\bba}}\lesssim\|u^2\|_{\bB^{s}_{\bbp;\bba}}.
$$
\el
\begin{proof}
Since $|a-b|\leq\sqrt{|a^2-b^2|}$ for $a,b\geq 0$, by \eqref{EC1}, we have
\begin{align*}
\|u\|_{\bB^{s/2}_{2\bbp;\bba}}&\asymp\|u\|_{\mL^{2\bbp}}+\sup_h\|u(\cdot+h)-u(\cdot)\|_{\mL^{2\bbp}}/|h|^{s/2}_\bba\\
&\leq \|u\|_{\mL^{2\bbp}}+\sup_h\|\sqrt{|u^2(\cdot+h)-u^2(\cdot)|}\|_{\mL^{2\bbp}}/|h|^{s/2}_\bba\\
&=\|u^2\|^{1/2}_{\mL^{\bbp}}+\sup_h\|u^2(\cdot+h)-u^2(\cdot)\|^{1/2}_{\mL^{\bbp}}/|h|^{s/2}_\bba\asymp \|u^2\|^{1/2}_{\bB^{s}_{\bbp;\bba}}.
\end{align*}
The proof is complete.
\end{proof}

\subsection{De-Giorgi's class}
In this subsection we present a general criterion for the local bound of a function in De-Giorgi's class.
Let $\sI\subset(1,\infty)^{N}$ be an open  multi-index set and $Q:=(Q_\tau)_{\tau\in[1,2]}$ 
be a family of increasing bounded open sets in $\mR^{N}$ with
\begin{align}\label{DD81}
Q_\tau\cap Q^c_{\sigma}=\emptyset\ \mbox{ for } \tau<\sigma\ \mbox{ and } \cap_{\sigma>\tau}Q_\sigma=\bar Q_{\tau}.
\end{align}

We introduce the following De-Giorgi class associated with $\sI$ and $Q$. 
\bd\label{Def6}
We call a function $u\in L^1(Q_2)$ being in the De-Giorgi class $\cD\cG^+_\sI(Q)$ 
if there are $\bbp_i\in\sI$, $i=1,\cdots,m$, $1\leq j<m$ and $\lambda, \cA\geq 0$ such that for any $\bbp\in\sI$, there is a constant $C_\bbp>0$ such that
for any $1\leq\tau<\sigma\leq 2$ and $\kappa\geq 0$,
\begin{align}\label{DX-1}
(\sigma-\tau)^{\lambda}\|\1_{Q_{\tau}}(u-\kappa)^+\|_{\mL^\bbp}\lesssim_{C_\bbp}
\sum_{i=1}^{j}\|\1_{Q_{\sigma}}(u-\kappa)^+\|_{\mL^{\bbp_i}}+\cA\sum_{i=j+1}^{m}\|\1_{\{u>\kappa\}\cap Q_{\sigma}}\|_{\mL^{\bbp_i}}.
\end{align}
\ed
The aim of this subsection is to give an upper bound estimate of $u\in\cD\cG^+_\sI(Q)$.
First of all, we recall the following two iteration lemmas (see \cite{HL97}).
\bl\label{Le15}
Let $(a_n)_{n\in\mN}$ be a sequence of nonnegative real numbers. Suppose that for some $C_0,\lambda>1$ and $\delta>0$,
$$
a_{n+1}\leq C_0\lambda^n a_n^{1+\delta}.
$$
Then under $a_1\leq (C_0\lambda^{(1+\delta)/\delta})^{-1/\delta}$, we have
$$
\lim_{n\to\infty} a_n=0.
$$
\el
\bl\label{Le26}
Let $h(\tau)\geq 0$ be bounded in $[\tau_1,\tau_2]$ with $\tau_1\geq 0$. Suppose that for some $\alpha\geq 0$, $\theta\in(0,1)$ and any
$\tau_1\leq\tau<\tau'\leq\tau_2$, 
$$
h(\tau)\leq\theta h(\tau')+(\tau'-\tau)^{-\alpha}A+B.
$$
Then there is a constant $C=C(\alpha,\theta)>0$ such that
$$
h(\tau_1)\lesssim_C(\tau_2-\tau_1)^{-\alpha}A+B.
$$
\el
Now we can show the following main result of this subsection (cf. \cite{ZZ21}).
\bt\label{TH20}
Let $u\in\cD\cG^+_\sI(Q)$. For any $p>0$, there are constants $\gamma, C>0$ only 
depending on $p$, ${\rm Vol}(Q_2)$ and the parameters $\bbp_i,\lambda, C_{\bbp_i}$ in the definition of $\cD\cG^+_\sI(Q)$ such that 
for all $1\leq\tau<\sigma\leq 2$,
\begin{align}\label{JK2}
\|u^+\1_{Q_\tau}\|_{\mL^\infty}\lesssim_C(\sigma-\tau)^{-\gamma}\|u^+\1_{Q_\sigma}\|_{\mL^p}+\cA.
\end{align}
\et
\begin{proof}
Fix $1\leq\tau<\sigma\leq 2$ and $\kappa>0$, which will be determined below. For $n\in\mN$, define
$$
\tau_n=\tau+(\sigma-\tau)2^{1-n},\ \ \kappa_n:=\kappa\left(1-2^{1-n}\right),
$$
and
$$
\ w_n:=(u-\kappa_n)^+.
$$
Clearly, by \eqref{DD81},
$$
\kappa_n\uparrow\kappa,\ \tau_n\downarrow\tau,\ \ Q_{\tau_{n+1}}\subset Q_{\tau_n}\downarrow \bar Q_\tau.
$$
Since $\sI$ is an open set, for any $\bbp\in\sI$, there is a $\bbp'\in\sI$ such that $\bbp'>\bbp$.
Let $\bbq=(q_{1},\cdots,q_{N})\in[1,\infty)^N$ be defined by
$$
\tfrac1{\bbp'}+\tfrac1{\bbq}=\tfrac1{\bbp}.
$$
By H\"older's inequality, there is a constant $C=C(Q_2,\bbq)>0$ such that
\begin{align}
\|\1_{Q_{\tau_{n+1}}}w_{n+1}\|_{\mL^{\bbp}}
&\leq\|\1_{Q_{\tau_{n+1}}}w_{n+1}\|_{\mL^{\bbp'}}\|\1_{\{w_{n+1}\not=0\}\cap Q_{\tau_{n+1}}}\|_{\mL^{\bbq}}\no\\
&\leq\|\1_{Q_{\tau_{n+1}}}w_{n+1}\|_{\mL^{\bbp'}}\|\1_{\{w_{n+1}\not=0\}\cap Q_{\tau_{n}}}\|_{\mL^{\bbq}}\no\\
&\lesssim_C\|\1_{Q_{\tau_{n+1}}}w_{n+1}\|_{\mL^{\bbp'}}\|\1_{\{w_{n+1}\not=0\}\cap Q_{\tau_{n}}}\|_{\mL^{1}}^{1/\max_i q_i},\label{EC2}
\end{align}
where the last step is due to $\1_A^p=\1_A$ for any $p>0$.
Thus, for $\bbp_i\in\sI$ being as in Definition \ref{Def6}, if we let
$$
a_n:=\sum_{i=1}^{m}\|\1_{Q_{\tau_{n}}}w_{n}\|_{\mL^{\bbp_i}},
$$
then by \eqref{EC2}, there are $\bbp_i'\in\sI$ and $C_1,\delta>0$ such that
\begin{align}\label{EB1}
a_{n+1}\lesssim_{C_1}\sum_{i=1}^m\|\1_{Q_{\tau_{n+1}}}w_{n+1}\|_{\mL^{\bbp'_i}}\|\1_{\{w_{n+1}\not=0\}\cap Q_{\tau_{n}}}\|_{\mL^{1}}^\delta.
\end{align}
On the other hand, noting that
$$
w_n|_{w_{n+1}\not=0}=(u-\kappa_{n+1}+\kappa_{n+1}-\kappa_n)^+|_{w_{n+1}\not=0}\geq \kappa_{n+1}-\kappa_n=\kappa 2^{-n},
$$
for any $\bbp\in[1,\infty)^N$, we have
\begin{align}\label{EB2}
\|\1_{Q_{\tau_n}}w_n\|_{\mL^\bbp}
\geq\|\1_{\{{w_{n+1}\not=0}\}\cap Q_{\tau_n}}w_n\|_{\mL^\bbp}\geq \kappa 2^{-n}\|\1_{\{w_{n+1}\not=0\}\cap Q_{\tau_n}}\|_{\mL^\bbp}.
\end{align}
Suppose
$$
\kappa\geq \cA.
$$
For each $i=1,\cdots,m$, since $\bbp_i'\in\sI$, by \eqref{DX-1}, \eqref{EB2} and $w_{n+1}\leq w_n$, we have
\begin{align*}
(\tau_{n}-\tau_{n+1})^\lambda\|\1_{Q_{\tau_{n+1}}}w_{n+1}\|_{\mL^{\bbp_i'}}
&\lesssim\sum_{i=1}^j\|\1_{Q_{\tau_{n}}}w_{n+1}\|_{\mL^{\bbp_i}}+\cA\sum_{i=j+1}^m\|\1_{\{w_{n+1}>0\}\cap Q_{\tau_n}}\|_{\mL^{\bbp_i}}\\
&\lesssim
\sum_{i=1}^j\|\1_{Q_{\tau_{n}}}w_{n+1}\|_{\mL^{\bbp_i}}+2^n\sum_{i=j+1}^m\|\1_{Q_{\tau_{n}}}w_{n}\|_{\mL^{\bbp_i}}\\
&\lesssim2^n\sum_{i=1}^m\|\1_{Q_{\tau_{n}}}w_{n}\|_{\mL^{\bbp_i}}=2^n a_n.
\end{align*}
Substituting this into \eqref{EB1} and by $\tau_n-\tau_{n+1}=(\sigma-\tau)2^{-n}$ and \eqref{EB2}, we obtain that for $\kappa\geq\cA$,
\begin{align*}
(\sigma-\tau)^\lambda 2^{-\lambda n}a_{n+1}&\lesssim2^na_n\|\1_{\{w_{n+1}\not=0\}\cap Q_{\tau_{n}}}\|_{\mL^{1}}^\delta\\
&\lesssim2^na_n\|\1_{\{w_{n+1}\not=0\}\cap Q_{\tau_{n}}}\|_{\mL^{\bbp_1}}^\delta\\
&\leq2^na_n(2^n\kappa^{-1})^\delta \|\1_{Q_{\tau_{n}}}w_{n}\|_{\mL^{\bbp_1}}^\delta\\
&\leq2^na_n(2^n\kappa^{-1})^\delta a_n^\delta.
\end{align*}
In particular, there is a $C_2>0$ such that for all $n\in\mN$,
$$
a_{n+1}\leq C_2(\sigma-\tau)^{-\lambda}\kappa^{-\delta} 2^{(\lambda+\delta+1)n} a_n^{1+\delta}.
$$
Now let us choose
$$
\kappa=\left(\left[C_2(\sigma-\tau)^{-\lambda}2^{(\lambda+\delta+1)(1+\delta)/\delta}\right]^{1/\delta}a_1\right)\vee\cA.
$$
Then by Fatou's lemma and Lemma \ref{Le15},
\begin{align*}
\|(u-\kappa)^+\1_{Q_\tau}\|_{\mL^{\bbp_1}}\leq \liminf_{n\to\infty}\|w_n\1_{Q_{\tau_n}}\|_{\mL^{\bbp_1}}
\leq\limsup_{n\to\infty}a_n=0,
\end{align*}
which in turn implies that $(u-\kappa)^+=0$ on $Q_{\tau}$ and by the choice of $\kappa$,
\begin{align}
\|u^+\1_{Q_\tau}\|_\infty&\leq\left(\left[C_2(\sigma-\tau)^{-\lambda}2^{(\lambda+\delta+1)(1+\delta)/\delta}\right]^{1/\delta}a_1\right)\vee\cA\no\\
&\leq C_3(\sigma-\tau)^{-\frac\lambda\delta}\sum_{i=0}^m\|u^+\1_{Q_{\sigma}}\|_{\mL^{\bbp_i}}+\cA.\label{DP9}
\end{align}
To show \eqref{JK2}, without loss of generality, we may assume 
$$
p\leq\gamma/2,\ \ \gamma:=\max_{i,j}p_{ij},\ \bbp_i=(p_{i1},\cdots, p_{iN}).
$$
By \eqref{DP9}, H\"older's inequality and Young's inequality, we have
\begin{align*}
\|u^+\1_{Q_\tau}\|_\infty&\leq C(\sigma-\tau)^{-\frac\lambda\delta}\|u^+\1_{Q_\sigma}\|_{\mL^\gamma}+\cA\\
&\leq C(\sigma-\tau)^{-\frac\lambda\delta}\|u^+\1_{Q_\sigma}\|^{1-\frac{p}{\gamma}}_\infty
\|u^+\1_{Q_\sigma}\|^{\frac{p}{\gamma}}_{\mL^p}+\cA\\
&\leq\tfrac12\|u^+\1_{Q_\sigma}\|_\infty+C(\sigma-\tau)^{-\frac{\gamma}{p\delta}}\|u^+\1_{Q_\sigma}\|_{\mL^p}+\cA,
\end{align*}
where $C$ is independent of $\sigma,\tau$.
Thus by Lemma \ref{Le26}, we obtain \eqref{JK2}.
\end{proof}

\section{Weak (sub)-solutions of kinetic equations}

In this section we present some basic properties about weak (sub)-solutions of linear kinetic FPKEs
and show the gain of regularity in $x$ via Duhamel's formula.
Consider the following kinetic equation of divergence form:
\begin{align}\label{PDE0}
\p_t u=\div_v(a\cdot\nabla_v u)+v\cdot\nabla_x u+b\cdot \nabla_v u+f,
\end{align}
where 
$$
a:\mR^{1+2d}\to \mM^d_{\rm sym},\ b:\mR^{1+2d}\to\mR^d,\ f:\mR^{1+2d}\to\mR
$$
are Borel measurable functions.
Suppose that for some $0<\kappa_0<\kappa_1$,
\begin{align}\label{C1}
\kappa_0\mI\leq a\leq\kappa_1\mI,
\end{align}
and for any bounded $Q\subset\mR^{1+2d}$,
\begin{align}\label{C2}
b\1_Q\in\mL^2,\ \ f\1_Q\in\mL^1,
\end{align}
where for $p,q\in[1,\infty]$, 
$$
\mL^q_t(\mL^p_z):=L^q(\mR; L^p(\mR^{2d})),\ \ \mL^p:=\mL^p_t(\mL^p_z).
$$
We introduce the following space of solutions: for an open set $Q\subset\mR^{1+2d}$,
$$
\sV_{Q}:=\Big\{ f\in\mL^1_{loc}: \|f\|_{\sV_Q}:=\|\1_Qf\|_{\mL^\infty_t(\mL^2_z)}+
\|\1_Q \nabla_v f\|_{\mL^2}<\infty\Big\}.
$$
For simplicity, we write for any $T>0$,
$$
\sV_T:=\sV_{[0,T]\times\mR^{2d}},\ \ \sV:=\sV_{\mR^{1+2d}}, \ \ \sV_{loc}:=\cap_{{\rm bounded }\, Q}\sV_Q,
$$
and for given $t\in\mR$ and $r>0$,
\begin{align}\label{DD069}
\cI_t:=\1_{(-\infty,t]},\ \ Q_r:=\big\{(t,x,v): |t|<r^2, |x|<r^3, |v|<r\big\}.
\end{align}
\subsection{Weak (sub)-solutions}
In this subsection we introduce the notion of weak (sub)-solutions and their basic properties for later use.
\bd\label{Def31}
Let $Q\subset\mR^{1+2d}$ be a bounded domain.
A function $u\in\sV_Q\cap\mL^\infty_Q$ is called a weak sub-solution of PDE \eqref{PDE0} in $Q$
if for any nonnegative $\varphi\in C^\infty_c(Q)$,
\begin{align}\label{Def0}
\begin{split}
-\int_{Q}u\p_t\varphi\leq &-\int_{Q}\<a\cdot\nabla_v u,\nabla_v\varphi\>
-\int_{Q}u(v\cdot\nabla_x \varphi)
+\int_{Q}(b\cdot\nabla_v u)\varphi+\int_{Q}f\varphi.
\end{split}
\end{align}
If both $u$ and $-u$ are weak sub-solutions, then we call $u$ a weak solution of PDE \eqref{PDE0} in $Q$.
If $u\in\sV_{loc}\cap\mL^\infty_{loc}$ and \eqref{Def0} holds for any bounded domain $Q$, then it is called a global weak solution. 
\ed
\br\label{Re32}\rm
By \eqref{C1} and \eqref{C2}, each term in \eqref{Def0} is well-defined.
Let $u$ be a global weak solution.
By a standard approximation, one sees that \eqref{Def0} is equivalent that for any $\varphi\in C^\infty_c(\mR^{2d})$ and Lebesgue almost all $t_0<t_1$,
\begin{align*}
\<\!\!\!\<u,\varphi\>\!\!\!\>|^{t_1}_{t_0}\leq &-\int^{t_1}_{t_0}
\<\!\!\!\<a\cdot\nabla_v u,\nabla_v\varphi\>\!\!\!\>
-\int^{t_1}_{t_0}\<\!\!\!\<v\cdot\nabla_x \varphi,u\>\!\!\!\>
+\int^{t_1}_{t_0}\<\!\!\!\<b\cdot\nabla_v u,\varphi\>\!\!\!\>+\int^{t_1}_{t_0}\<\!\!\!\<f,\varphi\>\!\!\!\>,
\end{align*}
where $\<\!\!\!\<u,\varphi\>\!\!\!\>=\int_{\mR^{2d}}u(z)\varphi(z)\dif z$.
In particular, let $u$ be a weak solution of the Cauchy problem of \eqref{PDE0} with initial value $u(0)=0$. Then 
it can be extended to be a global weak solution by
setting $u(t)=f(t)\equiv 0$ for $t\leq 0$.
\er

The following two lemmas are well known to experts. For readers' convenience, we provide detailed proofs here.
\bl\label{Le12}
Let $u\in\sV_{Q}\cap\mL^\infty_Q$ be a nonnegative weak sub-solution of PDE \eqref{PDE0} in  $Q$. Under \eqref{C1} and \eqref{C2},
for any nonnegative $\eta\in C^\infty_c(Q)$ and $t\in\mR$, it holds that
\begin{align}\label{ES5}
\begin{split}
\tfrac12\int_{\mR^{2d}}|(u\eta)(t)|^2&\leq\int_{\mR^{1+2d}}u^2\eta(\p_s\eta-v\cdot\nabla_x\eta)\cI_t-\int_{\mR^{1+2d}}
\<a\cdot\nabla_vu, \nabla_v(u\eta^2)\>\cI_t\\
&\quad+\int_{\mR^{1+2d}}(b\cdot\nabla_v u) u\eta^2\cI_t+\int_{\mR^{1+2d}}f u\eta^2\cI_t.
\end{split}
\end{align}
Moreover, if $u$ is a weak solution, then the above inequality becomes an equality.
\el
\begin{proof}
Let $\varGamma$ be a nonnegative symmetric smooth function in $\mR^{1+2d}$ with compact support.
For $\eps\in(0,1)$, we introduce the following mollifier:
\begin{align}\label{RR}
\varGamma_\eps(t,x,v):=\eps^{-(4d+2)}\varGamma\big(\tfrac{t}{\eps^2}, \tfrac{x}{\eps},\tfrac{v}{\eps^3}\big),\ \ u_\eps:=u*\varGamma_\eps.
\end{align}
Let $\phi\in C^\infty_c(Q)$ be nonnegative.
By the integration by parts, we have
$$
-\tfrac12\int_{\mR^{1+2d}} u^2_\eps\p_t\phi=\int_{\mR^{1+2d}}\p_tu_\eps u_\eps\phi=-\int_{\mR^{1+2d}} u_\eps\p_t(u_\eps\phi)
=-\int_{\mR^{1+2d}} u\p_t(\varGamma_\eps*(u_\eps\phi)).
$$
Thus by \eqref{Def0} with  $\varphi=\varGamma_\eps*(u_\eps\phi)\in C^\infty_c(Q)$ provided $\eps$ small enough, we get
\begin{align}
-\tfrac12\int_{\mR^{1+2d}} u^2_\eps\p_t\phi
&\leq-\int_{\mR^{1+2d}}\<a\cdot\nabla_vu, \varGamma_\eps*\nabla_v(u_\eps\phi)\>-\int_{\mR^{1+2d}}u(v\cdot\varGamma_\eps*\nabla_x(u_\eps\phi))\no\\
&\quad+\int_{\mR^{1+2d}}(\varGamma_\eps*(b\cdot\nabla_v u))u_\eps\phi+\int_{\mR^{1+2d}}(\varGamma_\eps*f) u_\eps\phi.\label{GG1}
\end{align}
Fix $t\in\mR$. Let $\chi_n$ be a family of smooth functions in $\mR$ so that
$$
\chi_n\to \cI_t=1_{(-\infty, t]},\ \ \chi'_n\to -\delta_{\{t\}}\ \ {\rm as }\ n\to\infty,
$$
where $\delta_{\{t\}}$ stands for the Dirac measure at point $t$.
Using $\chi_n\eta^2$ in place of $\phi$ in \eqref{GG1} and then taking limits $n\to\infty$, we obtain
\begin{align}
&\tfrac12\left(\int_{\mR^{2d}}|(u_\eps\eta)(t)|^2-\int_{\mR^{1+2d}}u^2_\eps\p_s\eta^2\cI_t\right)\no\\
&\quad\leq-\int_{\mR^{1+2d}}\<a\cdot\nabla_vu, \varGamma_\eps*\nabla_v(u_\eps\cI_t\eta^2)\>
-\int_{\mR^{1+2d}}u(v\cdot\varGamma_\eps*\nabla_x(u_\eps \cI_t\eta^2))\no\\
&\qquad+\int_{\mR^{1+2d}}(\varGamma_\eps*(b\cdot\nabla_v u))u_\eps\eta^2\cI_t+\int_{\mR^{1+2d}}(\varGamma_\eps*f) u_\eps\eta^2\cI_t.\label{GG2}
\end{align}
By the definition of convolutions, we have
\begin{align*}
&\left|\int_{\mR^{1+2d}}u(v\cdot\nabla_x\varGamma_\eps*(u_\eps\cI_t\eta^2))-v\cdot\nabla_x (u*\varGamma_\eps)(u_\eps\cI_t\eta^2)\right|\\
&=\left|\int_{\mR^{1+2d}}\left(\int_{\mR^{1+2d}}u(\bar s,\bar x,\bar v)(\bar v-v)\cdot\nabla_x\varGamma_\eps(s-\bar s, x-\bar x,v-\bar v)\right)(u_\eps\cI_t\eta^2)(s, x,v)\right|\\
&\lesssim\eps^2\int_{\mR^{1+2d}}\left(\int_{\mR^{1+2d}}|u(\bar s, \bar x,\bar v)|\,|\nabla_x\varGamma|
\big(\tfrac{s-\bar s}{\eps^2},\tfrac{x-\bar x}{\eps},\tfrac{v-\bar v}{\eps^3}\big)\right)|(u_\eps\cI_t\eta^2)(s,x,v)|,
\end{align*}
which converges to zero as $\eps\to 0$. Note that
\begin{align*}
\int_{\mR^{1+2d}}(v\cdot\nabla_x u_\eps)(u_\eps\cI_t\eta^2)=\tfrac12\int_{\mR^{1+2d}}(v\cdot\nabla_x u_\eps^2)\eta^2\cI_t
=-\tfrac12\int_{\mR^{1+2d}}(v\cdot\nabla_x \eta^2)u_\eps^2\cI_t.
\end{align*}
Substituting these into \eqref{GG2} and letting $\eps\to 0$, we get
\begin{align*}
&\tfrac12\int_{\mR^{2d}}|(u\eta)(t)|^2
\leq\int_{\mR^{1+2d}}u^2\eta\p_s\eta\cI_t
-\int_{\mR^{1+2d}}\<a\cdot\nabla_vu, \nabla_v(u\eta^2)\>\cI_t\\
&\qquad-\tfrac12\int_{\mR^{1+2d}}(v\cdot\nabla_x\eta^2)u^2\cI_t
+\int_{\mR^{1+2d}}(b\cdot\nabla_v u)u\eta^2\cI_t+\int_{\mR^{1+2d}} fu\eta^2\cI_t,
\end{align*}
where $u\in\mL^\infty_Q$ is used for taking limits for
$$
\lim_{\eps\to 0}\int_{\mR^{1+2d}}(\varGamma_\eps*(b\cdot\nabla_v u))u_\eps\eta^2\cI_t=\int_{\mR^{1+2d}}(b\cdot\nabla_v u)u\eta^2\cI_t
$$
and
$$
\lim_{\eps\to 0}\int_{\mR^{1+2d}}(\varGamma_\eps*f) u_\eps\eta^2\cI_t=\int_{\mR^{1+2d}} fu\eta^2\cI_t.
$$
Thus we obtain \eqref{ES5}. If $u$ is a weak solution, then from the above proof, one sees that \eqref{ES5} takes equality.
\end{proof}
\br\rm\label{Re12}
If $b=0$ and $\1_Q f\in\mL^1_t(\mL^2_z)$, it suffices to require $u\in\sV_{Q}$ in Lemma \ref{Le12}.
\er
The second part of the following lemma shall be used to deal with the distribution-valued inhomogeneous $f$.
\bl\label{Le213}
\begin{enumerate}[(i)]
\item Suppose that $u\in\sV_{loc}\cap \mL^\infty_{loc}$ is a weak sub-solution of \eqref{PDE0}. Then $u^+$ is still a weak sub-solution of \eqref{PDE0} with 
$f\1_{\{u>0\}}$ in place of $f$.

\item
Suppose that $u\in\sV_{loc}\cap \mL^\infty_{loc}$ is a  weak solution of \eqref{PDE0}. Then $(u^+)^2$ is still a weak solution of \eqref{PDE0} with 
$2f u^+-2\<a\cdot\nabla_vu^+,\nabla_vu^+\>$ in place of $f$.
\end{enumerate}
\el
\begin{proof}
Let $\varGamma_\eps$ be as in \eqref{RR}. By taking $\varphi=\varGamma_\eps(t-\cdot,x-\cdot,v-\cdot)$ in \eqref{Def0}, we obtain
$$
\p_t u_\eps\leq \div_v(a\cdot\nabla_vu)*\varGamma_\eps+(v\cdot\nabla_x u+b\cdot\nabla_v u)*\varGamma_\eps+f_\eps.
$$
First of all, we assume that $\beta:\mR\to\mR$ is a smooth non-decreasing function and $u$ is a weak sub-solution.
By the chain rule and $\beta'\geq 0$, we have
\begin{align*}
\p_t \beta(u_\eps)&\leq (\div_v(a\cdot\nabla_vu)*\varGamma_\eps+(v\cdot\nabla_x u+b\cdot\nabla_v u)*\varGamma_\eps+f_\eps)\beta'(u_\eps)\\
&=(\div_v(a\cdot\nabla_vu_\eps)+v\cdot\nabla_x u_\eps+b\cdot\nabla_v u_\eps+H_\eps+f_\eps)\beta'(u_\eps),
\end{align*}
where 
$$
H_\eps:=\div_v(a\cdot\nabla_vu)*\varGamma_\eps-\div_v(a\cdot\nabla_vu_\eps)
+[\varGamma_\eps*, v\cdot\nabla_x]u+[\varGamma_\eps*, b\cdot\nabla_v]u,
$$
and
$$
[\varGamma_\eps*, v\cdot\nabla_x]u:=\varGamma_\eps*(v\cdot\nabla_xu)-v\cdot\nabla_x(\varGamma_\eps*u).
$$
 For any $\varphi\in C^\infty_c(\mR^{1+2d})$, by the chain rule, we have
\begin{align}\label{Lim50}
\begin{split}
-\int \beta(u_\eps)\p_t\varphi
&\leq-\int\<a\cdot\nabla_vu_\eps,\nabla_v(\beta'(u_\eps)\varphi)\>
-\int (v\cdot\nabla_x \varphi)\beta(u_\eps)\\
&\quad+\int b\cdot\nabla_v \beta(u_\eps)\varphi+\int (H_\eps+f_\eps)\beta'(u_\eps)\varphi.
\end{split}
\end{align}
Note that $u$ is locally bounded and
$$
\<a\cdot\nabla_vu_\eps,\nabla_v(\beta'(u_\eps)\varphi)\>=\<a\cdot\nabla_v\beta(u_\eps),\nabla_v\varphi\>+\<a\cdot\nabla_vu_\eps,\nabla_vu_\eps\>\beta''(u_\eps)\varphi.
$$
By easy calculations, we have
$$
\lim_{\eps\to 0}\int H_\eps\beta'(u_\eps)\varphi=0,
$$
and by taking limits $\eps\downarrow 0$ for \eqref{Lim50},
\begin{align}\label{DM0}
\begin{split}
-\int \beta(u)\p_t\varphi
&\leq-\int\<a\cdot\nabla_v\beta(u),\nabla_v\varphi\>-\int \<a\cdot\nabla_vu,\nabla_vu\>\beta''(u)\varphi\\
&\quad-\int (v\cdot\nabla_x \varphi)\beta(u)+\int b\cdot\nabla_v \beta(u)\varphi+\int f\beta'(u)\varphi.
\end{split}
\end{align}
(i) For $\eps>0$, define
$$
\beta_\eps(r):=(\sqrt{(r-\eps)^2+\eps^3}+r-\eps)/2.
$$
It is easy to see that $\beta_\eps\geq 0$ is a smooth non-decreasing function and
$$
\beta_\eps''(r)\geq 0,\ \ \lim_{\eps\downarrow 0}\beta_\eps(r)=r^+,\ \ \lim_{\eps\downarrow 0}\beta'_\eps(r)=\1_{\{r>0\}},\ \ r\in\mR.
$$
Using $\beta_\eps$ in place of the $\beta$ in \eqref{DM0} and taking limits $\eps\downarrow 0$, by the dominated convergence theorem, one sees that
\begin{align*}
-\int u^+\p_t\varphi
&\leq-\int\<a\cdot\nabla_v u^+,\nabla_v\varphi\>
-\int (v\cdot\nabla_x \varphi)u^++\int b\cdot\nabla_v u^+\varphi+\int f\1_{\{u>0\}}\varphi.
\end{align*}
(ii) Note that
$$
\lim_{\eps\downarrow 0}(\beta^2_\eps(r))'=2r^+,\ \ \lim_{\eps\downarrow 0}(\beta^2_\eps(r))''=2\1_{r>0}.
$$
If $u$ is a weak solution, then \eqref{DM0} becomes an equality. 
Using $\beta^2_\eps$ in place of the $\beta$ in \eqref{DM0} and taking limits $\eps\downarrow 0$,  one sees that
\begin{align*}
-\int (u^+)^2\p_t\varphi
&=-\int\<a\cdot\nabla_v (u^+)^2,\nabla_v\varphi\>-2\int\<a\cdot\nabla_v u^+,\nabla_v u^+\>\varphi\\
&\quad-\int (v\cdot\nabla_x \varphi)(u^+)^2+\int b\cdot\nabla_v (u^+)^2\varphi+2\int fu^+\varphi.
\end{align*}
The proof is complete.
\end{proof}

\subsection{Gain of regularity via Duhamel's formula}
Below we take
\begin{align}\label{AA11}
\bba=\Big(\underbrace{3,\cdots,3}_{d},\underbrace{1,\cdots,1}_{d}\Big)\in\mR^{2d},\ \ |\bba|=4d.
\end{align}
Consider the following simple model equation:
$$
\p_t u=\Delta_vu+v\cdot\nabla_x u+f,\ \ u|_{t\leq 0}=0.
$$
By Duhamel's formula, it's unique solution can be represented by
\begin{align}\label{UU1}
u(t,z)=\int^t_{0} P_{t-s}f(s,z)\dif s,
\end{align}
where $P_t$ is the heat semigroup of operator $\Delta_v+v\cdot\nabla_x$ given by
\begin{align}\label{AA13}
P_t f(z)=\bE f\left(x+tv+\sqrt{2}\int^t_0W_s\dif s, v+\sqrt 2W_t\right),\ \ z=(x,v).
\end{align}
The following estimate is important for improving the regularity of $u$.
\bl\label{Le41}
Let $1\leq q\leq \nu\leq\infty$, $\gamma\in \mR$ and $\bbp\in[1,\infty]^{2d}$. For any $T>0$ and $\beta<\gamma+2(1+\frac1\nu-\frac1q)$,
there is a constant $C=C(T,d,\bbp,q,\nu,\gamma,\beta)>0$ such that
\begin{align}\label{ES0}
\|u\cI_t\|_{\mL^\nu_t(\bB^\beta_{\bbp;\bba})}\lesssim_C\|f\cI_t\|_{\mL^q_t(\bB^{\gamma}_{\bbp;\bba})},\ t\in(0,T],
\end{align}
where $u$ is defined by \eqref{UU1}.
Here we use the convention $f(t)|_{t\leq 0}=0$.
\el
\begin{proof}
Without loss of generality, we assume $\gamma\leq\beta$.
Since $\beta<\gamma+2(1+\frac1\nu-\frac1q)$, one can choose
$$
0<\kappa\in(1-\tfrac q\nu,(1-\tfrac{\beta-\gamma}{2})q).
$$ 
Noting that by \cite[Lemma 3.3]{ZZ21},
$$
\|P_tf\|_{\bB^{\beta}_{\bbp;\bba}}\lesssim t^{-\frac{\beta-\gamma}2}\|f\|_{\bB^{\gamma}_{\bbp;\bba}},\ \ t\in(0,T],
$$
by definition \eqref{UU1} and H\"older's inequality with respect to $(t-s)^{\kappa-1}\dif s$, we have
\begin{align*}
\|u(t)\|_{\bB^{\beta}_{\bbp;\bba}}&\lesssim\int^t_0 (t-s)^{-\frac{\beta-\gamma}2}\|f(s)\|_{\bB^{\gamma}_{\bbp;\bba}}\dif s\\
&\lesssim \left(\int^t_0(t-s)^{\kappa-1}\|f(s)\|^q_{\bB^{\gamma}_{\bbp;\bba}}\dif s\right)^{1/q}.
\end{align*}
Thus  by Minkowskii's inequality, for any $t\in[0,T]$,
\begin{align*}
\left(\int^{t}_0\|u(s)\|^\nu_{\bB^{\beta}_{\bbp;\bba}}\dif s\right)^{1/\nu}
&\lesssim \left(\int^{t}_0\left(\int^s_0(s-r)^{\kappa-1}\|f(r)\|^q_{\bB^{\gamma}_{\bbp;\bba}}\dif r\right)^{\nu/q}\dif s\right)^{1/\nu}\\
&\lesssim \left(\int^{t}_0\left(\int^{t}_r(s-r)^{\nu(\kappa-1)/q}\dif s\right)^{q/\nu}\|f(r)\|^q_{\bB^{\gamma}_{\bbp;\bba}}\dif r\right)^{1/q}\\
&\leq C_{q,\kappa}{t}^{\frac{1}{\nu}+\frac{\kappa-1}{q}} \left(\int^{t}_0\|f(r)\|^q_{\bB^{\gamma}_{\bbp;\bba}}\dif r\right)^{1/q},
\end{align*}
which gives \eqref{ES0}. The proof is complete.
\end{proof}
We also need the following simple interpolation lemma. 
\bl\label{Le22}
For any $\beta>0$, $1\leq q\leq \nu\leq\infty$ and $\bbr\in[2,\infty]^{2d}$ with 
\begin{align}\label{SB1}
\tfrac2\bbr\geq (1-\tfrac{q}\nu)\1,\ \ \bba\cdot(\tfrac1\bbr-\tfrac1{\bb2})+\tfrac{q\beta}{\nu}>0,
\end{align}
there is a constant $C=C(\beta,\nu,q,\bbr,\bba)>0$ such that
$$
\|u\|_{\mL^\nu_t(\mL^\bbr_z)}\lesssim_C \|u\|_{\mL^\infty_t(\mL^2_z)}+\|u\|_{\mL^q_t(\bB^{\beta}_{\bb2;a})}.
$$
\el
\begin{proof}
Let $\theta:=\frac q\nu\in[0,1]$. Since $\tfrac1\bbr\geq\tfrac{1-\theta}{\bb2}$, one can choose $\bbp\in[2,\infty]^{2d}$ so that
$$
\tfrac{\theta}{\bbp}+\tfrac{1-\theta}{\bb2}=\tfrac{1}{\bbr}.
$$
By H\"older's inequality and Sobolev's embedding \eqref{Sob}, we have
\begin{align*}
\|u\|_{\mL^\bbr_z}\leq\|u\|^{1-\theta}_{\mL^\bb2_z}\|u\|^{\theta}_{\mL^\bbp_z}
\lesssim\|u\|^{1-\theta}_{\mL^\bb2_z}\|u\|^{\theta}_{\bB^s_{\bbp;a}}\lesssim\|u\|^{1-\theta}_{\mL^\bb2_z}\|u\|^{\theta}_{\bB^{\beta}_{\bb2;a}},
\end{align*}
where the second inequality is due to $s=\beta-\bba\cdot(\tfrac1{\bb2}-\tfrac1{\bbp})>0$ by the assumption. Hence,
\begin{align*}
\|u\|_{\mL^\nu_t(\mL^\bbr_z)}\lesssim\|u\|^{1-\theta}_{\mL^\infty_t(\mL^\bb2_z)}\|u\|^{\theta}_{\mL^q_t(\bB^{\beta}_{\bb2;a})}.
\end{align*}
The proof is complete.
\end{proof}
\br\rm
The first condition in \eqref{SB1} means that for fixed $q$, 
if $\nu$ goes to infinity, then $\bbr$ has to be close to $\bb2$ whatever $\beta$ is how large and the second condition holds.
However, for $\beta\in(0,\frac12]$,
\begin{align}\label{SB4}
\bba\cdot(\tfrac1\bbr-\tfrac1{\bb2})+\tfrac{q\beta}{\nu}>0\Rightarrow \tfrac2\bbr\geq (1-\tfrac{q}\nu)\1.
\end{align}
Indeed, for $\bbr=(r_1,\cdots,r_{2d})\in[2,\infty]^{2d}$, and for each $i$, we have
\begin{align*}
0<\bba\cdot(\tfrac1\bbr-\tfrac1{\bb2})+\tfrac{q\beta}{\nu}\leq a_i\cdot(\tfrac1{r_i}-\tfrac1{2})+\tfrac{q}{2\nu}
\leq \tfrac1{r_i}-\tfrac1{2}+\tfrac{q}{2\nu}.
\end{align*}
\er
For any $T>0$, $q\in[1,\infty]$ and a Banach space $\mB$, we denote
$$
\mL^q_T(\mB):=L^q([0,T];\mB).
$$
Now we consider the following Cauchy problem of the kinetic equation
\begin{align}\label{DD5}
\p_t u=\div_v(a\cdot\nabla_v u)+v\cdot\nabla_x u+f,\ \ u|_{t\leq 0}\equiv0,
\end{align}
where the inhomogeneous $f$ satisfies that for any $T>0$,
\begin{enumerate}[{\bf (H$_f$)}]
\item $f=\div_v F+f_1+\cdots+f_m$, where $F\in\mL^2_T(\mL^2_z)$ and for each $i=1,\cdots,m$,
$f_i\in\mL^{q_i}_T(\mL^{\bbp_i}_z)$
for some $(q_i,\bbp_i)\in [1,2]^{1+2d}$ with
\begin{align}\label{DD95}
\tfrac1{\bbp_i}\leq(\tfrac32-\tfrac1{q_i})\1,\ \ \ \bba\cdot(\tfrac{1}{\bbp_i}-\tfrac1{\bb2})+\tfrac2{q_i}<2.
\end{align} 
\end{enumerate}

 We have the following solvability result for equation \eqref{DD5}.
\bt\label{Th3}
Under \eqref{C1} and {\bf (H$_f$)}, there is a unique weak solution $u$ to PDE \eqref{DD5}
so that for any $T>0$ and $\beta\in(0,1)$, there is a constant $C=C(T,\bba,\kappa_i,q_i,\bbp_i,\beta)>0$ such that
\begin{align}\label{DJ8}
\|u\cI_t\|_{\sV}+\|u\cI_t\|_{\mL^{2}_t(\bB^{\beta}_{2;\bba})}\lesssim_C\sum_{i=1}^m\|f_i\cI_t\|_{\mL^{q_i}_t(\mL^{\bbp_i}_z)}
+\|F\cI_t\|_{\mL^2},\ \ t\in[0,T].
\end{align}
Here we use the convention $f_i(t)=F(t)=0$ for $t\leq 0$.
\et
\begin{proof}
We only prove the apriori estimate \eqref{DJ8}. The existence follows by the standard continuity method.
Let $\chi\in C^\infty_c(\mR^{1+2d})$ with $\chi=1$ on $Q_1$ and $\chi=0$ on $Q^c_2$. For $R\geq 1$, 
let 
\begin{align}\label{Cut}
\chi_R(t,x,v):=\chi(t/R^2,x/R^3,v/R).
\end{align}
By taking $\eta=\chi_R$ in Lemma \ref{Le12}, we have for any $t\in(0,T]$,
\begin{align*}
\tfrac12\int_{\mR^{2d}}|(u\chi_R)(t)|^2&=\int_{\mR^{1+2d}}u^2\chi_R(\p_s\chi_R-v\cdot\nabla_x\chi_R)\cI_t\\
&-\int_{\mR^{1+2d}}\<a\cdot\nabla_vu, \nabla_v(u\chi_R^2)\>\cI_t
+\int_{\mR^{1+2d}}f u\chi_R^2\cI_t.
\end{align*}
Noting that for some $C>0$ independent of $R>0$,
$$
\|\p_s\chi_R-v\cdot\nabla_x\chi_R\|_\infty\leq C/R^2,
$$
by the dominated convergence theorem and letting $R\to\infty$, we obtain
\begin{align*}
\tfrac12\int_{\mR^{2d}}|u(t)|^2&=-\int_{\mR^{1+2d}}\<a\cdot\nabla_vu, \nabla_vu\>\cI_t-\int_{\cI^t\times\mR^d}F\cdot\nabla_v u\cI_t+\sum_{i=1}^m\int_{\mR^{1+2d}}f_iu\cI_t.
\end{align*}
By the uniform  ellipticity \eqref{C1} of $a$, we have
$$
\int_{\mR^{1+2d}}\<a\cdot\nabla_vu, \nabla_vu\>\cI_t\geq \kappa_0\|\nabla_v u\cI_t\|^2_{\mL^2}.
$$
By H\"older's inequality and Young's inequality, we have 
\begin{align*}
\left|\int_{\mR^{1+2d}}F\cdot \nabla_vu\cI_t\right|\leq\|F\cI_t\|_{\mL^2}\|\nabla_vu\cI_t\|_{\mL^2}
\leq\tfrac{\kappa_0}{4}\|\nabla_vu\cI_t\|_{\mL^2}^2+\tfrac{1}{\kappa_0}\|F\cI_t\|^2_{\mL^2},
\end{align*}
and for $\frac1{\bbr_i}+\frac1{\bbp_i}=1$, $\frac1{\nu_i}+\frac1{q_i}=1$ and any $\eps>0$,
$$
\left|\int_{\mR^{1+2d}}f_i u\cI_t\right|\leq\|u\cI_t\|_{\mL^{\nu_i}_t(\mL^{\bbr_i}_z)}\|f_i\cI_t\|_{\mL^{q_i}_t(\mL^{\bbp_i}_z)}
\leq \eps\|u\cI_t\|^2_{\mL^{\nu_i}_t(\mL^{\bbr_i}_z)}+\tfrac{1}{4\eps}\|f_i\cI_t\|^2_{\mL^{q_i}_t(\mL^{\bbp_i}_z)}.
$$
Combining the above calculations, we obtain that for any $\eps>0$,
\begin{align}\label{ED9}
\begin{split}
\|u\cI_t\|_{\sV}\lesssim_C 
\|F\cI_t\|_{\mL^2}+\eps\sum_{i=1}^m\|u\cI_t\|_{\mL^{\nu_i}_t(\mL^{\bbr_i}_z)}
+\eps^{-1}\sum_{i=1}^m\|f_i\cI_t\|_{\mL^{q_i}_t(\mL^{\bbp_i}_z)}.
\end{split}
\end{align}

On the other hand, we may write
$$
\p_t u=\Delta_v u+v\cdot\nabla_x u+\div_v((a-\mI)\cdot\nabla_v u+F)+\sum_{i=1}^mf_i.
$$
By Duhamel's formula, we have
\begin{align*}
u(t)&=\int^t_0P_{t-s}\div_v((a-\mI)\cdot\nabla_v u+F)\dif s+\sum_{i=1}^m\int^t_0P_{t-s}f_i\dif s
=:I_0(t)+\sum_{i=1}^mI_i(t).
\end{align*}
Below we fix $\beta\in(0,1)$ and estimate each $I_i(t)$ as following.
\begin{enumerate}[$\bullet$]
\item For $I_0(t)$,
by \eqref{ES0} with  $(\gamma,q,\nu,\bbp)=(-1,2,2,\bb2)$ and Lemma \ref{Bern}, we have
\begin{align*}
\|I_0\cI_t\|_{\mL^2_t(\bB^\beta_{2;\bba})}&\lesssim \|\div_v((a-\mI)\cdot\nabla_v u+F)\cI_t\|_{\mL^2_t(\bB^{-1}_{2;\bba})}\\
&\lesssim \|((a-\mI)\cdot\nabla_v u+F)\cI_t\|_{\mL^2}\\
&\lesssim \kappa_1\|\nabla_v u\cI_t\|_{\mL^2}+\|F\cI_t\|_{\mL^2}.
\end{align*}
\item For $I_i(t)$,  let $\beta_i>0$ be defined by
$$
\beta_i:=\beta+\bba\cdot(\tfrac1{\bbp_i}-\tfrac{1}{\bb2})<\beta+2-\tfrac2{q_i}<3-\tfrac2{q_i}.
$$
By \eqref{Sob} and \eqref{ES0} with   $(\gamma,\nu)=(0,2)$, we have for $ i=1,\cdots,m$,
\begin{align*}
\|I_i\cI_t\|_{\mL^2_t(\bB^\beta_{2;\bba})}\lesssim \|I_i\cI_t\|_{\mL^2_t(\bB^{\beta_i}_{\bbp_i;\bba})}
\lesssim \|f_i\cI_t\|_{\mL^{q_i}_t(\bB^0_{\bbp_i;\bba})}
\lesssim\|f_i\cI_t\|_{\mL^{q_i}_t(\mL^{\bbp_i}_z)}.
\end{align*}
\end{enumerate}
Combining the above calculations, we obtain
\begin{align*}
\|u\cI_t\|_{\mL^2_t(\bB^\beta_{2;\bba})}\lesssim\|\nabla_v u\cI_t\|_{\mL^2}
+\|F\cI_t\|_{\mL^2}+\sum_{i=1}^m\|f_i\cI_t\|_{\mL^{q_i}_t(\mL^{\bbp_i}_z)},
\end{align*}
which, together with \eqref{ED9},
yields that for any $\eps\in(0,1)$ and $\beta\in(0,1)$,
\begin{align}\label{DJ7}
\begin{split}
\|u\cI_t\|_{\sV}+\|u\cI_t\|_{\mL^2_t(\bB^\beta_{2;\bba})}
&\leq C\|F\cI_t\|_{\mL^2}
+\eps\sum_{i=1}^m\|u\cI_t\|_{\mL^{\nu_i}_t(\mL^{\bbr_i}_z)}
+C_\eps\sum_{i=1}^m\|f_i\cI_t\|_{\mL^{q_i}_t(\mL^{\bbp_i}_z)}.
\end{split}
\end{align}
Recalling  $\frac1{\bbr_i}+\frac1{\bbp_i}=1$ and $\frac1{\nu_i}+\frac1{q_i}=1$, by \eqref{DD95}, we have
$$
\tfrac{2}{\bbr_i}\geq (1-\tfrac2{\nu_i})\1,\ \ \bar\beta_i:=\tfrac{\nu_i}{2}[\bba\cdot(\tfrac1{\bb2}-\tfrac{1}{\bbr_i})]<1.
$$
Thus, for $\beta_i\in(\bar\beta_i,1)$, by Lemma \ref{Le22} and \eqref{DJ7}, we have for any $\eps\in(0,1)$,
\begin{align*}
&\sum_{i=1}^m\|u\cI_t\|_{\mL^{\nu_i}_t(\mL^{\bbr_i}_z)}
\lesssim \|u\cI_t\|_{\mL^\infty_t(\mL^2_z)}+\sum_{i=1}^m\|u\cI_t\|_{\mL^{2}_t(\bB^{\beta_i}_{2;\bba})}\\
&\qquad\leq C\|F\cI_t\|_{\mL^2}
+\eps\sum_{i=1}^m\|u\cI_t\|_{\mL^{\nu_i}_t(\mL^{\bbr_i}_z)}
+C_\eps\sum_{i=1}^m\|f_i\cI_t\|_{\mL^{q_i}_t(\mL^{\bbp_i}_z)}.
\end{align*}
Choosing $\eps$ small enough, we get
$$
\sum_{i=1}^m\|u\cI_t\|_{\mL^{\nu_i}_t(\mL^{\bbr_i}_z)}
\lesssim_C\|F\cI_t\|_{\mL^2}
+\sum_{i=1}^m\|f_i\cI_t\|_{\mL^{q_i}_t(\mL^{\bbp_i}_z)}.
$$
Substituting it into \eqref{DJ7}, we conclude the proof.
\end{proof}

\section{Local boundedness of weak solutions}

In this section we derive the local bounds of weak solutions for PDE \eqref{PDE0} by considering two cases: $L^p$-integrable inhomogeneous $f$ and distribution-valued
inhomogeneous $f$. 
\subsection{$L^p$-integrable inhomogeneous $f$}
We first establish the following local energy estimate.
\bl\label{Le32}
Let $Q\subset\mR^{1+2d}$ be a bounded open set and $u\in\sV_Q\cap\mL^\infty_Q$ a nonnegative weak sub-solution of PDE \eqref{PDE0} in $Q$. 
Let $(q_1,\bbp_1), (\nu_1,\bbr_1)\in(2,\infty)^{1+2d}$ satisfy
\begin{align}\label{DX5}
\bba\cdot\tfrac{1}{\bbp_1}+\tfrac2{q_1}<1,\ \tfrac{1}{\bbp_1}+\tfrac{1}{\bbr_1}=\tfrac1{\bb2},\ \tfrac{1}{q_1}+\tfrac{1}{\nu_1}=\tfrac1{2}.
\end{align}
Under \eqref{C1},
there is a constant $C=C(Q,\kappa_0,\kappa_1,\nu_1,\bbr_1)>0$ such that for any nonnegagive $\eta\in C^\infty_c(Q)$ and $t\in\mR$,
$$
\|u\eta\cI_t\|_{\sV}\lesssim_C
\Big(\|\1_Q b \|_{\mL^{q_1}_t(\mL^{\bbp_1}_z)}+\Xi_\eta\Big)\|\1_Qu\cI_t\|_{\mL^{\nu_1}_t(\mL^{\bbr_1}_z)}
+\|\<\!\!\!\<f\eta,u\eta\>\!\!\!\>\cI_t\|_{\mL^1_t}^{1/2},
$$
where $\cI_t=\1_{(-\infty,t]}$ and
\begin{align}\label{DP1}
\Xi_\eta:=\|\p_t\eta\|^{1/2}_\infty+\|v\cdot\nabla_x\eta\|^{1/2}_\infty+\|\nabla_v\eta\|_\infty.
\end{align}
\el
\begin{proof}
By \eqref{ES5}, we have for any $t\in\mR$,
\begin{align*}
\frac12\int_{\mR^{2d}}|(u\eta)(t)|^2&\leq\int_{\mR^{1+2d}}u^2\eta(\p_t\eta-v\cdot\nabla_x\eta)\cI_t-
\int_{\mR^{1+2d}}\<a\cdot\nabla_v u,\nabla_v(u\eta^2)\>\cI_t\\
&\quad+\int_{\mR^{1+2d}}(b\cdot\nabla_v u) u\eta^2\cI_t+\int_{\mR^{1+2d}}f\, u\eta^2\cI_t.
\end{align*}
Noting that by Young's inequality,
\begin{align*}
\<a\cdot\nabla_v u,\nabla_v(u\eta^2)\>&=\<a\cdot\nabla_v u,\nabla_vu\>\eta^2+
2\<a\cdot\nabla_v u,\eta\nabla_v\eta\>u\\
&\geq \kappa_0|\eta\nabla_v u|^2-\kappa_1|\eta\nabla u|\,|u\nabla_v\eta|
\geq \tfrac{\kappa_0}2|\eta\nabla_v u|^2-\tfrac{\kappa_1}{2\kappa_0}|u\nabla_v\eta|^2,
\end{align*}
we have
\begin{align*}
\int_{\cI^t\times\mR^{d}}\<a\cdot\nabla_v u,\nabla_v(u\eta^2)\>
&\geq \tfrac{\kappa_0}2\|\eta\nabla_v u\cI_t\|^2_{\mL^2}-\tfrac{\kappa_1}{2\kappa_0}\|u\nabla_v\eta\cI_t\|^2_{\mL^2}.
\end{align*}
By \eqref{DX5}, H\"older's inequality and Young's inequality, we have
\begin{align*}
\int_{\mR^{1+2d}}(b\cdot\nabla_v u) u\eta^2\cI_t
&\leq\|\eta\nabla_v u\cI_t\|_{\mL^2}\|\1_Qb\cI_t \|_{\mL^{q_1}_t(\mL^{\bbp_1}_z)}
\|\1_Qu\cI_t\|_{\mL^{\nu_1}_t(\mL^{\bbr_1}_z)}\\
&\leq\tfrac{\kappa_0}4\|\eta\nabla_v u\cI_t\|^2_{\mL^2}+\tfrac1{\kappa_0}\|\1_Qb\cI_t \|^2_{\mL^{q_1}_t(\mL^{\bbp_1}_z)}
\|\1_Qu\cI_t\|^2_{\mL^{\nu_1}_t(\mL^{\bbr_1}_z)}.
\end{align*}
Combining the above estimates, we obtain that for any $t\in\mR$,
\begin{align*}
\tfrac12\|(u\eta)(t)\|^2_{\mL^2_z}+\tfrac{\kappa_0}{4}\|\nabla_v u\eta\cI_t\|^2_{\mL^2}
&\leq\Big(\|\p_t\eta\|_\infty+\|v\cdot\nabla_x\eta\|_\infty+\tfrac{\kappa_1}{2\kappa_0}\|\nabla_v\eta\|^2_\infty\Big)\|\1_Q u\cI_t\|^2_{\mL^2}\\
&+\tfrac1{\kappa_0}\|\1_Qb\cI_t \|^2_{\mL^{q_1}_t(\mL^{\bbp_1}_z)}\|\1_Qu\cI_t\|^2_{\mL^{\nu_1}_t(\mL^{\bbr_1}_z)}
+\|\<\!\!\!\<f\eta,u\eta\>\!\!\!\>\cI_t\|_{\mL^1_t},
\end{align*}
which in turn implies the desired estimate by noting that for some $C=C(Q,\nu_1,\bbr_1)>0$,
$$
\|\1_Q u\cI_t\|^2_{\mL^2}\leq C \|\1_Qu\cI_t\|^2_{\mL^{\nu_1}_t(\mL^{\bbr_1}_z)}.
$$
The proof is complete.
\end{proof}

Recalling $Q_r$ being defined by \eqref{DD069}, we make the following assumption about the drift $b$:
\begin{enumerate}[{\bf (H$_b$)}]
\item  Suppose that $\|\1_{Q_2}b\|_{\mL^{q_1}_t(\mL^{\bbp_1}_z)}\leq\kappa_2$ for some $(q_1,\bbp_1)\in(2,\infty)^{1+2d}$ with 
\begin{align}\label{Con77}
\bba\cdot\tfrac{1}{\bbp_1}+\tfrac2{q_1}<1,\ \ \tfrac1{\bbp_1}<(\tfrac12-\tfrac1{q_1})\1.
\end{align}
\end{enumerate}
\br\rm
The second condition comes from the first condition in \eqref{SB1}.
Note that if $\bbp_1=(p_{1},\cdots,p_{1})\in(2,\infty)^{2d}$, then 
$$
\bba\cdot\tfrac{1}{\bbp_1}+\tfrac2{q_1}<1\Rightarrow\tfrac1{\bbp_1}<(\tfrac12-\tfrac1{q_1})\1.
$$
\er

Now we can show the following main result of this subsection.
\bt\label{Th33}
Under \eqref{C1} and {\bf (H$_b$)}, for any $(q_0,\bbp_0)\in(1,\infty)^{1+2d}$ with 
\begin{align}\label{Con7}
\bba\cdot\tfrac{1}{\bbp_0}+\tfrac2{q_0}<2,\ \ \tfrac1{\bbp_0}<(1-\tfrac1{q_0})\1,
\end{align}
there is a constant $C=C(\kappa_i,q_i,\bbp_i)>0$ such that for any weak sub-solution $u\in\sV_{Q_2}\cap\mL^\infty_{Q_2}$,
\begin{align}\label{LOC01}
\|\1_{Q_1}u^+\cI_t\|_{\mL^\infty}+\|\1_{Q_1}\nabla_v  u^+\cI_t\|_{\mL^2}\lesssim_C \|\1_{Q_2}u^+\cI_t\|_{\mL^2}+\|\1_{Q_2}f\cI_t\|_{\mL^{q_0}_t(\mL^{\bbp_0}_z)},\ \ 
\forall |t|\leq 1.
\end{align}
\et
\begin{proof}
Let $u\in\sV_{Q_2}\cap\mL^\infty_{Q_2}$ be a weak sub-solution and for fixed $t\in[-1,1]$,
\begin{align}\label{WQ}
\widetilde Q^t_\tau:=Q_\tau\cap((-\infty,t)\times\mR^{2d}).
\end{align}
To show the local upper bound estimate in \eqref{LOC01}, by Theorem \ref{TH20}, it suffices to show that
$$
u\in\cD\cG^+_{\sI_0}(\widetilde Q^t_\cdot)\mbox{ with $\cA=\|\1_{Q_2}f\cI_t\|_{\mL^{q_0}_t(\mL^{\bbp_0}_z)}$,}
$$ 
where $\sI_0$ is an open index subset defined by (see \eqref{SB1})
$$
\sI_0:=\Big\{(\nu,\bbr)\in(2,\infty)^{1+2d}: \tfrac2\bbr>(1-\tfrac2\nu)\1,\ \bba\cdot(\tfrac{1}\bbr-\tfrac1{\bb2})+\tfrac2{\nu}>0\Big\}.
$$
More precisely, we want to show that for any 
$(\nu,\bbr)\in\sI_0$, there is a constant $C>0$ only depending on $\kappa_0, \kappa_1, \kappa_2$, $\nu,\bbr$ and $q_i,\bbp_i, i=0,1$ such that
for any $1\leq\tau<\sigma\leq 2$ and $|t|\leq1$, $\kappa\geq 0$,
\begin{align}\label{DJ180}
\begin{split}
&(\sigma-\tau)^{2}\|\1_{Q_\tau}(u-\kappa)^+\cI_t\|_{\mL^\nu_t(\mL^\bbr_z)}
\lesssim_C\sum_{i=0,1}\|\1_{Q_\sigma}(u-\kappa)^+\cI_t\|_{\mL^{\nu_i}_t(\mL^{\bbr_i}_z)}\\
&\qquad\qquad\quad+\|\1_{Q_2}f\cI_t\|_{\mL^{q_0}_t(\mL^{\bbp_0}_z)}\|\1_{\{u>\kappa\}\cap Q_\sigma}\cI_t\|_{\mL^{\nu_0}_t(\mL^{\bbr_0}_z)},
\end{split}
\end{align}
where for $i=0,1$, $(\nu_i,\bbr_i)\in\sI_0$ is determined by $( q_i,\bbp_i)$.

\medskip

We follow the idea of \cite{GIMV19} and divide the proof into three steps.

({\sc Step 1}) Let $1\leq \tau<\sigma\leq 2$ and $\bar\tau:=\frac{\tau+\sigma}2$.  Let 
$\eta_0\in C^\infty_c(Q_{\bar\tau}; [0,1])$
and $\eta_1\in C^\infty_c(Q_{\sigma};[0,1])$ with $\eta_0|_{Q_\tau}\equiv 1$ and $\eta_1|_{Q_{\bar\tau}}\equiv 1$ and satisfy that
for some universal constant $C>0$,
\begin{align}\label{DX87}
\Xi_{\eta_i}=\|\p_t\eta_i\|^{1/2}_\infty+\|v\cdot\nabla_x\eta_i\|^{1/2}_\infty+\|\nabla_v\eta_i\|_\infty\leq C(\sigma-\tau)^{-1},\ i=0,1.
\end{align}
Since by Lemma \ref{Le213}, $(u-\kappa)^+$ is a nonnegative weak sub-solution of PDE \eqref{PDE0} in $Q_2$, without loss of generality, we may assume $u$ itself being
a nonnegative weak sub-solution of PDE \eqref{PDE0}, and prove \eqref{DJ180} for $u$ in place of $(u-\kappa)^+$. 
Note that 
\begin{align}
\p_t (u\eta_0)&\leq \div_v(a\cdot\nabla_v(u\eta_0))+v\cdot\nabla_x(u\eta_0)-\div_v(a\cdot\nabla_v\eta_0 u)\no\\
&-\<a\cdot\nabla_vu,\nabla_v\eta_0\>+(b\cdot\nabla_v u)\eta_0+(\p_t\eta_0-v\cdot\nabla_x\eta_0)u+f\1_{\{u>0\}}\eta_0\no\\
&=\div_v(a\cdot\nabla_v(u\eta_0))+v\cdot\nabla_x(u\eta_0)+\div_v F+f_1+f_2+f_3,\label{Eq1}
\end{align}
where
$$
F:=-(a\cdot\nabla_v\eta_0)u,\ f_1:=(\p_t\eta_0-v\cdot\nabla_x\eta_0)u,
$$
and
$$
f_2:=f\1_{\{u>0\}}\eta_0,\ \ f_3:=(b\eta_0-a\cdot\nabla_v\eta_0)\cdot \nabla_vu.
$$
Fix $t\in[-1,1]$. By supp($\eta_0$)$\subset Q_\sigma$, it is easy to see that
\begin{align}\label{EE2}
\|F\cI_t\|_{\mL^{2}}&=\|(a\cdot\nabla_v\eta_0)u\cI_t\|_{\mL^{2}}
\leq \kappa_1\Xi_{\eta_0}\|\1_{Q_\sigma}u\cI_t\|_{\mL^{2}},
\end{align}
and
\begin{align}\label{EE1}
\|f_1\cI_t\|_{\mL^2}
&=\|(\p_t\eta_0-v\cdot\nabla_x\eta_0)u\cI_t\|_{\mL^2}\lesssim\Xi_{\eta_0}^2\|\1_{Q_\sigma}u\cI_t\|_{\mL^2}.
\end{align}
For $f_2$, let $(\nu_0,\bbr_0)\in(2,\infty)^{1+2d}$ and $(\bar q_0,\bar\bbp_0)\in(1,2)^{1+2d}$ be defined by
\begin{align}\label{EE34}
\tfrac1{q_0}+\tfrac{2}{\nu_0}=1,\ \tfrac1{\bbp_0}+\tfrac{2}{\bbr_0}=\1,\  \tfrac{1}{\bar q_0}=\tfrac{1}{ q_0}+\tfrac1{\nu_0},\ \ \tfrac{1}{\bar\bbp_0}=\tfrac{1}{\bbp_0}+\tfrac1{\bbr_0}.
\end{align}
By H\"older's inequality, we have
\begin{align}\label{EE4}
\|f_2\cI_t\|_{\mL^{\bar q_0}_t(\mL^{\bar\bbp_0}_z)}=
\|f\1_{\{u>0\}}\eta_0\cI_t\|_{\mL^{\bar q_0}_t(\mL^{\bar\bbp_0}_z)}\leq\|f\eta_0\cI_t\|_{\mL^{ q_0}_t(\mL^{\bbp_0}_z)}\|\1_{\{u\cI_t>0\}\cap Q_\sigma}\|_{\mL^{\nu_0}_t(\mL^{\bbr_0}_z)}.
\end{align}
For $f_3$, let $(\bar q_1,\bar\bbp_1)\in(1,2)^{1+2d}$ be defined by
$$
\tfrac{1}{ q_1}+\tfrac12=\tfrac{1}{\bar q_1},\ \ \tfrac{1}{\bbp_1}+\tfrac1{\bb2}=\tfrac{1}{\bar\bbp_1}.
$$
Since $\eta_1=1$ on the support of $\eta_0$, by  H\"older's inequality and Lemma \ref{Le32},  we have
\begin{align}\label{EE3}
\begin{split}
\|f_3\cI_t\|_{\mL^{\bar q_1}_t(\mL^{\bar\bbp_1}_z)}
&=\|(b\eta_0-a\cdot\nabla_v\eta_0)\cdot \nabla_vu\eta_1\cI_t\|_{\mL^{\bar q_1}_t(\mL^{\bar\bbp_1}_z)}\\
&\leq\|b\eta_0-a\cdot\nabla_v\eta_0\|_{\mL^{ q_1}_t(\mL^{\bbp_1}_z)}\|\nabla_vu\eta_1\cI_t\|_{\mL^{2}}\\
&\lesssim\big(\kappa_2+\kappa_1\Xi_{\eta_0}\big)\Big((\kappa_2+\Xi_{\eta_1})\|\1_{Q_\sigma}u\cI_t\|_{\mL^{\nu_1}_t(\mL^{\bbr_1}_z)}
+\|fu\eta_1^2\cI_t\|_{\mL^1}^{1/2}\Big),
\end{split}
\end{align}
where $(\nu_1,\bbr_1)\in(2,\infty)^{1+2d}$ is defined by \eqref{DX5}. 
Moreover, by \eqref{EE34} and H\"older's inequality,
\begin{align*}
\|fu\eta_1^2\cI_t\|_{\mL^1}
\lesssim \|f\eta_1\cI_t\|_{\mL^{ q_0}_t(\mL^{\bbp_0}_z)}\|\1_{\{u\cI_t>0\}\cap Q_\sigma}\|_{\mL^{\nu_0}_t(\mL^{\bbr_0}_z)}\|u\eta_1\cI_t\|_{\mL^{\nu_0}_t(\mL^{\bbr_0}_z)}.
\end{align*}
Substituting this into \eqref{EE3} and by Young's inequality, we get
\begin{align}\label{EE30}
\begin{split}
\|f_3\cI_t\|_{\mL^{\bar q_1}_t(\mL^{\bar\bbp_1}_z)}
&\lesssim\big(\kappa_2+\kappa_1\Xi_{\eta_0}\big)\Big((\kappa_2+\Xi_{\eta_1})\|\1_{Q_\sigma}u\cI_t\|_{\mL^{\nu_1}_t(\mL^{\bbr_1}_z)}\\
&+\|\1_{Q_\sigma}u\cI_t\|_{\mL^{\nu_0}_t(\mL^{\bbr_0}_z)}+\|f\eta_1\cI_t\|_{\mL^{ q_0}_t(\mL^{\bbp_0}_z)}\|\1_{\{u\cI_t>0\}\cap Q_\sigma}\|_{\mL^{\nu_0}_t(\mL^{\bbr_0}_z)}\Big).
\end{split}
\end{align}
Note that by \eqref{Con7} and \eqref{Con77}, for $i=0,1$, $(\nu_i,\bbr_i)\in\sI_0$, i.e.,
$$
\tfrac2{\bbr_i}>(1-\tfrac2{\nu_i})\1,\ \bba\cdot(\tfrac{1}{\bbr_i}-\tfrac1{\bb2})+\tfrac2{\nu_i}>0.
$$

({\sc Step 2})  For $\beta\in(0,1)$ and $|t|\leq 1$, we define
$$
\sA^\beta_t:=\|w\cI_t\|_{\sV}+\|w\cI_t\|_{\mL^{2}_t(\bB^{\beta}_{2;\bba})}.
$$
Let $(\bar q_i,\bar\bbp_i)\in(1,2)^{1+2d}$, $i=0,1$ be defined as above. By \eqref{Con7} and \eqref{Con77}, one sees that
$$
\tfrac1{\bar\bbp_i}<(\tfrac32-\tfrac1{\bar q_i})\1,\ \ \ \bba\cdot(\tfrac{1}{\bar\bbp_i}-\tfrac1{\bb2})+\tfrac2{\bar q_i}<2, \ i=0,1.
$$
Thus, by Theorem \ref{Th3}, there is a unique weak solution $w$ to the following PDE:
\begin{align}\label{Eq2}
\p_t w=\div_v(a\cdot\nabla_vw)+v\cdot\nabla_x w+\div_v F+f_1+f_2+f_3,\ \ w|_{t\leq-\bar\tau^2}=0
\end{align}
so that for any $\beta\in(0,1)$ and $|t|\leq 1$,
\begin{align*}
\sA^\beta_t\lesssim
\|F\cI_t\|_{\mL^2}+\|f_1\cI_t\|_{\mL^2}+\|f_2\cI_t\|_{\mL^{\bar q_0}_t(\mL^{\bar\bbp_0}_z)}
+\|f_3\cI_t\|_{\mL^{\bar q_1}_t(\mL^{\bar\bbp_1}_z)}.
\end{align*}
Furthermore, by \eqref{EE2}, \eqref{EE1}, \eqref{EE4} and \eqref{EE30}, we get
\begin{align}\label{SZ1}
\begin{split}
\sA^\beta_t&\lesssim
(1+\Xi^2_{\eta_0})\|u\eta_1\cI_t\|_{\mL^2}+\|f\eta_0\cI_t\|_{\mL^{ q_0}_t(\mL^{\bbp_0}_z)}\|\1_{\{u\cI_t>0\}\cap Q_\sigma}\|_{\mL^{\nu_0}_t(\mL^{\bbr_0}_z)}\\
&\quad+\big(1+\Xi_{\eta_0}\big)\Big((1+\Xi_{\eta_1})\|\1_{Q_\sigma}u\cI_t\|_{\mL^{\nu_1}_t(\mL^{\bbr_1}_z)}
+\|\1_{Q_\sigma}u\cI_t\|_{\mL^{\nu_0}_t(\mL^{\bbr_0}_z)}\Big)\\
&\quad+\big(1+\Xi_{\eta_0}\big)\|f\eta_1\cI_t\|_{\mL^{ q_0}_t(\mL^{\bbp_0}_z)}\|\1_{\{u\cI_t>0\}\cap Q_\sigma}\|_{\mL^{\nu_0}_t(\mL^{\bbr_0}_z)}.
\end{split}
\end{align}
Using \eqref{DX87} and $\|u\eta_1\cI_t\|_{\mL^2}\lesssim\|\1_{Q_\sigma}u\cI_t\|_{\mL^{\nu_1}_t(\mL^{\bbr_1}_z)}$, 
we find that for any $\beta\in(0,1)$ and $1\leq\tau<\sigma\leq 2$,
\begin{align}\label{DJ18}
(\sigma-\tau)^{2}\sA^\beta_t
&\lesssim
\sum_{i=0,1}\|\1_{Q_\sigma}u\cI_t\|_{\mL^{\nu_i}_t(\mL^{\bbr_i}_z)}
+\|f\eta_1\cI_t\|_{\mL^{ q_0}_t(\mL^{\bbp_0}_z)}\|\1_{\{u\cI_t>0\}\cap Q_\sigma}\|_{\mL^{\nu_0}_t(\mL^{\bbr_0}_z)},
\end{align}
where the implicit constant is independent of $\tau,\sigma$ and $t$, but may depend on $\beta$.

({\sc Step 3})  By \eqref{Eq1} and \eqref{Eq2}, one sees that $\bar w:=u\eta_0-w\in\sV$ is a weak sub-solution of 
$$
\p_t \bar w\leq \div_v(a\cdot\nabla_vw)+v\cdot\nabla_x\bar w,\ \ \bar w|_{t\leq-\bar\tau^2}=0.
$$
By Lemma \ref{Le213}, $\bar w^+$ is still a weak sub-solution of the above equation.
By Lemma \ref{Le32} with $b=f=0$, $\eta=\chi_R$, where $\chi_R$ is the same as in \eqref{Cut}, 
and letting $R\to\infty$, it is easy to see that
$$
\bar w^+\equiv0\Rightarrow 0\leq u\eta_0\leq w.
$$
For given $(\nu,\bbr)\in\sI_0$, let
$$
\bar\beta:=\tfrac{\nu}2[\bba\cdot(\tfrac1{\bb2}-\tfrac1{\bbr})]\in(0,1).
$$
By Lemma \ref{Le22} with $q=2$ and $\beta\in(\bar\beta,1)$, we have
\begin{align*}
\|\1_{Q_\tau}u\cI_t\|_{\mL^\nu_t(\mL^\bbr_z)}
&\leq\|u\eta_0\cI_t\|_{\mL^\nu_t(\mL^\bbr_z)}\lesssim \|w\cI_t\|_{\mL^\nu_t(\mL^\bbr_z)}
\lesssim  \sA^\beta_t,
\end{align*}
which combining with \eqref{DJ18} yields \eqref{DJ180}. Thus, by Theorem \ref{TH20},
$$
\|\1_{Q_{3/2}}u^+\cI_t\|_{\mL^\infty}\lesssim_C \|\1_{Q_2}u^+\cI_t\|_{\mL^2}+\|\1_{Q_2}f\cI_t\|_{\mL^{q_0}_t(\mL^{\bbp_0}_z)}.
$$
Finally, let $\chi\in C^\infty_c(Q_{3/2})$ be nonnegative and $\chi=1$ on $Q_1$. By \eqref{LOC01} with $\eta=\chi$ and the above estimate, we also have
\begin{align*}
\|u^+\chi\cI_t\|_{\sV}&\lesssim_C\|\1_{Q_{3/2}}u^+\cI_t\|_{\mL^{\nu_1}_t(\mL^{\bbr_1}_z)}
+\|\<\!\!\!\<f\chi,u^+\chi\>\!\!\!\>\cI_t\|_{\mL^1_t}^{1/2}\\
&\lesssim_C\|\1_{Q_{3/2}}u^+\cI_t\|_{\mL^\infty}
+\|f\chi\cI_t\|_{\mL^1}^{1/2}\|\1_{Q_{3/2}}u^+\cI_t\|_{\mL^\infty}^{1/2}\\
&\lesssim_C \|\1_{Q_2}u^+\cI_t\|_{\mL^2}+\|\1_{Q_2}f\cI_t\|_{\mL^{q_0}_t(\mL^{\bbp_0}_z)},
\end{align*}
where we have used that $\|f\chi\cI_t\|_{\mL^1}\lesssim\|\1_{Q_2}f\cI_t\|_{\mL^{q_0}_t(\mL^{\bbp_0}_z)}$.
The proof is complete.
\end{proof}
\subsection{Distribution-valued inhomogeneous $f$}
In this subsection, we consider the case of $f$ being a distribution and suppose that 
\begin{enumerate}[{\bf (H$'_f$)}]
\item $f=g+\div_v G$, where $g\in\mL^{ q_0}_t(\bB^{-\alpha_0}_{\bbp_0;\bba})\cap\mL^1$
for some $\alpha_0\in[0,1)$ and $( q_0,\bbp_0)\in(1,\infty)^{1+2d}$ with
\begin{align}\label{Con6}
1-\alpha_0<\tfrac2{q_0}, \ \tfrac2{q_0}+\bba\cdot\tfrac1{\bbp_0}<2-2\alpha_0,
\end{align}
and $G\in\mL^{q_2}_t(\mL^{\bbp_2}_z)$ for some $q_2\in(2,4)$ and $\bbp_2\in(2,\infty)^{2d}$ with 
\begin{align}\label{Con60}
\tfrac2{ q_2}+\bba\cdot\tfrac1{\bbp_2}<1.
\end{align}
\end{enumerate}
\begin{enumerate}[{\bf (H$'_b$)}]
\item  $\|\1_{Q_2}b\|_{\mL^{q_1}_t(\mL^{\bbp_1}_z)}\leq\kappa_2$ for some $q_1\in(2,4)$ and $\bbp_1\in(2,\infty)^{2d}$ with 
\begin{align}\label{Con61}
\tfrac2{q_1}+\bba\cdot\tfrac{1}{\bbp_1}<1.
\end{align}
\end{enumerate}
\br\rm\label{Re44}
Condition {\bf (H$'_b$)} does not imply {\bf (H$_b$)} due to the extra requirement $\tfrac1{\bbp_1}+\frac{\1}{q_1}<\tfrac1{\bb2}$ in \eqref{Con77}.
In other words, when $q_1\in(2,4)$, we can drop the assumption $\tfrac1{\bbp_1}+\frac{\1}{q_1}<\tfrac1{\bb2}$ in {\bf (H$_b$)} in the following theorem
because we shall use Lemma \ref{Le22} with $q=4$ and $\beta\in(0,\frac12)$ (see \eqref{SB4}).
\er
\bt\label{Th34} 
Under \eqref{C1}, {\bf (H$'_b$)}  and {\bf (H$'_f$)},  
there is a constant $C>0$ only depending on $\alpha_0, \kappa_i, q_i,\bbp_i$ such that for any weak solution $u\in\sV_{Q_2}\cap\mL^\infty_{Q_2}$ of PDE \eqref{PDE0}
and $|t|\leq 1$,
\begin{align}\label{LOC11}
\begin{split}
\|\1_{Q_1}u\cI_t\|_{\mL^\infty}+\|\1_{Q_1}\nabla_v  u\cI_t\|_{\mL^2}
\lesssim_C \|\1_{Q_2}u\cI_t\|_{\mL^2}+\|g\chi_2\cI_t\|_{\mL^{ q_0}_t(\bB^{-\alpha_0}_{\bbp_0;\bba})}
+\|G\chi_2\cI_t\|_{\mL^{ q_2}_t(\mL^{\bbp_2}_z)},
\end{split}
\end{align}
and for any $\beta\in(0,1)$,
\begin{align}\label{LOC31}
\|u\chi_1\cI_t\|_{\mL^4_t(\bB^{\beta/2}_{\boldsymbol{2};\bba})}
\lesssim_C \|\1_{Q_2}u\cI_t\|_{\mL^2}+\|g\chi_2\cI_t\|_{\mL^{ q_0}_t(\bB^{-\alpha_0}_{\bbp_0;\bba})}+\|G\chi_2\cI_t\|_{\mL^{ q_2}_t(\mL^{\bbp_2}_z)},
\end{align}
where $\chi_1$ and $\chi_2$ are defined by \eqref{Cut}.
\et
\begin{proof}
Let $u\in\sV_{Q_2}\cap\mL^\infty_{Q_2}$ be a weak solution of PDE \eqref{PDE0}.
For proving \eqref{LOC11}, by Theorem \ref{TH20}, it suffices to show that for fixed $t\in[-1,1]$,
$$
u\in\cD\cG^+_{\sI_1}(\widetilde Q^t_\cdot)\mbox{ with $\cA=\|g\chi_2\cI_t\|_{\mL^{ q_0}_t(\bB^{-\alpha_0}_{\bbp_0;\bba})}+\|G\chi_2\cI_t\|_{\mL^{ q_2}_t(\mL^{\bbp_2}_z)}$,}
$$ 
where $\widetilde Q^t$ is defined by \eqref{WQ} and $\sI_1$ is an open index subset defined by
\begin{align}\label{II1}
\sI_1:=\Big\{(\nu,\bbr)\in(4,\infty)\times(2,\infty)^{2d}: \bba\cdot(\tfrac{1}\bbr-\tfrac1{\bb2})+\tfrac2{\nu}>0\Big\}.
\end{align}
More precisely, we want to show that for any 
$(\nu,\bbr)\in\sI_1$, there is a constant $C>0$ only depending on $\nu,\bbr$ and $\alpha_0$, $\kappa_i$, $ q_i,\bbp_i, i=0,1,2$ such that
for any $1\leq\tau<\sigma\leq 2$ and $|t|\leq1$, $\kappa\geq 0$,
\begin{align}\label{DX90}
\begin{split}
&(\sigma-\tau)\|\1_{Q_\tau}w\cI_t\|_{\mL^\nu_t(\mL^\bbr_z)}
\lesssim_C\|\1_{Q_\sigma}w\cI_t\|_{\mL^{\nu_1}_t(\mL^{\bbr_1}_z)}\\
&\quad+\cA\Big(\|\1_{\{w>0\}\cap Q_\sigma}\cI_t\|_{\mL^{\nu_0}_t(\mL^{\bbr_0}_z)}
+\|\1_{\{w>0\}\cap Q_\sigma}\cI_t\|_{\mL^{\nu_2}_t(\mL^{\bbr_2}_z)}\Big),
\end{split}
\end{align}
where $(\nu_i,\bbr_i)\in\sI_1$ are determined by $( q_i,\bbp_i)$, $i=0,1,2$ and 
$$
w:=(u-\kappa)^+.
$$

We divide the proof into four steps.

({\sc Step 1}) Let $1\leq \tau<\sigma\leq 2$ and $\eta\in C^\infty_c(Q_\sigma; [0,1])$ with $\eta|_{Q_\tau}\equiv 1$ and satisfy that
for some universal constant $C>0$,
\begin{align}\label{DX877}
\Xi_{\eta}:=1+\|\p_t\eta\|^{1/2}_\infty+\|v\cdot\nabla_x\eta\|^{1/2}_\infty+\|\nabla_v\eta\|_\infty\leq C(\sigma-\tau)^{-1}.
\end{align}
Note that by (ii) of Lemma \ref{Le213}, $w^2$ is a weak solution of PDE
$$
\p_tw^2=\div_v(a\cdot\nabla_v w^2)
+v\cdot\nabla_x w^2+b\cdot\nabla_v w^2-2\tr(a\cdot\nabla_v w\otimes \nabla_vw)+2fw.
$$
Let 
$$
\bar w:=w\eta.
$$
By definition, one sees that $\bar w^2$ is a nonnegative weak solution of the following PDE
\begin{align}\label{AA8}
\p_t\bar w^2=\div_v(a\cdot\nabla_v\bar w^2)+v\cdot\nabla_x\bar w^2+h,
\end{align}
where, thanks to $f=g+\div_v G$, $h$ is given by
\begin{align*}
h&:=-\div_v(a\cdot\nabla_v\eta^2 w^2)-\tr(a\cdot\nabla_v w^2\otimes\nabla_v\eta^2)-2\tr(a\cdot\nabla_v w\otimes \nabla_vw)\eta^2\\
&\qquad+(\p_t\eta^2-v\cdot\nabla_x\eta^2)w^2+(b\cdot\nabla_v w^2)\eta^2+2(g+\div_v G)w\eta^2.
\end{align*}
Integrating both sides of \eqref{AA8} over $\mR^{2d}$, by the integration by parts, we obtain
\begin{align}\label{AA09}
\p_t\int\bar w^2=\int h=\int f_1+\int f_2+\int f_3+2\<\!\!\!\<g\eta,\bar w\>\!\!\!\>,
\end{align}
where $\<\!\!\!\<\cdot,\cdot\>\!\!\!\>$ stands for the duality between $\mL^1_z$ and $\mL^\infty_z$, and for $\bar b:=b\eta-2a\cdot\nabla\eta$,
\begin{align*}
f_1&:=-2\tr(a\cdot\nabla_vw\otimes \nabla_vw)\eta^2,\\
f_2&:=(\p_t\eta^2-v\cdot\nabla_x\eta^2)w^2+(\bar b\cdot\nabla_v w^2)\eta,\\
f_3&:=-2G\cdot(\nabla_vw\eta^2+\nabla\eta^2 w).
\end{align*}
For $f_1$, by the uniform ellipticity of $a$, we have
\begin{align*}
\int f_1&=-\int\tr(a\cdot\nabla_v w\otimes\nabla_v w)\eta^2
\leq-\kappa_0\int|\nabla_vw\eta|^2=-\kappa_0\|\nabla_vw\eta\|_{\mL^2_z}^2.
\end{align*}
For $f_2$, since $\eta(t,z)\in[0,1]$ has support in $Q_\sigma$, by H\"older and Young's inequalities, we have
\begin{align*}
\int |f_2|&\leq 2\|(|\p_t\eta|+|v\cdot\nabla_x\eta|)w^2\|_{\mL^1_z}+2\int |(\bar b\cdot\nabla_v w) w\eta|\\
&\leq\Xi^2_\eta \|\1_{Q_\sigma}w\|_{\mL^2_z}^2+\tfrac{\kappa_0}4\|\nabla_vw\eta\|_{\mL^2_z}^2+\tfrac{2}{\kappa_0}\|\bar b w\|^2_{\mL^2_z}.
\end{align*}
For $f_3$, by $\nabla_v w=\nabla_v w\1_{w\not=0}$ (see  \cite[Lemma 7.6]{Gi}) and H\"older's inequality, we have
\begin{align}\label{SZ2}
\begin{split}
\int |f_3|&\leq 2\|G\1_{w\not=0}\eta\|_{\mL^2_z}(\|\nabla_v w\eta\|_{\mL^2_z}+2\|w\nabla_v\eta\|_{\mL^2_z})\\
&\leq \tfrac{\kappa_0}4\|\nabla_vw\eta\|^2_{\mL^2_z}+\tfrac{\kappa_0}2\|w\nabla_v\eta\|_{\mL^2_z}^2+\tfrac{4}{\kappa_0}\|G\1_{w\not=0}\eta\|_{\mL^2_z}^2.
\end{split}
\end{align}
Integrating both sides of \eqref{AA09} from $-4$ to $t$ and combining the above calculations, we obtain
\begin{align*}
\|\bar w(t)\|^2_{\mL^{2}_{z}}&\leq\int^t_{-4}\Big(-\tfrac{\kappa_0}2\|\nabla_v w\eta\|^2_{\mL^{2}_{z}}
+C\Xi^2_\eta\|\1_{Q_\sigma} w\|_{\mL^2_z}^2+\tfrac{2}{\kappa_0}\|\bar b w\|^2_{\mL^2_z}+\tfrac{4}{\kappa_0}\|G\1_{w\not=0}\eta\|_{\mL^2_z}^2+2|\<\!\!\!\<g\eta,\bar w\>\!\!\!\>|\Big)\dif s,
\end{align*}
which yields by taking supremum in $t$ that for fixed $t\in[-1,1]$,
\begin{align}\label{DD7}
\|\bar w\cI_t\|^2_{\mL^\infty_t(\mL^2_z)}+\|\nabla_v  w\cI_t\eta\|^2_{\mL^2}
\lesssim\Xi^2_\eta\|\1_{Q_\sigma} w\cI_t\|^2_{\mL^2}+\|\bar b w\cI_t\|^2_{\mL^2}+\|G\1_{w\not=0}\eta\|_{\mL^2}^2+\|\<\!\!\!\<g\eta,\bar w\>\!\!\!\>\cI_t\|_{\mL^1_t}.
\end{align}

({\sc Step 2}) Note that in the distributional sense,
\begin{align*}
\p_t\bar w^2&=\Delta_v\bar w^2+v\cdot\nabla_x\bar w^2+\div_vF+f_1+f_2+f_3+2gw\eta^2,
\end{align*}
where $f_1,f_2,f_3$ are defined as above and
$$
F:=(a-\mI)\cdot\nabla_v\bar w^2-a\cdot\nabla_v\eta^2 w^2+Gw\eta^2.
$$ 
By Duhamel's formula, we have
\begin{align*}
\bar w(t)^2&=\int^t_{-4}P_{t-s}\big(\div_v F+f_1+f_2+f_3+2g\eta\bar w\big)\dif s\\
&=:I_0(t)+I_1(t)+I_2(t)+I_3(t)+I_4(t).
\end{align*}
We estimate each term by Lemma \ref{Le41} as following.
\begin{enumerate}[$\bullet$]
\item For $I_0(t)$, by \eqref{ES0} with $\beta\in(0,1)$ and $(\nu,q,\gamma,\bbp)=(2,2,-1,\boldsymbol{1})$,  
\begin{align*}
\|I_0\cI_t\|_{\mL^2_t(\bB^{\beta}_{\boldsymbol{1};\bba})}&\lesssim\|\div_v((a-\mI)\cdot\nabla_v\bar w^2-a\cdot\nabla_v\eta^2 w^2+Gw\eta^2)
\cI_t\|_{\mL^2_t(\bB^{-1}_{\boldsymbol{1};\bba})}\\
&\lesssim\|((a-\mI)\cdot\nabla_v\bar w^2-a\cdot\nabla_v\eta^2 w^2+Gw\eta^2)\cI_t\|_{\mL^2_t(\mL^1_z)}\\
&\lesssim\|\bar w\nabla_v\bar w\cI_t\|_{\mL^2_t(\mL^1_z)}+\|w\nabla_v\eta \bar w\cI_t\|_{\mL^2_t(\mL^1_z)}+\|G\bar w\eta\cI_t\|_{\mL^2_t(\mL^1_z)}\\
&\lesssim\big(\|\nabla_v\bar w\cI_t\|_{\mL^2}+\|w\nabla_v\eta\cI_t\|_{\mL^2}+\|G\1_{w\not=0}\eta\|_{\mL^2}\big)\|\bar w\cI_t\|_{\mL^\infty_t(\mL^2_z)}\\
&\lesssim\|\nabla_v w\eta\cI_t\|^2_{\mL^2}+\|w\nabla_v\eta\cI_t\|^2_{\mL^2}+\|G\1_{w\not=0}\eta\|_{\mL^2}^2+\|\bar w\cI_t\|_{\mL^\infty_t(\mL^2_z)}^2.
\end{align*}
\item For $I_1(t)$, by \eqref{ES0} with $\beta\in(0,1)$ and $(\nu,q,\gamma,\bbp)=(2,1,0,\boldsymbol{1})$, 
\begin{align*}
\|I_1\cI_t\|_{\mL^2_t(\bB^{\beta}_{\boldsymbol{1};\bba})}&\lesssim\|\tr(a\cdot\nabla_v w\otimes \nabla_vw)\eta^2\cI_t\|_{\mL^1_t(\bB^0_{\boldsymbol{1};\bba})}\\
&\lesssim\|\tr(a\cdot\nabla_v w\otimes \nabla_vw)\eta^2\cI_t\|_{\mL^1}
\lesssim\|\nabla_vw\eta\cI_t\|_{\mL^2}^2.
\end{align*}
\item For $I_2(t)$, by \eqref{ES0} with $\beta\in(0,1)$ and $(\nu,q,\gamma,\bbp)=(2,1,0,\boldsymbol{1})$, 
\begin{align*}
\|I_2\cI_t\|_{\mL^2_t(\bB^{\beta}_{\boldsymbol{1};\bba})}&\lesssim\|((\p_t\eta^2-v\cdot\nabla_x\eta^2)w^2-(\bar b\cdot\nabla_v w^2)\eta)\cI_t\|_{\mL^1}\\
&\lesssim\Xi^2_\eta \|\1_{Q_\sigma}w\cI_t\|^2_{\mL^2}+\|\bar b w\cI_t\|_{\mL^2}\|\nabla_vw\eta\cI_t\|_{\mL^2}.
\end{align*}
\item For $I_3(t)$, by \eqref{ES0}  with $\beta\in(0,1)$ and  $(\nu,q,\gamma,\bbp)=(2,1,0,1)$, as in \eqref{SZ2},
\begin{align*}
\|I_3\cI_t\|_{\mL^2_t(\bB^{\beta}_{\boldsymbol{1};\bba})}&\lesssim\|G\cdot(\nabla_vw\eta^2+\nabla\eta^2 w)\cI_t\|_{\mL^1}\\
&\leq\|G\1_{w\not=0}\eta\|_{\mL^2}(\|\nabla_v w\eta\|_{\mL^2}+2\|w\nabla_v\eta\|_{\mL^2}).
\end{align*}
\item For $I_4(t)$, by \eqref{ES0}  with $\beta\in(0,1)$ and  $(\nu,q,\gamma,\bbp)=(2,\frac{2}{2-\alpha_0},-\alpha_0,\1)$, 
\begin{align*}
\|I_4\cI_t\|_{\mL^2_t(\bB^{\beta}_{\boldsymbol{1};\bba})}\lesssim\|g\eta\bar w\cI_t\|_{\mL^{q}_t(\bB^{-\alpha_0}_{\boldsymbol{1};\bba})}.
\end{align*}
\end{enumerate}
Combining the above calculations, we obtain that for $q=\frac2{2-\alpha_0}$ and any $\beta\in(0,1)$,
\begin{align}\label{EE6}
\begin{split}
\|\bar w^2\cI_t\|_{\mL^2_t(\bB^{\beta}_{\boldsymbol{1};\bba})}
&\lesssim \|\bar w\cI_t\|^2_{\mL^\infty_t(\mL^2_z)}
+\|\nabla_v  w\eta\cI_t\|^2_{\mL^2}+\|G\1_{w\not=0}\eta\cI_t\|_{\mL^2}^2\\
&+\Xi^2_\eta \|\1_{Q_\sigma}w\cI_t\|^2_{\mL^2}+\|\bar b w\cI_t\|_{\mL^2}\|\nabla_vw\eta\cI_t\|_{\mL^2}
+\|g\eta\bar w\cI_t\|_{\mL^{q}_t(\bB^{-\alpha_0}_{\boldsymbol{1};\bba})}.
\end{split}
\end{align}
Note that by  Lemma \ref{Le24},
\begin{align}\label{EE06}
\|\bar w\cI_t\|^2_{\mL^4_t(\bB^{\beta/2}_{\boldsymbol{2};\bba})}\lesssim \|\bar w^2\cI_t\|_{\mL^2_t(\bB^{\beta}_{\boldsymbol{1};\bba})}.
\end{align}
If we define
$$
\sA^\beta_t:=\|\bar w\cI_t\|^2_{\mL^4_t(\bB^{\beta/2}_{\boldsymbol{2};\bba})}+\|\bar w\cI_t\|^2_{\mL^\infty_t(\mL^2_z)}+\|\nabla_v  w\eta\cI_t\|^2_{\mL^2},
$$
then by \eqref{EE06}, \eqref{EE6} and \eqref{DD7}, we get  for $q=\frac2{2-\alpha_0}$ and any $\beta\in(0,1)$, $|t|<1$,
\begin{align}\label{AA0}
\sA^\beta_t\lesssim\Xi^2_\eta \|\1_{Q_\sigma}w\cI_t\|^2_{\mL^2}+\|\bar b w\cI_t\|^2_{\mL^2}+\|G\1_{w\not=0}\eta\cI_t\|_{\mL^2}^2
+\|\<\!\!\!\<g\eta,\bar w\>\!\!\!\>\cI_t\|_{\mL^1_t}+\|g\eta\bar w\cI_t\|_{\mL^{q}_t(\bB^{-\alpha_0}_{\boldsymbol{1};\bba})}.
\end{align}

({\sc Step 3}) In this step we estimate the last two terms by Lemma \ref{Le23} and H\"older's inequality.
Let $ q_2\in(0,\infty)$ and $\bbp_2\in(1,2]^{2d}$ be defined by
\begin{align}\label{DD8}
\tfrac1{ q_0}+\tfrac1{ q_2}=\tfrac1{q}=1-\tfrac{\alpha_0}2,\ \ \tfrac1{\bbp_0}+\tfrac1{\bbp_2}=1.
\end{align}
By Lemma \ref{Le23} and H\"older's inequality, we have for any $\alpha_2>\alpha_0$,
\begin{align}\label{DD6}
\|\<\!\!\!\<g\eta,\bar w\>\!\!\!\>\cI_t\|_{\mL^1_t}+\|g\eta\bar w\cI_t\|_{\mL^{q}_t(\bB^{-\alpha_0}_{\boldsymbol{1};\bba})}\lesssim
\|g\eta\cI_t\|_{\mL^{ q_0}_t(\bB^{-\alpha_0}_{\bbp_0;\bba})}\|\bar w\cI_t\|_{\mL^{ q_2}_t(\bB^{\alpha_2}_{\bbp_2;\bba})}.
\end{align}
Now due to \eqref{Con6}, one can choose $\theta\in(2\alpha_0,1)$ being close to $2\alpha_0$ so that
\begin{align}\label{DD91}
\tfrac{2}{ q_0}+\bba\cdot\tfrac1{\bbp_0}<2-\alpha_0-\tfrac\theta2.
\end{align}
Let $ q_3\in(0,\infty)$ and $\bbp_3\in(1,2)^{2d}$ be defined by
\begin{align}\label{DD9}
\tfrac1{ q_2}=\tfrac{1-\theta}{ q_3}+\tfrac{\theta}{4},\ \ \tfrac1{\bbp_2}=\tfrac{1-\theta}{\bbp_3}+\tfrac{\theta}{\bb2}.
\end{align}
In particular, since $\frac\theta2>\alpha_0$, one can choose $\beta\in(0,1)$ so that
$$
\alpha_2:=\bba\cdot\tfrac{1}{\bbp_2}+(1-\theta)(0-\bba\cdot\tfrac{1}{\bbp_3})+\theta(\tfrac\beta 2-\tfrac{|\bba|}{2})=\tfrac{\theta\beta}2>\alpha_0.
$$
Thus by Gagliado-Nirenberge's inequality \eqref{Sob} and H\"older's inequality, we have
\begin{align*}
\|\bar w\cI_t\|_{\mL^{ q_2}_t(\bB^{\alpha_2}_{\bbp_2;\bba})}
&\lesssim \|\bar w\cI_t\|^{1-\theta}_{\mL^{ q_3}_t(\bB^0_{\bbp_3;\bba})}\|\bar w\cI_t\|^{\theta}_{\mL^{4}_t(\bB^{\beta/2}_{2;\bba})}
\lesssim\|\bar w\cI_t\|^{1-\theta}_{\mL^{ q_3}_t(\mL^{\bbp_3}_z)}\|\bar w\cI_t\|^{\theta}_{\mL^{4}_t(\bB^{\beta/2}_{2;\bba})}.
\end{align*}
On the other hand, for $\nu_0:=\frac{(2-\theta) q_3}{1-\theta}$ and $\bbr_0:=\frac{(2-\theta)\bbp_3}{1-\theta}$, by H\"older's inequality,
\begin{align*}
\|\bar w\cI_t\|^{1-\theta}_{\mL^{ q_3}_t(\mL^{\bbp_3}_z)}
&\leq \|\1_{\bar w\not=0}\cI_t\|^{1-\theta}_{\mL^{\nu_0(1-\theta)}_t(\mL^{\bbr_0(1-\theta)}_z)}
\|\bar w\cI_t\|^{1-\theta}_{\mL^{\nu_0}_t(\mL^{\bbr_0}_z)}\\
&=\|\1_{\bar w\not=0}\cI_t\|_{\mL^{\nu_0}_t(\mL^{\bbr_0}_z)}\|\bar w\cI_t\|^{1-\theta}_{\mL^{\nu_0}_t(\mL^{\bbr_0}_z)}.
\end{align*}
Hence,
\begin{align}\label{AA1}
\|\bar w\cI_t\|_{\mL^{ q_2}_t(\bB^{\alpha_2}_{\bbp_2;\bba})}
\lesssim \|\1_{\bar w\not=0}\cI_t\|_{\mL^{\nu_0}_t(\mL^{\bbr_0}_z)}
\|\bar w\cI_t\|^{1-\theta}_{\mL^{\nu_0}_t(\mL^{\bbr_0}_z)}\|\bar w\cI_t\|^{\theta}_{\mL^{4}_t(\bB^{\beta/2}_{2;\bba})}.
\end{align}
Since $1-\alpha_0<\tfrac2{q_0}$, by \eqref{DD91}, it is easy to see that $(\nu_0,\bbr_0)\in\sI_1$.

({\sc Step 4})  Combining \eqref{AA0}, \eqref{DD6} and \eqref{AA1} and by H\"older's inequality, we have 
\begin{align*}
\sA^\beta_t&\lesssim\Xi^2_\eta \|\1_{Q_\sigma} w\cI_t\|^2_{\mL^2}+\|\bar b w\cI_t\|^2_{\mL^2}+\|G\1_{w\not=0}\eta\cI_t\|_{\mL^2}^2
+\|g\eta\cI_t\|_{\mL^{ q_0}_t(\bB^{-\alpha_0}_{\bbp_0;\bba})}\|\bar w\cI_t\|_{\mL^{ q_2}_t(\bB^{\alpha_1}_{\bbp_2;\bba})}\\
&\lesssim \Xi^2_\eta \|\1_{Q_\sigma} w\cI_t\|^2_{\mL^2}+\|\bar b\cI_t\|^2_{\mL^{ q_1}_t(\mL^{\bbp_1}_z)}\|\1_{Q_\sigma}w\cI_t\|^2_{\mL^{\nu_1}_t(\mL^{\bbr_1}_z)}
+\|G\eta\cI_t\|^2_{\mL^{ q_2}_t(\mL^{\bbp_2}_z)}\|\1_{\bar w\not=0}\cI_t\|^2_{\mL^{\nu_2}_t(\mL^{\bbr_2}_z)}\\
&\quad+\|g\eta\cI_t\|_{\mL^{ q_0}_t(\bB^{-\alpha_0}_{\bbp_0;\bba})}
\|\1_{\bar w\not=0}\cI_t\|_{\mL^{\nu_0}_t(\mL^{\bbr_0}_z)}
\|\bar w\cI_t\|^{1-\theta}_{\mL^{\nu_0}_t(\mL^{\bbr_0}_z)}\|\bar w\cI_t\|^{\theta}_{\mL^{4}_t(\bB^{\beta/2}_{2;\bba})},
\end{align*}
where we have used $\1_{w\not=0}\eta=\1_{\bar w\not=0}\eta$, and $(\nu_i,\bbr_i)\in(2,\infty)^{1+2d}$, $i=1,2$, are defined by 
$$
\tfrac{1}{\bbp_i}+\tfrac{1}{\bbr_i}=\tfrac1{\bb2},\ \tfrac{1}{q_i}+\tfrac{1}{\nu_i}=\tfrac1{2}.
$$
Note that by \eqref{Con60} and \eqref{Con61},
$$
(\nu_i,\bbr_i)\in\sI_1,\ \ i=1,2.
$$
Recalling $\bar b=b\eta-2a\cdot\nabla\eta$,
for any $\eps\in(0,1)$, by Young's inequality, we arrive at
\begin{align}\label{DX99}
\sA^\beta_t&\lesssim \Xi^2_\eta \|\1_{Q_\sigma}w\cI_t\|^2_{\mL^2}+\Xi^2_\eta\|\1_{Q_\sigma} w\cI_t\|^2_{\mL^{\nu_1}_t(\mL^{\bbr_1}_z)}
+\|G\eta\cI_t\|^2_{\mL^{ q_2}_t(\mL^{\bbp_2}_z)}\|\1_{\bar w\not=0}\cI_t\|^2_{\mL^{\nu_2}_t(\mL^{\bbr_2}_z)}\no\\
&\quad+C_\eps\|g\eta\cI_t\|^2_{\mL^{ q_0}_t(\bB^{-\alpha_0}_{\bbp_0;\bba})} \|\1_{\bar w\not=0}\cI_t\|^2_{\mL^{\nu_0}_t(\mL^{\bbr_0}_z)}
+\eps\Big(\|\bar w\cI_t\|^2_{\mL^{\nu_0}_t(\mL^{\bbr_0}_z)}+\|\bar w\cI_t\|^2_{\mL^{4}_t(\bB^{\beta/2}_{2;\bba})}\Big).
\end{align}
For any $(\nu,\bbr)\in\sI_1$, one can choose $\beta\in(0,1)$ so that $\bba\cdot(\tfrac1\bbr-\tfrac1{\bb2})+\tfrac{2\beta}{\nu}>0$,
and thus, by Lemma \ref{Le22} with $q=4$ and \eqref{SB4},
\begin{align}\label{DX909}
\|\bar w\cI_t\|^2_{\mL^\nu_t(\mL^\bbr_z)}\lesssim \sA^{\beta}_t,
\end{align}
which, together with $(\nu_0,\bbr_0)\in\sI_1$, implies that 
the last term in \eqref{DX99} can be absorbed by the left hand side if one lets $\eps$ being small enough.
Moreover,
$$
\|G\eta\cI_t\|_{\mL^{ q_2}_t(\mL^{\bbp_2}_z)}
=\|G\chi_2\eta\cI_t\|_{\mL^{ q_2}_t(\mL^{\bbp_2}_z)}\leq \|G\chi_2\cI_t\|_{\mL^{ q_2}_t(\mL^{\bbp_2}_z)},
$$
and
$$
\|g\eta\|_{\mL^{ q_0}_t(\bB^{-\alpha_0}_{\bbp_0;a})}
=\|g\chi_2\eta\|_{\mL^{ q_0}_t(\bB^{-\alpha_0}_{\bbp_0;a})}\lesssim\|g\chi_2\|_{\mL^{ q_0}_t(\bB^{-\alpha_0}_{\bbp_0;a})}
\|\eta\|_{\mL^{\infty}_t(\bB^1_{\infty;a})}.
$$
Substituting these into \eqref{DX99} and by \eqref{DX877}, \eqref{DX909} and 
$\|\1_{Q_\sigma}w\cI_t\|_{\mL^2}\lesssim\|\1_{Q_\sigma}w\cI_t\|_{\mL^{\nu_1}_t(\mL^{\bbr_1}_z)}$, we obtain \eqref{DX90} as well as \eqref{LOC31}.
\end{proof}
\br\label{Re46}\rm
Since $\nu_0, \nu_1,\nu_2$ in \eqref{DX90} are required to be greater than $4$ (see \eqref{II1}), we have to require 
$1-\alpha_0<\tfrac2{ q_0}$ in  \eqref{Con6} and $ q_1,q_2\in(2,4)$ in \eqref{Con61} and \eqref{Con60}. 
Note that $ q_i\in(2,4)$ naturally leads to $\bba\cdot\tfrac1{\bbp_i}<\frac12$, and when $\alpha_0=0$, the condition $1<\tfrac2{ q_0}$
leads to $\bba\cdot\tfrac1{\bbp_0}<1$. 
\er

\br\label{Re45}\rm
If $g=\sum_{i=0}^m g_i$, where $g_i\in\mL^{ q_i}_t(\bB^{-\alpha_i}_{\bbp_i;\bba})$ for some $\alpha_i\in[0,1)$ and $(q_i,\bbp_i)\in(1,\infty)^{1+2d}$ satisfying \eqref{Con6}, then
we still have \eqref{LOC11} with $\|g\chi_2\cI_t\|_{\mL^{ q_0}_t(\bB^{-\alpha_0}_{\bbp_0;\bba})}$ replaced by 
$\sum_{i=0}^m\|g_i\chi_2\cI_t\|_{\mL^{ q_i}_t(\bB^{-\alpha_i}_{\bbp_i;\bba})}$.
\er

\section{Global boundness and stability of weak solutions}

Fix $r>0$. For $\bbp\in[1,\infty]^{2d}$, we introduce the following localized $L^\bbp$-space:
\begin{align}\label{LOC101}
\nor f\nor_{\widetilde\mL^\bbp_z}:=\sup_{z_0}\|f\1_{B^\bba_r(z_0)}\|_{\mL^\bbp_z}=\sup_{z_0}\|f(\cdot-z_0)\1_{B^\bba_r}\|_{\mL^\bbp_z},
\end{align}
and for $(q,\bbp)\in[1,\infty]^{1+2d}$, 
\begin{align}\label{LOC1}
\nor f\nor_{\widetilde\mL^{q}_t(\widetilde\mL^\bbp_z)}
:=\sup_{t_0,z_0}\|f(\cdot-t_0,\cdot-z_0)\1_{Q_r}\|_{\mL^q_t(\mL^\bbp_z)},
\end{align}
and for $\alpha\in\mR$,
$$
\nor f\nor_{\widetilde\mL^{q}_t(\widetilde\bB^{\alpha}_{\bbp;\bba})}:=\sup_{t_0,z_0}\|f(\cdot-t_0,\cdot-z_0)\chi_2\|_{\mL^{q}_t(\bB^{\alpha}_{\bbp;\bba})},
$$
where $\chi_2$ is a cutoff function as in \eqref{Cut}.
By a finitely covering technique, it is easy to see that for different $r,r'>0$ (cf. \cite{ZZ21}),
\begin{align}\label{DS1}\sup_{z_0}\|f\1_{B^\bba_r(z_0)}\|_{\mL^\bbp_z}\asymp \sup_{z_0}\|f\1_{B^\bba_{r'}(z_0)}\|_{\mL^\bbp_z}\end{align}and
$$
\sup_{t_0,z_0}\|f(\cdot-t_0,\cdot-z_0)\1_{Q_r}\|_{\mL^q_t(\mL^\bbp_z)}\asymp \sup_{t_0,z_0}\|f(\cdot-t_0,\cdot-z_0)\1_{Q_{r'}}\|_{\mL^q_t(\mL^\bbp_z)}.
$$
In particular, for any $T>0$,
\begin{align}\label{HD8}
\nor \1_{[0,T]}f\nor_{\widetilde\mL^{q}_t(\widetilde\mL^\bbp_z)}\asymp \sup_{z_0}\|\1_{[0,T]}f\1_{B^\bba_r(z_0)}\|_{\mL^q_t(\mL^\bbp_z)}
\leq \left(\int^T_0\nor f(s)\nor_{\widetilde\mL^\bbp_z}^q\dif s\right)^{1/q},
\end{align}
and for any bounded $Q\subset\mR^{1+2d}$, there is a constant $C=C(Q,d,q,\bbp)>0$ such that
\begin{align}\label{HD9}
\|\1_Qf\|_{\mL^q_t(\mL^\bbp_z)}\lesssim_C \nor f\nor_{\widetilde\mL^{q}_t(\widetilde\mL^\bbp_z)}.
\end{align}
We also introduce the following localized energy space for later use:
$$
\widetilde\sV:=\Big\{ f\in\mL^1_{loc}: \nor f\nor_{\widetilde\sV}:=\nor f\nor_{\widetilde\mL^\infty_t(\widetilde\mL^2_z)}+
\nor \nabla_v f\nor_{\widetilde\mL^2_t(\widetilde\mL^2_z)}<\infty\Big\}.
$$

Using the above notations, we make the following global assumption about drift $b$:
\begin{enumerate}[{\bf ($\widetilde {\bf H}_b$)}]
\item  Suppose that $\nor b\nor_{\widetilde\mL^{q_1}_t(\widetilde\mL^{\bbp_1}_z)}\leq\kappa_2$
for some $(q_1,\bbp_1)\in(2,\infty)^{1+2d}$ satisfying \eqref{Con77}.
\end{enumerate}

First of all we derive the following existence and uniqueness result of global weak solutions for FPKE \eqref{PDE0} with initial value $u(t)|_{t\leq 0}=0$. 
\bt\label{Th41}
Under \eqref{C1} and {\bf ($\widetilde{\bf H}_b$)}, for any $f\in \widetilde\mL^{q_0}_t(\widetilde\mL^{\bbp_0}_z)$
with $(q_0,\bbp_0)\in(1,\infty)^{1+2d}$ satisfying \eqref{Con7},
and $T>0$, there is a unique global weak solution $u$ to PDE \eqref{PDE0} in the sense of Definition \ref{Def31} 
with initial value $u(t)|_{t\leq 0}=0$ such that for some $C=C(T,\kappa_i,q_i,\bbp_i)>0$ and any $t\in[0,T]$,
\begin{align}\label{EA80}
\|u\cI_t\|_{\mL^\infty}+\nor u\cI_t\nor_{\widetilde\sV}\lesssim_C\nor f\cI_t\nor_{\widetilde\mL^{q_0}_t(\widetilde\mL^{\bbp_0}_z)}.
\end{align}
\et
\begin{proof}
We divide the proof into two steps.

({\sc Step 1}) We first show \eqref{EA80} for any weak solution $u$. 
Without loss of generality, we assume $T=1$ and 
$$
u(t,x,v)=b(t,x,v)=f(t,x,v)\equiv 0, \ \  \forall t\leq 0. 
$$
Fix $z_0=(x_0,v_0)\in\mR^{2d}$. Let $\Gamma_tz_0:=(x_0+tv_0,v_0)$ and define
$$
u_{z_0}(t,z):=u(t,z-\Gamma_tz_0),\ a_{z_0}(t,z):=a(t,z-\Gamma_tz_0),\ 
$$
and
$$
b_{z_0}(t,z):=b(t,z-\Gamma_tz_0),\ f_{z_0}(t,z):=f(t,z-\Gamma_tz_0).
$$
By definition, it is easy to see that $u_{z_0}$ is a weak solution of PDE \eqref{PDE0} 
with coefficients $(a_{z_0}, b_{z_0}, f_{z_0})$.
By {\bf ($\widetilde {\bf H}^{q_1}_b$)} and \eqref{HD9}, it is easy to see that
$$
\sup_{z_0}\|\1_{Q_2} b_{z_0}\|_{\mL^{q_1}_t(\mL^{\bbp_1}_z)}\leq C(\kappa_2).
$$
By applying Theorem \ref{Th33} to nonnegative weak sub-solution $u^+_{z_0}$ and  $u^-_{z_0}=(-u_{z_0})^+$ separately, 
there is a constant $C=C(\kappa_i,q_i,\bbp_i)>0$ such that for each $t\in[0,1]$ and $z_0\in\mR^{2d}$,
\begin{align}\label{DD50}
\|\1_{Q_1}u_{z_0}\cI_t\|_{\mL^\infty}+\|\1_{Q_1}\nabla_v  u_{z_0}\cI_t\|_{\mL^2}\lesssim_C\|\1_{Q_2}u_{z_0}\cI_t\|_{\mL^2}+
\|\1_{Q_2}f_{z_0}\cI_t\|_{\mL^{q_0}_t(\mL^{\bbp_0}_z)}.
\end{align}
In particular, for each $t\in[0,1]$,
\begin{align*}
\|\1_{Q_1} u_{z_0}(t)\|^2_{\mL^2_z}&\leq {\rm Vol}(Q_1)\|\1_{Q_1}u_{z_0}(t)\|_{\mL^\infty_z}^2
\lesssim\int^t_0\|\1_{Q_2}u_{z_0}(s)\|_{\mL^2_z}^2\dif s+
\|\1_{Q_2}f_{z_0}\cI_t\|^2_{\mL^{q_0}_t(\mL^{\bbp_0}_z)}.
\end{align*}
Taking supremum in $z_0\in\mR^{2d}$ and by \eqref{DS1} and \eqref{HD8}, we obtain
\begin{align*}
\nor u(t)\nor^2_{\mL^2_z}\lesssim_C\int^t_0\nor u(s)\nor_{\mL^2_z}^2\dif s+
\nor f\cI_t\nor^2_{\widetilde\mL^{q_0}_t(\widetilde\mL^{\bbp_0}_z)},
\end{align*}
which implies by Gronwall's inequality that for any $t\in[0,1]$,
$$
\sup_{s\in[0,t]}\nor u(s)\nor_{\mL^2_z}\lesssim \nor f\cI_t\nor_{\widetilde\mL^{q_0}_t(\widetilde\mL^{\bbp_0}_z)}.
$$
By this estimate and taking supremum in $z_0$ for both sides of \eqref{DD50}, we obtain \eqref{EA80}. In particular, \eqref{EA80} implies the uniqueness of weak solutions.

({\sc Step 2}). In this step we show the existence of weak solutions. Let $\varGamma_\eps$ be as in \eqref{RR}. Define
$$
a_\eps:=a*\varGamma_\eps,\ \ b_\eps:=b*\varGamma_\eps,\ \ f_\eps:=f*\varGamma_\eps.
$$
Under \eqref{C1} and {\bf ($\widetilde{\bf H}^{q_1}_b$)}, we have
\begin{align}\label{DS3}
\kappa_0\mI\leq a_\eps\leq \kappa_1\mI,\ \ \sup_{\eps\in(0,1)}\nor b_\eps\nor_{\widetilde\mL^{q_1}_t(\widetilde\mL^{\bbp_1}_z)}\leq C(\kappa_2),
\end{align}
and for each $\eps\in(0,1)$,
$$
a_\eps, b_\eps, f_\eps\in C^\infty_b(\mR^{1+2d}).
$$
It is well known that there is a unique smooth solution $u_\eps$ solving the following PDE:
\begin{align}\label{PDE10}
\p_t u_\eps=\div_v(a_\eps \cdot\nabla_v u_\eps)+v\cdot\nabla_x u_\eps+b_\eps\cdot \nabla_v u_\eps+f_\eps,\ \ u_\eps(0)=0.
\end{align}
In particular, by Remark \ref{Re32}, for any $\varphi\in C^\infty_c(\mR^{1+2d})$,
\begin{align}\label{PDE11}
\begin{split}
\int_{\mR^{1+2d}}u_\eps\p_t\varphi&=-\int_{\mR^{1+2d}}\<a_\eps\cdot\nabla_v u_\eps,\nabla_v\varphi\>
-\int_{\mR^{1+2d}}u_\eps(v\cdot\nabla_x \varphi)\\
&\quad+\int_{\mR^{1+2d}}(b_\eps\cdot\nabla_v u_\eps)\varphi+\int_{\mR^{1+2d}}f_\eps\varphi.
\end{split}
\end{align}
Moreover, by \eqref{DS3} and \eqref{EA80}, we also have
\begin{align}\label{SD1}
\sup_\eps(\|u_\eps\cI_1 \|_{\mL^\infty}+\nor u_\eps\cI_1\nor_{\widetilde\sV})<\infty.
\end{align}
By the weak compactness of $\widetilde\sV$, there are a sequence $\eps_k\to 0$ and $u\in\widetilde\sV\cap \mL^\infty$ such that
for each $\varphi\in C^\infty_c(\mR^{1+2d})$ and $\Phi\in L^2(\mR^{1+2d};\mR^d)$,
$$
\int_{\mR^{1+2d}}u_{\eps_k}\varphi\to\int_{\mR^{1+2d}}u\varphi,\ \ \int_{\mR^{1+2d}}(\varphi\Phi)\cdot\nabla_v u_{\eps_k}\to\int_{\mR^{1+2d}}(\varphi\Phi)\cdot\nabla_vu.
$$
Taking limits for both sides of \eqref{PDE11} along $\eps_k$, one sees that $u$ is a weak solution of PDE \eqref{PDE0}. 
Let us only show the following limit since the others are completely the same: for $\varphi\in C^\infty_c(\mR^{1+2d})$,
$$
\lim_{k\to\infty}\int_{\mR^{1+2d}}(b_{\eps_k}\cdot\nabla_v u_{\eps_k})\varphi=\int_{\mR^{1+2d}}(b\cdot\nabla_v u)\varphi.
$$
Let $Q$ be the support of $\varphi$. By H\"older's inequality and \eqref{SD1}, we have
\begin{align*}
\left|\int_{\mR^{1+2d}}((b_{\eps_k}-b)\cdot\nabla_v u_{\eps_k})\varphi\right|&\leq \|(b_{\eps_k}-b)\1_Q\|_{\mL^2}
\|\nabla_v u_{\eps_k}\1_Q\|_{\mL^2}\|\varphi\|_{\mL^\infty_z}\\
&\leq C\|(b_{\eps_k}-b)\1_Q\|_{\mL^2}\to 0,\ k\to\infty,
\end{align*}
and by the weak convergence of $\nabla_v u_{\eps_k}$,
$$
\lim_{k\to\infty}\int_{\mR^{1+2d}}b\cdot\nabla_v (u_{\eps_k}-u)\,\varphi=0.
$$
Thus we complete the proof.
\end{proof}

Now we present a stability result which shall be used to show the well-posedness of generalized martingale problem when the coefficients are bounded measurable.
\bt\label{Th42}
Let $(a_n, b_n,f_n)$ be a sequence of coefficients that satisfy the following assumptions:
\begin{enumerate}[(i)]
\item $a_n$ satisfies \eqref{C1} uniformly in $n$ and $b_n$, $f_n$ are uniformly bounded.
\item  $(a_n,b_n,f_n)$ converges to $(a,b,f)$ in Lebesgue measure as $n\to\infty$. 
\end{enumerate}
Let $u_n$ and $u$ be the respective weak solutions of PDE \eqref{PDE0} corresponding to $(a_n,b_n,f_n)$ and $(a,b,f)$ with zero initial value.
Then for any bounded domain $Q\subset\mR_+\times\mR^{2d}$,
\begin{align}\label{CC41}
\lim_{n\to\infty}\sup_{(t,z)\in Q}|u_{n}(t,z)-u(t,z)|=0.
\end{align}
\et
\begin{proof}
First of all, by (i) and \eqref{EA80} we have for any $T>0$,
\begin{align}\label{CC4}
\sup_n(\|u_n\cI_T \|_{\mL^\infty}+\nor u_n\cI_T\nor_{\widetilde\sV})<\infty.
\end{align}
By Remark \ref{Re32}, one can extend $u_n$ to a global weak solution in $\mR^{1+2d}$.
Thus by (i) and \cite[Theorem 3]{GIMV19}, there are $\alpha\in(0,1)$ and constant $C>0$ such that for any $n\in\mN$ and all $(t,z), (t',z')\in Q_{1/2}$,
$$
|u_n(t,z)-u_n(t',z')|\leq C|(t,z)-(t',z')|^\alpha(\|\1_{Q_1}u_n\|_{\mL^2}+\|\1_{Q_1}f_n\|_{\mL^\infty}),
$$
which, together with \eqref{CC4} and a shifting argument as used in {\sc Step 1} of Theorem \ref{Th41}, yields that for any bounded domain $Q\subset\mR^{1+2d}$, 
there is a constant $C>0$ such that for all $n\in\mN$ and $(t,z), (t',z')\in Q$,
$$
|u_n(t,z)-u_n(t',z')|\leq C|(t,z)-(t',z')|^\alpha.
$$
In particular, by \eqref{CC4}, the weak compactness of $\widetilde\sV$ and Ascolli-Arzela's theorem, 
there are a subsequence $n_k$ and $\bar u\in\widetilde\sV\cap C(\mR^{1+2d})$ such that for any $\Phi\in L^2(\mR^{1+2d};\mR^d)$ with compact support,
$$
\int_{\mR^{1+2d}}\Phi\cdot\nabla_v u_{n_k}\to\int_{\mR^{1+2d}}\Phi\cdot\nabla_v\bar u,
$$
and for any bounded domain $Q\subset\mR^{1+2d}$,
$$
\lim_{k\to\infty}\sup_{(t,z)\in Q}|u_{n_k}(t,z)-\bar u(t,z)|=0.
$$
As in {\sc Step 2} of Theorem \ref{Th41}, by (ii) one sees that $\bar u=u$ is the unique weak solution of PDE \eqref{PDE0}.
By a contradiction method, one sees that \eqref{CC41} holds for the whole sequence.
\end{proof}

Below we assume that $a\in L^\infty_T(C^\infty_b(\mR^{2d}))$ satisfies \eqref{C1} and
\begin{align}\label{KB3}
b_1,b_2\in L^\infty_T(C^\infty_b(\mR^{2d})),\ g\in L^\infty_T(C^\infty_c(\mR^{2d})),\ \varphi\in C^\infty_b(\mR^{2d}).
\end{align}
For applications in SDEs with density dependent coefficients, we need the following  apriori $L^\infty$-estimate for the Cauchy problem of the following kinetic FPKE:
$$
\p_t u=\div_v(a\cdot\nabla_v u)+v\cdot\nabla_x u+b_1\cdot \nabla_v u+\div_v(b_2 u)+g,\ \ u(0)=\varphi.
$$
Note that under \eqref{KB3}, it is well known that there is a unique smooth solution to the above PDE.
\bt\label{Th44}
Suppose that for some $q_i\in(2,4)$ and $\bbp_i\in(2,\infty)^{2d}$ with $\tfrac2{q_i}+\bba\cdot\tfrac{1}{\bbp_i}<1$, $i=1,2$,
$$
\nor b_1\nor_{\widetilde\mL^{q_1}_t(\widetilde\mL^{\bbp_1}_z)}
+\nor b_2\nor_{\widetilde\mL^{q_2}_t(\widetilde\mL^{\bbp_2}_z)}\leq\kappa_2.
$$
If $(q_0,\bbp_0)\in\mR^{1+2d}$ satisfies \eqref{Con6} for some $\alpha_0\in[0,1)$, then for any $T>0$,
there is a constant $C>0$ only depending on $T, \alpha_0,\kappa_i, q_i,\bbp_i$ such that for all $t\in[0,T]$,
\begin{align}\label{EA81}
\|u\cI_t\|_{\mL^\infty}\lesssim_C\nor g\cI_t\nor_{\widetilde\mL^{q_0}_t(\widetilde\bB^{-\alpha_0}_{\bbp_0;\bba})}+\|\varphi\|_{C^1_b}.
\end{align}
\et
\begin{proof}
First of all, let $\bar u:=P_t\varphi$, where $P_t$ is defined by \eqref{AA13}. Note that
$$
\p_t \bar u=\Delta_v\bar u+v\cdot\nabla_x\bar u,\ \ \bar u(0)=\varphi.
$$
It is easy to see that $\tilde u:=u-\bar u$ solves the following PDE:
\begin{align}\label{PDE98}
\p_t \tilde u=\div_v(a\cdot\nabla_v\tilde u)+v\cdot\nabla_x \tilde u+b_1\cdot \nabla_v\tilde u+\div_v G+\tilde g+g,\ \ \tilde u(0)=0,
\end{align}
where
$$
G:=(a-\mI)\cdot\nabla_v\bar u+b_2  u,\ \ \tilde g:=b_1\cdot\nabla_v\bar u.
$$
Below we consider PDE \eqref{PDE98} with zero initial value. 
Since $q_1\in(2,4)$ and $\bbp_1\in(2,\infty)^{2d}$ satisfy $\tfrac2{q_1}+\bba\cdot\tfrac{1}{\bbp_1}<1$, one can choose 
$q_3\in(1,\frac43)$ and $\bbp_3\in(1,\infty)^{2d}$ so that
$$
\bbp_3\leq\bbp_1,\ \ \tfrac2{ q_3}+\bba\cdot\tfrac1{\bbp_3}<2.
$$
Suppose $T=1$. As in the {\sc Step 1} of Theorem \ref{Th41}, by applying Theorem \ref{Th34} to $\tilde u^+_{z_0}$ and  $\tilde u^-_{z_0}$ (see Remark \ref{Re45}), 
there is a constant $C=C(\kappa_i,q_i,\bbp_i)>0$ such that for each $t\in[0,1]$ and $z_0\in\mR^{2d}$,
\begin{align}\label{DD55}
\begin{split}
\|\1_{Q_1}\tilde u_{z_0}\cI_t\|_{\mL^\infty}+\|\1_{Q_1}\nabla_v \tilde u_{z_0}\cI_t\|_{\mL^2}
&\lesssim_C\|\1_{Q_2}\tilde u_{z_0}\cI_t\|_{\mL^2}+\cG^t_{z_0},
\end{split}
\end{align}
where 
$$
\cG^t_{z_0}:=\|g_{z_0}\chi_2\cI_t\|_{\mL^{q_0}_t(\bB^{-\alpha_0}_{\bbp_0;\bba})}
+\|\tilde g_{z_0}\chi_2\cI_t\|_{\mL^{ q_3}_t(\mL^{\bbp_3}_z)}+\|G_{z_0}\chi_2\cI_t\|_{\mL^{ q_2}_t(\mL^{\bbp_2}_z)}.
$$
In particular, for each $t\in[0,1]$,
\begin{align*}
\|\1_{Q_1}\tilde  u_{z_0}(t)\|^2_{\mL^2_z}&\leq {\rm Vol}(Q_1)\|\1_{Q_1}\tilde u_{z_0}(t)\|_{\mL^\infty_z}^2
\lesssim\int^t_0\|\1_{Q_2}\tilde u_{z_0}(s)\|_{\mL^2_z}^2\dif s+(\cG^t_{z_0})^2.
\end{align*}
Taking supremum in $z_0\in\mR^{2d}$ and by \eqref{DS1}, we obtain
\begin{align*}
\nor \tilde u(t)\nor^2_{\mL^2_z}\lesssim\int^t_0\nor\tilde  u(s)\nor_{\mL^2_z}^2\dif s+\sup_{z_0}(\cG^t_{z_0})^2,
\end{align*}
which implies by Gronwall's inequality that for any $t\in[0,1]$,
$$
\sup_{s\in[0,t]}\nor\tilde  u(s)\nor_{\mL^2_z}\lesssim \sup_{z_0}\cG^t_{z_0}.
$$
Substituting it into \eqref{DD55} and taking supremum in $z_0\in\mR^{2d}$, we get
\begin{align}\label{DD56}
\|u\cI_t\|_{\mL^\infty}\leq\sup_{s\in[0,t]}\|P_s\varphi\|_\infty+\|\tilde u\cI_t\|_{\mL^\infty}\lesssim_C\|\varphi\|_\infty+\sup_{z_0}\cG^t_{z_0}.
\end{align}
By $\|\nabla_v\bar u\|_{\mL^\infty}=\sup_t\|\nabla_v P_t \varphi\|_{\mL^\infty_z}\lesssim \|\varphi\|_{C^1_b}$, and $\bbp_3\leq\bbp_1$ and $q_3\leq q_1$,
$$
\sup_{z_0}\|\tilde g_{z_0}\chi_2\cI_t\|_{\mL^{ q_3}_t(\mL^{\bbp_3}_z)}\leq\nor b_1\nor_{\widetilde\mL^{ q_3}_t(\widetilde\mL^{\bbp_3}_z)}\|\nabla_v\bar u\|_{\mL^\infty}
\lesssim\nor b_1\nor_{\widetilde\mL^{ q_1}_t(\widetilde\mL^{\bbp_1}_z)}\|\varphi\|_{C^1_b},
$$
and
\begin{align*}
\sup_{z_0}\|G_{z_0}\chi_2\cI_t\|_{\mL^{ q_2}_t(\mL^{\bbp_2}_z)}&\lesssim \|\nabla_v\bar u\|_{\mL^\infty}+\sup_{z_0}\|(b_2  u)_{z_0}\chi_2\cI_t\|_{\mL^{ q_2}_t(\mL^{\bbp_2}_z)}\\
&\lesssim \|\varphi\|_{C^1_b}+\left(\int^t_0\nor b_2(s)\nor^{q_2}_{\mL^{\bbp_2}_z} \| u(s)\|^{q_2}_{\mL^\infty_z}\dif s\right)^{1/q_2}.
\end{align*}
Substituting the above estimates into \eqref{DD56}, we get for any $t\in[0,1]$,
\begin{align*}
\|u(t)\|_{\mL^\infty_z}\lesssim_C\|\varphi\|_{C^1_b}+\nor g\cI_t\nor_{\widetilde\mL^{q_0}_t(\widetilde\bB^{-\alpha_0}_{\bbp_0;\bba})}
+\left(\int^t_0\nor b_2(s)\nor^{q_2}_{\mL^{\bbp_2}_z} \| u(s)\|^{q_2}_{\mL^\infty_z}\dif s\right)^{1/q_2},
\end{align*}
which implies the desired estimate \eqref{EA81} by Gronwall's inequality.
\end{proof}
\section{Weak solutions of second order McKean-Vlasov SDEs}\label{sec5}

In this section we devote to showing the existence of weak solutions to the following second order McKean-Vlasov SDE:
\begin{align}\label{SDE00}
\dif X_t=V_t\dif t,\ \ \dif V_t=b_Z(t,Z_t)\dif t+\sigma_Z(t,X_t)\dif W_t,
\end{align}
where $Z_t=(X_t, V_t)$ and $b_Z, \sigma_Z$ are defined by \eqref{BB8} and \eqref{BB9}, respectively.

Let $\mC$ be the space of all continuous functions from $\mR_+$ to $\mR^{2d}$, which is endowed with 
the locally uniformly convergence topology. Let $\omega_t$ be the canonical process on $\mC$, and $\cB_t$ the natural filtration associated with $\omega_t$.
Let $\cP(\mC)$ be the space of all probability measures over $\mC$. 

The notion of weak solutions has been introduced in Definition \ref{Def2}.
Now we introduce the equivalent notion of classical martingale solutions (cf. \cite{St-Va}).
The advantage of using martingale solutions is that it avoids the use of stochastic integrals.

\bd\label{Def3}
Let $\nu\in\cP(\mR^{2d})$. We call a probability measure $\mP\in\cP(\mC)$ a classical martingale solution of DDSDE \eqref{SDE00} with initial distribution $\nu$ if 
\begin{enumerate}[(i)]
\item $\mP\circ \omega^{-1}_0=\nu$, and for Lebesgue almost all $t\geq 0$, $\mP\circ \omega_t^{-1}(\dif z)=\rho_t(z)\dif z$.

\item For each $f\in C^2_c(\mR^{2d})$, the process
\begin{align}\label{MM1}
M_t(\omega):=f(\omega_t)-f(\omega_0)-\int^t_0\left(\tr(a_{\rho}\cdot\nabla^2_vf)+v\cdot\nabla_xf+b_{\rho}\cdot\nabla_vf\right)(r,\omega_r)\dif r
\end{align}
is a $\cB_t$-martingale with respect to $\mP$, where for $\<\rho_t\>(x):=\int_{\mR^d}\rho_t(x,v)\dif v$,
\begin{align}
a_{\rho}(t,x)&:=\int_{\mR^d}a(t,x, \<\rho_t\>(x), z')\rho_t(z')\dif z',\label{MM6}\\
b_\rho(t,z)&:=\int_{\mR^{2d}}b(t,z, \rho_t(z), z')\rho_t(z')\dif z'.\label{MM66}
\end{align}
\end{enumerate}
The set of all the martingale solutions $\mP\in\cP(\mC)$  with initial distribution $\nu$ is denoted by $\sM^{a,b}_{\nu}$.
\ed

\br\label{Pr53}\rm
Note that $a_\rho$ and $b_\rho$ are only well defined except for some Lebesgue null set. 
This does not affect the definition of \eqref{MM1} due to the point (i).
Under \eqref{ELL}, it is well known that weak solution in Definition \ref{Def2} is equivalent to the above martingale solution (see \cite{St-Va}).
\er

 Our main result of this section is to show the following existence result.
 \bt\label{Th63}
 Under {\bf (H$_1$)} and {\bf (H$_2$)}, for any initial probability measure $\nu\in\cP(\mR^{2d})$, there exists at least one classical martingale solution $\mP\in\sM^{a,b}_\nu$
 in the sense of Definition \ref{Def3}. Moreover, $\rho$ enjoys the regularity \eqref{DL131}.
 \et
 
To show the existence of a martingale solution, by mollifying the coefficients and the fact that $a_\rho$ does not depend on $v$, 
the key point is to use \eqref{EA80} and \eqref{EA81} to establish the uniform Krylov estimate for approximating DDSDEs.
For $N\in\mN$, let $\varGamma\in C^\infty_c(\mR^N)$ be a probability density function and define 
$$
\varGamma_n(\cdot)=n^N\varGamma(n\cdot),\ \  n\in\mN.
$$ 
With a little confusion, we shall use the same notation $\varGamma_n$ to denote different mollifiers for different dimension $N$. 
For $n\in\mN$,  define
\begin{align*}
b_n(t,z,r,z')&:=(b(t,\cdot,\cdot,\cdot)*\varGamma_n)(z,r,z'),\\ 
a_n(t,x,r,z')&:=(a(t,\cdot,\cdot,\cdot)*\varGamma_n)(x,r,z'),
\end{align*}
and for a probability measure $\mu\in\cP(\mR^{2d})$,
\begin{align*}
b^n_\mu(t,z)&:=\int_{\mR^{2d}} b_n(t,z,\mu*\varGamma_n(z), z')\mu(\dif z'),\\  
a^n_\mu(t,x)&:=\int_{\mR^d}a_n(t,x,\mu*\varGamma_n(x),z')\mu(\dif z').
\end{align*}
For two probability measures $\mu_1,\mu_2\in\cP_2(\mR^{2d})$
with finite second order moments, the Wasserstein metric is defined by
$$
\cW_2(\mu_1,\mu_2):=\inf\left\{\int_{\mR^{4d}}|z-z'|^2\pi(\dif z,\dif z'): \pi\in\cP(\mR^{4d}),\ \pi(\cdot,\mR^{2d})=\mu_1,\pi(\mR^{2d},\cdot)=\mu_2\right\}^{\frac12}.
$$

The following lemma provides necessary estimates for $b^n_\mu$ and $a^n_\mu$. 
\bl\label{Le71}
Under {\bf (H$_1$)} and {\bf (H$_2$)},  for each $n\in\mN$, $T>0$ and $\mu\in\cP(\mR^{2d})$, we have
\begin{align}\label{DR121}
b^n_\mu\in L^\infty_TC^\infty_b(\mR^{2d}),\  \ a^n_\mu\in L^\infty_TC^\infty_b(\mR^{d}),
\end{align}
and  there is a constant $C=C(T,n)>0$ such that for any $\mu_1,\mu_2\in\cP_2(\mR^{2d})$,
\begin{align}\label{DR12}
\|b^n_{\mu_1}-b^n_{\mu_2}\|_\infty+\|\sqrt{a^n_{\mu_1}}-\sqrt{a^n_{\mu_2}}\|_\infty\lesssim_C \cW_2(\mu_1,\mu_2).
\end{align}
Moreover, we also have
\begin{align}\label{DR2}
\sup_n\nor b^n_\mu\nor_{\widetilde\mL^{q_1}_t(\widetilde\mL^{\bbp_1}_z)}\leq \nor h\nor_{\widetilde\mL^{q_1}_t(\widetilde\mL^{\bbp_1}_z)}\leq\kappa_2,
\end{align}
and for all $\xi\in\mR^d$ and $(t,x)\in\mR_+\times\mR^{d}$,
\begin{align}\label{ELL2}
\kappa_0|\xi|^2\leq\xi\cdot a^n_\mu(t,x)\xi\leq\kappa_1|\xi|^2.
\end{align}
\el
\begin{proof}
Estimates \eqref{DR121} and \eqref{ELL2} are obvious by definition. For \eqref{DR2}, noting that by \eqref{DR1},
\begin{align*}
|b^n_\mu(t,z)|&\leq \int_{\mR^{2d}}\left( \int_{\mR^{4d}}h(t,z-z_1-z'+z'_1)\varGamma_n(z_1,z_1')\dif z_1\dif z_1'\right)\mu(\dif z'),
\end{align*}
by Minkowskii's inequality, we have
\begin{align*}
\nor b^n_\mu\nor_{\widetilde\mL^{q_1}_t(\widetilde\mL^{\bbp_1}_z)}&\leq \int_{\mR^{2d}}\left( \int_{\mR^{4d}}
\nor h\nor_{\widetilde\mL^{q_1}_t(\widetilde\mL^{\bbp_1}_z)}
\varGamma_n(z_1,z_1')\dif z_1\dif z_1'\right)\mu(\dif z')=\nor h\nor_{\widetilde\mL^{q_1}_t(\widetilde\mL^{\bbp_1}_z)}.
\end{align*}

Below we  prove the second estimate in \eqref{DR12}, and divide the proofs into three steps.

({\sc Step 1}) In this step we show the following claim: Suppose that $A(r):\mR\to \mM^d_{\rm sym}$ is a $C^1$-positive 
definite matrix-valued function so that $A(r)\geq \kappa_0\mI$ for all $r$. Then  (cf. \cite[Theorem 5.2.2]{St-Va})
\begin{align}\label{XC1}
\|\p_r A^{\frac12}(r)\|_{HS}\leq \|\p_r A(r)\|_{HS}/(2\sqrt{\kappa_0}).
\end{align}
Indeed, by the chain rule, we have
$$
\p_r A(r)=\p_r A^{\frac12}(r)\cdot A^{\frac12}(r)+A^{\frac12}(r)\cdot\p_r A^{\frac12}(r).
$$
Let $P$ be an orthogonal matrix so that $P A^{\frac12}(r) P^*=:D$ is an diagonal matrix, where the asterisk stands for the transpose of a matrix.
Then
$$
P\p_r A(r)P^*=(P\p_r A^{\frac12}(r)P^*) D+D(P\p_r A^{\frac12}(r)P^*)
$$
and
$$
(P\p_r A^{\frac12}(r)P^*)_{ij}=(P\p_r A(r)P^*)_{ij}/(D_{jj}+D_{ii}).
$$
From this we derive \eqref{XC1}.

({\sc Step 2}) Let $\mH$ be a separable Hilbert space and $A: \mH\to \mM^d_{\rm sym}$  a $C^1$-positive 
definite matrix-valued function. Suppose that for all $h,h'\in\mH$,
$$
A(h)\geq\kappa_0\mI,\ \|A(h)-A(h')\|_{HS}\leq\kappa_1\|h-h'\|_\mH.
$$
 Then it holds that
\begin{align}\label{XC11}
\|A^{\frac12}(h)-A^{\frac12}(h')\|_{HS}\leq \kappa_1\|h-h'\|_\mH/(2\sqrt{\kappa_0}).
\end{align}
Indeed, without loss of generality, we may assume $\mH=\ell^2$, where $\ell^2$ is the usual sequence Hilbert space.
For each $h=(h_1,h_2,\cdots)\in\ell^2$, the assumption implies that
$$
\sum_{j\in\mN}\|\p_{h_j}A(h)\|^2_{HS}\leq\kappa_1^2.
$$
Thus, by H\"older's inequality and \eqref{XC1}, we have for some ${\bf h}^*_j\in\ell^2$,
\begin{align*}
\|A^{\frac12}(h)-A^{\frac12}(h')\|_{HS}&\leq\sum_{j\in\mN}|h_j-h_j'|\|(\p_{h_j}A^{\frac12})({\bf h}^*_j)\|_{HS}\\
&\leq\|h-h'\|_\mH\left(\sum_{j\in\mN}\|(\p_{h_j}A^{\frac12})({\bf h}^*_j)\|^2_{HS}\right)^{1/2}\\
&\leq\frac{\|h-h'\|_\mH}{2\sqrt{\kappa_0}}\left(\sum_{j\in\mN}\|(\p_{h_j}A)({\bf h}^*_j)\|^2_{HS}\right)^{1/2}\leq\frac{\kappa_1\|h-h'\|_\mH}{2\sqrt{\kappa_0}}.
\end{align*}

({\sc Step 3}) For simplicity, we drop the variable $(t,x)$. Suppose that $X\in L^2(\Omega,\bP)$ has distribution $\mu$.
Then we can write
$$
a^n_\mu=\int_{\mR^{2d}}a_n(\mu*\varGamma_n,z')\mu(\dif z')=\bE a_n(\mu*\varGamma_n,X)=:A(\mu*\varGamma_n,X).
$$
Thus for any $X, Y\in L^2(\Omega,\bP)$ with distributions $\mu_1$ and $\mu_2$, we have
\begin{align*}
\sqrt{a^n_{\mu_1}}-\sqrt{a^n_{\mu_2}}=\sqrt{A(\mu_1*\varGamma_n,X)}-\sqrt{A(\mu_2*\varGamma_n,Y)},
\end{align*}
and by \eqref{XC1} and \eqref{XC11}, 
\begin{align*}
\|\sqrt{a^n_{\mu_1}}-\sqrt{a^n_{\mu_2}}\|_{HS}&\leq C_n(|\mu_1*\varGamma_n-\mu_2*\varGamma_n|+\|X-Y\|_{L^2})\\
&\leq C_n(\bE|\varGamma_n(\cdot-X)-\varGamma_n(\cdot-Y)|+\|X-Y\|_{L^2})\\
&\leq C_n\|X-Y\|_{L^2},
\end{align*}
which in turn implies that
$$
\|\sqrt{a^n_{\mu_1}}-\sqrt{a^n_{\mu_2}}\|_{HS}\leq C_n\cW_2(\mu_1,\mu_2).
$$
The proof is complete.
\end{proof}

Now let $Z_0$ be an $\sF_0$-measurable random variable with distribution $\nu$. Let us consider the following approximating DDSDE:
\begin{align}\label{SDE19}
\dif X^n_t=V^n_t\dif t,\ \ \dif V^n_t=b^n_{\mu_{Z^n_t}}(t,Z^n_t)\dif t+\sqrt{a^n_{\mu_{Z^n_t}}}(t,X^n_t)\dif W_t,
\end{align}
subject to the initial value
$$
Z^n_0:=(-n)\vee Z_0\wedge n,
$$
where for a vector $z\in\mR^{2d}$, $((-n)\vee z\wedge n)_i:=(-n)\vee z_i\wedge n$.

By Lemma \ref{Le71} and standard Picard's iteration, 
it is well known that the above SDE admits a unique strong solution $Z^n_{t}$.
Now we use Theorems \ref{Th41} and \ref{Th44} to derive the following crucial Krylov estimates.
\bl [Krylov's estimates]\label{Kry0}
Let $\alpha_0\in[0,1)$ and $(q_0,\bbp_0)\in(1,\infty)^{1+2d}$ satisfy \eqref{Con6}.
For any $T>0$, there are $\theta\in(0,1)$ and constant $C>0$ depending only on $\kappa_i$ 
and $q_i,\bbp_i$ such that for any $\delta\in(0,1)$, stopping time $\tau\leq T$ and  $f\in C^\infty_c(\mR^{1+2d})$,
\begin{align}\label{Kr0}
\sup_{n\in\mN}\bE\left(\int^{\tau+\delta}_{\tau}f(r,Z^n_r)\dif r\Big|\sF_{\tau}\right)\lesssim_C\delta^\theta
\nor\1_{[0,T]}f\nor_{\widetilde\mL^{q_0}_t(\widetilde\bB^{-\alpha_0}_{\bbp_0;\bba})}
\lesssim_C\delta^\theta\|\1_{[0,T]}f\|_{\mL^{q_0}_t(\bB^{-\alpha_0}_{\bbp_0;\bba})}.
\end{align}
\el
\begin{proof}
Fix $T>0$. By Lemma \ref{Le71}, we have
$$
b^n_{\mu_{Z^n_t}}\in L^\infty_TC^\infty_b(\mR^{2d}),\  \ a^n_{\mu_{Z^n_t}}\in L^\infty_TC^\infty_b(\mR^{d}).
$$
Fix $0\leq t_0<t_1\leq T$ and $f\in C^\infty_c(\mR^{1+2d})$. Let $u_n\in L^\infty_{t_1}(C^\infty_b(\mR^{2d}))$ solve the following backward Kolmogorov equation:
\begin{align}\label{DC10}
\p_t u_n+\tr\big(a^n_{\mu_{Z^n_t}}\cdot\nabla^2_v u_n\big)+v\cdot\nabla_x u_n+b^n_{\mu_{Z^n_t}}\cdot\nabla_v u_n=f,\ \ u_n(t_1)=0.
\end{align}
Since $\bba\cdot\tfrac1{\bbp_0}+\tfrac2{q_0}<2$, one can choose $\bar q_0< q_0$ so that
$$
\bba\cdot\tfrac1{\bbp_0}+\tfrac2{\bar q_0}<2.
$$
Thus by  \eqref{DR2}, \eqref{ELL2}, \eqref{DZ8} and \eqref{EA81}, there is a constant $C=C(\kappa_i, q_i,\bbp_i)>0$ such that
\begin{align}\label{DC0}
\sup_{n\in\mN}\|\1_{[t_0,t_1]}u_n\|_\infty\lesssim_C\nor\1_{[t_0,t_1]}f\nor_{\widetilde\mL^{\bar q_0}_t(\widetilde\bB^{-\alpha_0}_{\bbp_0;\bba})}
\lesssim_C(t_1-t_0)^{\frac1{\bar q_0}-\frac1{ q_0}}\nor\1_{[t_0,t_1]}f\nor_{\widetilde\mL^{q_0}_t(\widetilde\bB^{-\alpha_0}_{\bbp_0;\bba})},
\end{align}
where the second step is due to H\"older's inequality.
Now by \eqref{DC10} and It\^o's formula, we have
\begin{align*}
0=u_n(t_1, Z^n_{t_1})=u_n(t_0, Z^n_{t_0})+\int^{t_1}_{t_0}f(r,Z^n_r)\dif r+\int^{t_1}_{t_0}\Big\<\nabla_vu_n(r,Z^n_r),\sqrt{a^n_{\mu_{Z^n_r}}}(r,X^n_r)\dif W_r\Big\>.
\end{align*}
Hence, by taking conditional expectation with respect to $\sF_{t_0}$ and \eqref{DC0},
\begin{align*}
\bE\left(\int^{t_1}_{t_0} f(r,Z^n_r)\dif r\Big|\sF_{t_0}\right)=-u_n(t_0, Z^n_{t_0})
\lesssim(t_1-t_0)^{\frac1{\bar q_0}-\frac1{ q_0}}\nor\1_{[t_0,t_1]}f\nor_{\widetilde\mL^{q_0}_t(\widetilde\bB^{-\alpha_0}_{\bbp_0;\bba})}.
\end{align*}
By discretization stopping time approximation (see \cite[Remark 1.2]{ZZ18}), we obtain the first estimate in \eqref{Kr0}. 
Note that by Lemma \ref{Le23},
$$
\nor \1_{[0,T]}f\nor_{\widetilde\mL^{q_0}_t(\widetilde\bB^{-\alpha_0}_{\bbp_0;\bba})}\lesssim\sup_{z_0}\|\1_{[0,T]}f(\cdot-z_0)\|_{\mL^{q_0}_t(\bB^{-\alpha_0}_{\bbp_0;\bba})}
\|\chi_2\|_{\mL^{\infty}_t(\bB^1_{\infty;\bba})}\lesssim \|\1_{[0,T]}f\|_{\mL^{q_0}_t(\bB^{-\alpha_0}_{\bbp_0;\bba})}.
$$
The second estimate  in \eqref{Kr0} then follows.
\end{proof}
\br\rm
Under {\bf ($\widetilde {\bf H}_b$)},  one can use Theorem \ref{Th41} to obtain the following Krylov estimate: for any $(q_0,p_0)\in(1,\infty)^{1+2d}$
satisfying \eqref{Con7},
\begin{align}\label{Kr1}
\sup_{n\in\mN}\bE\left(\int^{\tau+\delta}_{\tau}f(r,Z^n_r)\dif r\Big|\sF_{\tau}\right)\leq C\delta^\theta\nor\1_{[0,T]}f\nor_{\mL^{q_0}_t(\mL^{\bbp_0}_z)}.
\end{align}
\er

The following corollary is direct by \eqref{Kr0}.
\bc\label{Cor1}
For Lebesgue almost all $t\geq 0$, $Z^n_{t}$ admits a density $\rho^n_{t}(z)$ with the following regularity: for any $\alpha\in[0,1)$ and $(q,\bbp)\in(1,\infty)^{1+2d}$ satisfying
$$
\tfrac2{q}<1+\alpha,\ \  \ \tfrac2{q}+\bba\cdot(\tfrac1{\bbp}-\1)>2\alpha,
$$
it holds that
\begin{align}\label{DL11}
\sup_{n\in\mN}\|\rho^n\1_{[0,T]}\|_{\mL^{q}_t(\bB^{\alpha}_{\bbp;\bba})}<\infty,\ T>0.
\end{align}
\ec
\begin{proof}
Let $\mu^n_t(\dif z)$ be the probability distribution of $Z^n_t$. By Krylov's estimate \eqref{Kr0},
for any $\alpha_0\in[0,1)$ and $(q_0,\bbp_0)\in(1,\infty)^{1+2d}$ satisfying \eqref{Con6}, there is a constant $C_T>0$ such that for any $f\in C^\infty_c(\mR^{2d})$,
$$
\sup_{n\in\mN}\int^T_0\int_{\mR^{2d}}f(r,z)\mu^n_t(\dif z)\dif r\leq C_T
\|\1_{[0,T]}f\|_{\mL^{q_0}_t(\bB^{-\alpha_0}_{\bbp_0;\bba})}.
$$
This  implies that $\mu^n_t$ has a density $\rho^n_{t}(z)$, which belongs to the dual space of $L^{ q_0}([0,T]; \bB^{-\alpha_0}_{\bbp_0;\bba})$. In particular, we have
 \eqref{DL11}.
\end{proof}

By Krylov's estimate \eqref{Kr0}, we can also show the following tightness result.
\bl\label{Le56}
The law $\mP_n$ of $Z^n$ in $\mC$ is tight. Moreover, for any accumulation point $\mP$, it holds that $\mP\circ\omega_t^{-1}(\dif z)=\rho_t(z)\dif z$,
where $\rho$ satisfies the same estimate  \eqref{DL11}.
\el
\begin{proof}
Fix $T>0$. Let $\tau$ be any stopping time bounded by $T$. By SDE \eqref{SDE19}, we have
\begin{align*}
V^n_{\tau+t}-V^n_{\tau}=\int^{\tau+t}_\tau b^n_{\mu_{Z^n_r}}(r,Z^n_r)\dif r
+\int^{\tau+t}_\tau\sqrt{a^n_{\mu_{Z^n_r}}}(r,X^n_r)\dif W_r,
\end{align*}
and
\begin{align}\label{DA0}
X^n_{\tau+t}-X^n_{\tau}=\int^{\tau+t}_\tau V^n_{r}\dif r.
\end{align}
By Krylov's estimate \eqref{Kr0} and \eqref{DR2}, there is a $\theta\in(0,1)$ such that for all $\delta\in(0,1)$ and $n\in\mN$,
\begin{align*}
\bE\left(\int^{\tau+\delta}_\tau |b^n_{\mu_{Z^n_t}}(r,Z^n_r)|\dif r\right)
\lesssim\delta^\theta\nor b^n_{\mu_{Z^n_{r}}}\nor_{\widetilde\mL^{ q_1}_t(\widetilde\mL^{\bbp_1}_z)}\lesssim\delta^\theta\nor h\nor_{\widetilde\mL^{ q_1}_t(\widetilde\mL^{\bbp_1}_z)},
\end{align*}
and by BDG's inequality and \eqref{ELL2}, 
\begin{align*}
\bE\left|\sup_{t\in[0,\delta]}\int^{\tau+t}_\tau\sqrt{a^n_{\mu_{Z^n_r}}}(r,X^n_r)\dif W_r\right|\lesssim 
\bE\left(\int^{\tau+\delta}_\tau\big\|a^n_{\mu_{Z^n_r}}(r,X^n_r)\big\|_{HS}\dif r\right)^{1/2}\lesssim \delta^{1/2},
\end{align*}
where the implicit constant is independent of $n$, $\tau$ and $\delta\in(0,1)$.
Hence,
\begin{align}\label{BV1}
\sup_n\bE\left(\sup_{t\in[0,\delta]}|V^n_{\tau+t}-V^n_{\tau}|\right)\lesssim \delta^{\theta\wedge(1/2)}.
\end{align}
Noting that for any $R>0$,
$$
\bP\left(\sup_{t\in[0,T+\delta]}|V^n_t|\geq R\right)\leq \bP\left(\sup_{t\in[0,T+\delta]}|V^n_t-V^n_0|\geq \tfrac R2\right)
+\bP\left(|V^n_0|\geq \tfrac R2\right),
$$
and by $|V^n_0|\leq |V_0|$, \eqref{BV1} and Chebyschev's inequality, we get
\begin{align}\label{BV2}
\lim_{R\to\infty}\sup_n\bP\left(\sup_{t\in[0,T+\delta]}|V^n_t|\geq R\right)=0.
\end{align}
Thus, by \eqref{DA0} and \eqref{BV2}, we derive that for any $\eps>0$,
\begin{align*}
\lim_{\delta\to 0}\sup_{n}\bP\left(\sup_{t\in[0,\delta]}|X^n_{\tau+t}-X^n_{\tau}|>\eps\right)=0,
\end{align*}
which, together with \eqref{BV1} and \eqref{BV2}, and by \cite[Theorem 1.3.2]{St-Va}, yields the tightness of $\mP_n$.

Let $\mP$ be any accumulation point. By Prohorov's theorem, there is a subsequence $n_m$ such that  $\mP_{n_m}$ weakly converges to $\mP$ as $m\to\infty$.
By taking weak limits, one sees that for any $T>0$, $f\in C^\infty_c(\mR^{1+2d})$ and  any $\alpha_0\in[0,1)$ and $(q_0,\bbp_0)\in(1,\infty)^{1+2d}$ satisfying \eqref{Con6}, 
$$
\mE^\mP\left(\int^T_0f(r,\omega_r)\dif r\right)\lesssim\|f\|_{\mL^{q_0}_t(\bB^{-\alpha_0}_{\bbp_0;\bba})}.
$$
As in Corollary \ref{Cor1}, $\mP\circ\omega_t^{-1}(\dif z)=\rho_t(z)\dif z$ and $\rho$ satisfies \eqref{DL11}.
\end{proof}

The following lemma provides the strong convergence of the density.
\bl\label{Le57}
For fixed accumulation point $\mP$ of $(\mP_n)_{n\in\mN}$, there is a subsequence $n_k$ such that $\mP_{n_k}$ weakly converges to $\mP$ 
as $k\to\infty$, and for any bounded domain $Q\subset\mR_+\times\mR^{2d}$,
\begin{align}\label{DL1}
\lim_{k\to\infty}\|(\rho^{n_k}-\rho)\1_Q\|_{\mL^{2}_t(\mL^1_z)}=0,
\end{align}
where $\mP^{n_k}\circ\omega_t^{-1}(\dif z)=\rho^{n_k}_t(z)\dif z$ and $\mP\circ\omega_t^{-1}(\dif z)=\rho_t(z)\dif z$.
\el
\begin{proof}
First of all, by Lemma \ref{Le56}, there is a subsequence $n_m$ so that $\mP_{n_m}$ weakly converges to $\mP$ as $m\to\infty$.
Next we look at \eqref{DL1}. Note that $\rho^n_{t}(z)$ solves the following FPKE:
$$
\p_t\rho^n_{t}=\p_{v_i}\p_{v_j}\big((a^n_{\mu_{Z^n_t}})_{ij}\rho^n_{t}\big)-v\cdot\nabla_x\rho^n_{t}+\div_v\big(b^n_{\mu_{Z^n_t}}\rho^n_{t}\big).
$$
Given $R\geq 1$, let $\chi_R$ be the cutoff function defined in \eqref{Cut}. Since $\chi_{2R}\equiv 1$ on the support of $\chi_R$,
by Bernstein's inequality \eqref{Ber}, it is easy to see that
\begin{align}
\|\p_t(\rho^n\chi_R)\|_{\mL^2_t(\bB^{-3}_{1;a})}
&\leq \|\p_{v_i}\p_{v_j}\big(a^n_{\mu_{Z^n_t}}\rho^n\chi_{2R}\big)\chi_R\|_{\mL^2_t(\bB^{-3}_{1;a})}+
\|v\cdot\nabla_x(\rho^n\chi_{2R})\chi_R\|_{\mL^2_t(\bB^{-3}_{1;a})}\no\\
&\quad+\|\div_v\big(b^n_{\mu_{Z^n_t}}\rho^n\chi_{2R}\big)\chi_R\|_{\mL^2_t(\bB^{-3}_{1;a})}\no\\
&\lesssim\|a^n_{\mu_{Z^n_t}}\rho^n\chi_{2R}\|_{\mL^2_t(\mL^1_z)}
+\|\rho^n\chi_{2R}\|_{\mL^2_t(\mL^1_z)}+\|b^n_{\mu_{Z^n_t}}\rho^n\chi_{2R}\|_{\mL^2_t(\mL^1_z)},\label{DQ1}
\end{align}
where the implicit constant may depend on $R$, but is independent of $n$.
For the last term, by H\"older's inequality, we have
\begin{align}\label{DQ2}
\|b^n_{\mu_{Z^n_t}}\rho^n\chi_{2R}\|_{\mL^2_t(\mL^1_z)}
\leq\|b^n_{\mu_{Z^n_t}}\1_{Q_{4R}}\|_{\mL^{ q_1}_t(\mL^{\bbp_1}_z)}\|\rho^n\1_{Q_{4R}}\|_{\mL^{\bar q_1}_t(\mL^{\bar\bbp_1}_z)},
\end{align}
where $\frac1{\bar  q_1}+\frac1{ q_1}=\frac12$ and $\frac1{\bar\bbp_1}+\frac1{\bbp_1}=\1$.
Since $\frac2{ q_1}+\bba\cdot\frac1{\bbp_1}<1$, we have $\frac2{\bar q_1}+\bba\cdot(\frac1{\bar\bbp_1}-\1)>0$.
Thus by \eqref{DL11}, \eqref{DR2}, \eqref{DQ1} and \eqref{DQ2}, we obtain
$$
\sup_n\|\p_t(\rho^n\chi_R)\|_{\mL^2_t(\bB^{-3}_{1;a})}<\infty.
$$
On the other hand, by \eqref{DL11} we also have for some $\alpha>0$,
$$
\sup_n\|\rho^n\chi_R\|_{\mL^2_t(\bB^\alpha_{1;a})}<\infty.
$$
Since $\bB^\alpha_{1;a}$ is locally and compactly embedded in $\mL^1_z$,  by Aubin-Lions' lemma (see \cite{Si}) and a further diagolization method, there is a subsubsequence $n_k$ of $n_m$
so that $\mP_{n_k}$ weakly converges to $\mP$ as $k\to\infty$ and \eqref{DL1} holds. The proof is complete.
\end{proof}

{\bf Convention:} For simplicity of notations, we shall assume that $(\mP_n)_{n\in\mN}$ weakly converges to $\mP$ and \eqref{DL1} holds for the whole sequence below.

\medskip

We need the following technical lemma for taking limits below.
\bl\label{Le511}
Fix $N\in\mN$. Let $f\in L^1_{loc}(\mR^{N+1})$ and $\rho_n,\rho\in L^1_{loc}(\mR^N)$. Suppose that for any $m\in\mN$ and bounded $Q\subset \mR^N$,
$$
\lim_{h\to 0}\int_Q\sup_{|r|,|r'|\leq m,|r-r'|\leq h} |f(z, r)-f(z,r')|\dif z=0,\ \ \lim_{n\to\infty}\|\rho_n-\rho\|_{L^1(Q)}=0.
$$
Let $f_n(z,r):=(f*\varGamma_n)(z,r)$. Then for any $\eps>0$ and bounded $Q\subset\mR^N$,
$$
\lim_{n\to\infty}\cL\Big\{z\in Q: |f_n(z,\rho_n(z))-f(z,\rho(z))|\geq\eps\Big\}=0,
$$
where $\cL$ stands for the Lebesgue measure in $\mR^N$. 
\el
\begin{proof}
Fix bounded $Q\subset\mR^N$. For any $m,n\in\mN$ and $h>0$, we define
$$
A^m_n:=\big\{z\in Q: |\rho_n(z)|\vee|\rho(z)|\leq m\big \},\ \ D^h_n:=\big\{z\in Q: |\rho_n(z)-\rho(z)|\leq h\big \},
$$
and for $\eps>0$,
$$
Q^\eps_n:=\big\{z\in Q: |f_n(z,\rho_n(z))-f(z,\rho(z))|\geq\eps\big\}.
$$
By the assumption and Chebyschev's inequality, we clearly have
\begin{align}\label{Lim112}
\lim_{m\to\infty}\sup_n\cL((A^m_n)^c)=0,\ \ \lim_{n\to\infty}\cL((D^h_n)^c)=0,\ \ \forall h>0.
\end{align}
Note that
\begin{align*}
Q^\eps_n&\subset\big\{z\in Q: |f(z,\rho_n(z))-f(z,\rho(z))|\geq\eps/2\big\}\\
&\quad\cup\big\{z\in Q: |f_n(z,\rho_n(z))-f(z,\rho_n(z))|\geq\eps/2\big\}=:Q^{\eps,1}_n\cup Q^{\eps,2}_n.
\end{align*}
Thus, to show $\lim_{n\to\infty}\cL(Q^\eps_n)=0$, by \eqref{Lim112}, it suffices to show that for each $m\in\mN$,
\begin{align}\label{Lim12}
\lim_{h\to 0}\sup_n\cL(Q^{\eps,1}_n\cap A^m_n\cap D^h_n)=0,
\end{align}
and
\begin{align}\label{Lim11}
\lim_{n\to\infty}\cL(Q^{\eps,2}_n\cap A^m_n)=0.
\end{align}
For limit \eqref{Lim12}, it follows by the following observation and the assumption
$$
Q^{\eps,1}_n\cap A^m_n\cap D^h_n\subset\Big\{z\in Q: \sup_{|r|,|r'|\leq m, |r-r'|\leq h}|f(z,r)-f(z,r')|\geq\eps/2\Big\}.
$$
For limit \eqref{Lim11}, note that by the definition of $f_n$ and the assumption,
\begin{align}\label{Lim0}
\begin{split}
&\lim_{h\to 0}\sup_{n\in\mN}\int_{Q}\sup_{|r|\leq m}|f_n(z,r+h)-f_n(z,r)|\dif z\\
&\quad\leq \lim_{h\to 0}\int_{Q'}\sup_{|r|\leq m+1}|f(z,r+h)-f(z,r)|\dif z=0,
\end{split}
\end{align}
where $Q\Subset Q'$.
Since for each $r$,
$$
\lim_{n\to\infty}\int_Q|f_n(z,r)-f(z,r)|\dif z=0,
$$
by \eqref{Lim0} and a finitely covering technique, we have for each $m\in\mN$,
$$
\lim_{n\to\infty}\int_Q\sup_{|r|\leq m}|f_n(z,r)-f(z,r)|\dif z=0,
$$
which in turn implies \eqref{Lim11} by noting that
$$
Q^{\eps,2}_n\cap A^m_n\subset\Big\{z\in Q: \sup_{|r|\leq m}|f_n(z,r)-f(z,r)|\geq\eps/2\Big\}.
$$
The proof is complete.
\end{proof}

Now define
\begin{align}\label{DS88}
\bar a_n(t,x):=\int_{\mR^{2d}}a_n\Big(t,x,(\<\rho^{n}_{t}\>*\varGamma_n)(x),z'\Big)\rho^{n}_{t}(z')\dif z',
\end{align}
where $\<\rho^{n}_{t}\>(x):=\int_{\mR^d}\rho^{n}_{t}(x,v)\dif v$ and
\begin{align}\label{DS89}
\bar b_n(t,z):=\int_{\mR^{2d}}b_n\Big(t,z,(\rho^{n}_t*\varGamma_n)(z),z'\Big)\rho^{n}_t(z')\dif z'.
\end{align}
We have
\bl\label{Le58}
For any bounded domain $Q\subset\mR_+\times\mR^{2d}$, it holds that for any $( q_0,\bbp_0)\in[1,\infty)^{1+2d}$,
\begin{align}\label{DS8}
\lim_{n\to\infty}\|(\bar a_n-a_\rho)\1_Q\|_{\mL^{ q_0}_t(\mL^{\bbp_0}_z)}=0,
\end{align}
and for any $q_0<q_1$ and $\bbp_0<\bbp_1$, where $(q_1,\bbp_1)$ is from \eqref{DR1},
\begin{align}\label{DS9}
\lim_{n\to\infty}\|(\bar b_n-b_\rho)\1_Q\|_{\mL^{ q_0}_t(\mL^{\bbp_0}_z)}=0,
\end{align}
where $a_\rho$ and $b_\rho$ are defined by \eqref{MM6} and \eqref{MM66} in terms of $\rho$ in Lemma \ref{Le57}.
\el
\begin{proof}
We only prove \eqref{DS9} since \eqref{DS8} is completely the same. Note that by \eqref{DR1} and \eqref{HD9}, 
$$
\sup_n\|(\bar b_n-b_\rho)\1_Q\|_{\mL^{ q_1}_t(\mL^{\bbp_1}_z)}\lesssim\nor h\nor_{\mL^{ q_1}_t(\mL^{\bbp_1}_z)}.
$$
It suffices to prove that
$$
\lim_{n\to\infty}\|(\bar b_n-b_\rho)\1_Q\|_{\mL^1}=0.
$$
For $R\geq 1$, let $B_R:=\{z\in\mR^{2d}: |z|\leq R\}$.
By \eqref{DS89} and \eqref{MM66}, we make the following decomposition:
\begin{align*}
(\bar b_n-b_\rho)(t,z)&=\int_{\mR^{2d}}b_n\big(t,z,(\rho^{n}_t*\varGamma_n)(z),z'\big)\rho_{t}^{n}(z')\dif z'
-\int_{\mR^{2d}}b\big(t,z,\rho_t(z),z\big)\rho_{t}(z')\dif z'\\
&=\left\{\int_{B^c_R}b_n\big(t,z,(\rho^{n}_t*\varGamma_n)(z),z'\big)\rho_{t}^{n}(z')\dif z'
-\int_{B^c_R}b\big(t,z,\rho_t(z),z'\big)\rho_{t}(z')\dif z'\right\}\\
&+\left\{\int_{B_R}b_n\big(t,z,(\rho^{n}_t*\varGamma_n)(z),z'\big)\rho_{t}^{n}(z')\dif z'
-\int_{B_R}b\big(t,z,\rho_t(z),z'\big)\rho_{t}(z')\dif z'\right\}\\
&=:I^n_R(t,z)+J^n_R(t,z).
\end{align*}
For $I^n_R$, noting that by \eqref{DR1} and Minkowskii's inequality,
\begin{align*}
\nor I^n_R\nor_{\widetilde\mL^1}
\leq\nor h\nor_{\widetilde\mL^1}\sup_{t\in[0,T]}\left(\int_{B^c_R}\rho_{t}^{n}(z')\dif z'
+\int_{B^c_R}\rho_{t}(z')\dif z'\right),
\end{align*}
by the tightness of $(\mP_n)_{n\in\mN}$, we get
$$
\lim_{R\to\infty}\sup_{n}\nor I^n_R\nor_{\widetilde\mL^1}=0.
$$
Thus it remains to show that for each $R\geq 1$,
\begin{align}\label{Lim88}
\lim_{n\to\infty}\|J^n_R\1_Q\|_{\mL^1}=0.
\end{align}
We make the following decomposition for $J^n_R$,
\begin{align*}
J^n_R(t,z)&=\int_{B_R}\Big(b_n\big(t,z,(\rho^{n}_t*\varGamma_n)(z),z'\big)-b\big(t,z,\rho_t(z),z'\big)\Big)\rho_{t}(z')\dif z'\\
&+\int_{B_R}b_n\big(t,z,(\rho^{n}_t*\varGamma_n)(z),z'\big)\Big(\rho_{t}^{n}(z')-\rho_{t}(z')\Big)\dif z'
=:J^{1,n}_R(t,z)+J^{2,n}_R(t,z).
\end{align*}
For $J^{1,n}_R$, noting that by \eqref{DL1}, 
$$
\lim_{n\to\infty}\|(\rho^{n}*\varGamma_n-\rho)\1_Q\|_{\mL^1}=0,
$$
by \eqref{CON1} and Lemma \ref{Le511}, we have for each $\eps>0$,
$$
\lim_{n\to\infty}\cL\Big\{(t,z,z')\in Q: |b_n\big(t,z,(\rho^{n}_t*\varGamma_n)(z),z'\big)-b\big(t,z,\rho_t(z),z'\big)|\geq\eps\Big\}=0.
$$
Thus, by the uniform integrability of $J^{1,n}_R\1_Q$,
\begin{align}\label{Lim87}
\lim_{n\to\infty}\|J^{1,n}_R\1_Q\|_{\mL^1}=0.
\end{align}
For $J^{2,n}_R$, by \eqref{DR1}, \eqref{DL1} and H\"older's inequality, we have
$$
\lim_{n\to\infty}\nor J^{2,n}_R\nor_{\widetilde\mL^1}\leq 
\nor h\nor_{\widetilde\mL^{2}_t(\widetilde\mL^1_z)}\lim_{n\to\infty}\|(\rho^n-\rho)\1_{B_R}\|_{\mL^{2}_t(\mL^{1}_z)}=0,
$$
which together with \eqref{Lim87} yields \eqref{Lim88}. The proof is complete.
\end{proof}

Now we can give the proof of Theorem \ref{Th63}.

\begin{proof}[Proof of Theorem \ref{Th63}] 
Let $\mP$ be the accumulation point in Lemma \ref{Le57}. We aim to show $\mP\in\sM^{a,b}_\nu$.
More precisely, we need to show that for each $t_1>t_0\geq 0$ and bounded continuous $\cB_{t_0}$-measurable functional $G_{t_0}$,
\begin{align}\label{MM5}
\mE^{\mP}(M_{t_1} G_{t_0})=\mE^{\mP}(M_{t_0} G_{t_0}),
\end{align}
where $M_t$ is defined by \eqref{MM1}.
Below we simply write $\mE:=\mE^{\mP}$ and $\mE_n:=\mE^{\mP_n}$. 
By SDE \eqref{SDE19} and It\^o's formula, it is easy to see that
\begin{align}\label{MM55}
\mE_n(M^n_{t_1} G_{t_0})=\bE(M^n_{t_1}(Z^n_\cdot) G_{t_0}(Z^n_\cdot))=\bE(M^n_{t_0}(Z^n_\cdot) G_{t_0}(Z^n_\cdot))=\mE_n(M^n_{t_0} G_{t_0}),
\end{align}
where
$$
M^n_t(\omega):=f(\omega_t)-f(\omega_0)-\int^t_0\left(\tr(\bar a_n\cdot\nabla^2_vf)+v\cdot\nabla_xf+\bar b_n\cdot\nabla_vf\right)(r,\omega_r)\dif r,
$$
where $\bar a_n$ and $\bar b_n$ are defined by \eqref{DS88} and \eqref{DS89}.
Note that \eqref{MM55} is equivalent to
\begin{align}\label{MM4}
\mE_n\left( G_{t_0}\int^{t_1}_{t_0}\left(\tr(\bar a_n \cdot\nabla^2_vf)+v\cdot\nabla_xf+\bar b_n\cdot\nabla_vf\right)(r,\omega_r)\dif r\right)=0.
\end{align}
We use Lemma \ref{Kry0} and Lemma \ref{Le58} to show the following limits
\begin{align}\label{MM13}
\lim_{n\to\infty}\mE_n\left(G_{t_0}\int^{t_1}_{t_0}\tr(\bar a_n\cdot\nabla^2_vf)(r,\omega_r)\dif r \right)
=\mE\left(G_{t_0}\int^{t_1}_{t_0}\tr(a_{\rho}\cdot\nabla^2_vf)(r,\omega_r)\dif r\right)
\end{align}
and
\begin{align}\label{MM33}
\lim_{n\to\infty}\mE_n\left(G_{t_0}\int^{t_1}_{t_0}(\bar b_n\cdot\nabla_vf)(r,\omega_r)\dif r \right)
=\mE\left(G_{t_0}\int^{t_1}_{t_0}(b_\rho\cdot\nabla_vf)(r,\omega_r)\dif r\right).
\end{align}
We  only show \eqref{MM33} since \eqref{MM13} is completely the same. For \eqref{MM33}, it follows from the following two limits:
\begin{align}\label{MM34}
\lim_{n\to\infty}\mE_n\left(G_{t_0}\int^{t_1}_{t_0}((\bar b_n-b_\rho)\cdot\nabla_vf)(r,\omega_r)\dif r \right)=0,
\end{align}
and
\begin{align}\label{MM32}
\lim_{n\to\infty}\mE_n\left(G_{t_0}\int^{t_1}_{t_0}(b_\rho\cdot\nabla_vf)(r,\omega_r)\dif r \right)
=\mE\left(G_{t_0}\int^{t_1}_{t_0}(b_\rho\cdot\nabla_vf)(r,\omega_r)\dif r\right).
\end{align}
We first show \eqref{MM34}. Let
$$
q_0:=q_1/2,\ \ \bbp_0:=\bbp_1/2.
$$
For any $R\geq 1$, since $\frac 2{q_0}+\bba\cdot\frac1{\bbp_0}<2$,
 by Krylov's estimate \eqref{Kr0} and \eqref{DS9}, we have
\begin{align*}
\mE_n\left(\int^{t_1}_{t_0}\big|\1_{B_R}(\bar b_n-b_\rho)\big|(r,\omega_r)\dif r \right)
&\lesssim \big\|\1_{[t_0,t_1]\times B_R}(\bar b_n-b_\rho)\big\|_{\mL^{ q_0}_t(\mL^{\bbp_0}_z)}\stackrel{n\to\infty}\to 0,
\end{align*}
and by H\"older's inequality,
\begin{align*}
\mE_n\left(\int^{t_1}_{t_0}\big|\1_{B^c_R}(\bar b_n-b_\rho)\big|(r,\omega_r)\dif r \right)
&\lesssim\left(\int^{t_1}_{t_0}\mP_n(|\omega_r|\geq R)\dif r \right)^{\frac12}\left(\mE_n\int^{t_1}_{t_0}\big|\bar b_n-b_\rho\big|^2(r,\omega_r)\dif r \right)^{\frac12}\\
&\lesssim \sup_{r\in[t_0,t_1]}(\mP_n(|\omega_r|\geq R))^{\frac12}\left(\nor|\bar b_n-b_\rho|^2\nor_{\widetilde\mL^{ q_0}_t(\widetilde\mL^{\bbp_0}_z)} \right)^{\frac12}\\
&\leq \sup_{r\in[t_0,t_1]}(\mP_n(|\omega_r|\geq R))^{\frac12}\nor\bar b_n-b_\rho\nor_{\widetilde\mL^{ q_1}_t(\widetilde\mL^{\bbp_1}_z)}\\
&\lesssim \sup_{r\in[t_0,t_1]}(\bP(|Z^n_r|\geq R))^{\frac12},
\end{align*}
which converges to zero uniformly in $n$ as $R\to\infty$ by the proof of Lemma \ref{Le56}.
Thus we obtain \eqref{MM34}. Next we look at \eqref{MM32}. Let
$$
b^m_\rho(t,z):=b_\rho(t,\cdot)*\varGamma_m(z).
$$
As above, by Krylov's estimate \eqref{Kr0}, we have
\begin{align}\label{MM2}
\lim_{m\to \infty}\sup_n\mE_n\left(\int^{t_1}_{t_0}|(b^m_\rho-b_\rho)\cdot\nabla_vf|(r,\omega_r)\dif r\right)=0.
\end{align}
Moreover, for fixed $m\in\mN$, since $\omega\to G_{t_0}(\omega)\int^{t_1}_{t_0}(b^m_\rho\cdot\nabla_vf)(r,\omega_r)\dif r$ 
is a bounded continuous functional, we have
$$
\lim_{n\to\infty}\mE_n\left(G_{t_0}\int^{t_1}_{t_0}(b^m_\rho\cdot\nabla_vf)(r,\omega_r)\dif r\right)
=\mE\left( G_{t_0}\int^{t_1}_{t_0}(b^m_\rho\cdot\nabla_vf)(r,\omega_r)\dif r\right),
$$
which together with \eqref{MM2} yields \eqref{MM32}. 
Thus, by taking limits $n\to\infty$ for \eqref{MM4}, we obtain
$$
\mE\left( G_{t_0}\int^{t_1}_{t_0}\left(\tr(a_\rho \cdot\nabla^2_vf)+v\cdot\nabla_xf+b_\rho\cdot\nabla_vf\right)(r,\omega_r)\dif r\right)=0,
$$
which implies \eqref{MM5}. As for the regularity \eqref{DL131}, it follows by \eqref{DL11} and \eqref{DL1}. The proof is complete.
\end{proof}
\br\label{Re611}\rm
Suppose that the coefficients do not depend on the density, that is, $b(t,z,r,z')$ and $a(t,x,r,z')$ are independent of $r$, 
then we can drop the assumption $q_1\in(2,4)$ in {\bf (H$_1$)} but require $\tfrac1{\bbp_1}+\frac{\1}{q_1}<\tfrac1{\bb2}$ (see Remark \ref{Re44}). 
In this case, we can use the Kyrlov estimate \eqref{Kr1} and follow the same argument as in \cite{RZ21} to take the limits.
\er
\section{Well-posedness of generalized martingale problems}
In this section we show the well-posedness for a class of generalized martingale problems of SDE \eqref{SDE00}
when diffusion matrix does not depend on the distribution, but may be discontinuous in position variable $x$.
Instead of {\bf (H$_1$)} and {\bf (H$_2$)}, we suppose that {\bf (H$_3$)} holds. 
Fix $T>0$ and $f\in  C^\infty_c(\mR^{2d})$. Consider the following backward Kolmogorov equation:
\begin{align}\label{BPDE}
\p_t u+\tr(a\cdot \nabla^2_v u)+v\cdot\nabla_x u+b_{\rho}
\cdot\nabla_v u=f,\ \ u(T)=0,
\end{align}
where for a family of  density functions $\rho_t(z)$ in $\mR^{2d}$,
$$
b_{\rho}(t,z):=\int_{\mR^{2d}}b(t,z,\rho_t(z), z')\rho_t(z')\dif z'.
$$
Since $b_\rho$ is bounded measurable and $a$ does not depend on $v$, by reversing the time variable and Theorem \ref{Th42},
there is a unique weak solution $u^f_\rho$ to \eqref{BPDE} with
$$
u^f_\rho\in \widetilde\sV_T\cap C_b([0,T]\times\mR^{2d}).
$$
Here the continuity of $u^f_\rho$ is crux for the following notion of generalized martingale solutions.
\bd[Generalized martingale problem]\label{MP81}
Let $s\geq 0$ and $\nu\in\cP(\mR^{2d})$. A probability measure $\mP\in\cP(\mC)$ is called a generalized martingale solution of DDSDE \eqref{SDE00} starting from $\nu$ 
at time $s$ if
\begin{enumerate}[(i)]
\item $\mP\circ \omega^{-1}_s=\nu$, and for Lebesgue almost all $t\geq s$, $\mP\circ \omega_t^{-1}(\dif z)=\rho_t(z)\dif z$.

\item For any $T>s$ and $f\in C^\infty_c(\mR^{2d})$, the process
\begin{align}\label{MM17}
M_t:=u^f_\rho(t,\omega_t)-u^f_\rho(s,\omega_s)-\int^t_sf(\omega_r)\dif r,\ \ t\in[s,T],
\end{align}
is a $\cB_t$-martingale with respect to $\mP$, where $u^f_\rho$ is the unique solution of \eqref{BPDE}.
\end{enumerate}
The set of all the generalized martingale solutions $\mP\in\cP(\mC)$  with initial distribution $\nu$ at time $s$ is denoted by $\widetilde\cM^{a,b}_{s,\nu}$.
\ed
\br\rm
The above notion of martingale solutions was introduced in its most general form in \cite[Chapter 4]{EK86}.
It should be noticed that if \eqref{BPDE} has a $C^2$-solution $u$, then by It\^o's formula, any weak solution of \eqref{SDE00} must be a generalized martingale solution. In general, 
these two notions are not equivalent.
\er

Now we can give 
\begin{proof}[Proof of Theorem \ref{Th15}]
{\bf (Existence)} Without loss of generality we assume $s=0$.
As in Section \ref{sec5}, we consider the approximating SDE \eqref{SDE19}. Let $\mP\in\cP(\mC)$ be an accumulation point of $(\mP_n)_{n\in\mN}$.
In particular, by Lemma \ref{Le56},
$$
\mP_n\circ \omega_t^{-1}(z)=\rho^n_t(z)\dif z,\ \ \mP\circ \omega_t^{-1}(z)=\rho_t(z)\dif z.
$$
In the following we use the convention before Lemma \ref{Le511}.
Recall $a_n(t,x)=a(t,\cdot)*\varGamma_n(x)$ and $\bar b_n$ being defined by \eqref{DS89}.
Since for each $\varphi\in C^\infty_c(\mR^{2d})$, by It\^o's formula,
$$
\p_t \int \varphi\rho^n_t=\int \tr (a_n\cdot\nabla^2_v\varphi)\rho^n_t+\int (v\cdot\nabla_x\varphi) \rho^n_t+\int (\bar b_n \cdot\nabla_v \varphi) \rho^n_t,
$$
by Lemmas \ref{Le57} and \ref{Le58} and taking limits, it is easy to see that
$$
\p_t \int \varphi\rho_t=\int\tr (a\cdot\nabla^2_v\varphi)\rho_t+\int (v\cdot\nabla_x\varphi) \rho_t+\int (b_{\rho} \cdot\nabla_v \varphi) \rho_t.
$$
In other words, $\rho_t$ solves \eqref{FPK1} in the distributional sense.
Fix $T>0$ and $f\in C^\infty_c(\mR^{2d})$. Let $u_n\in L^\infty_T(C^\infty_b(\mR^{2d}))$ be the unique smooth solution of the following backward PDE:
\begin{align}\label{CC10}
\p_t u_n+\tr(a_n\cdot \nabla^2_v u_n)+v\cdot\nabla_x u_n+\bar b_n\cdot\nabla_v u_n=f,\ \ u_n(T)=0.
\end{align}
By Lemma \ref{Le58} and Theorem \ref{Th42}, we have
for any bounded $Q\subset(0,T)\times\mR^{2d}$,
\begin{align}\label{CC1}
\lim_{n\to\infty}\sup_{(t,z)\in Q}|u_{n}(t,z)-u^f_\rho(t,z)|=0,
\end{align}
where $u^f_\rho$ is the unique weak solution of PDE \eqref{BPDE}.
Now we show that for each $T\geq t_1>t_0\geq 0$ and bounded continuous $\cB_{t_0}$-measurable functional $G_{t_0}$,
\begin{align}\label{MM35}
\mE^{\mP}(M_{t_1} G_{t_0})=\mE^{\mP}(M_{t_0} G_{t_0}),
\end{align}
where $M_t$ is defined by \eqref{MM17}.
By SDE \eqref{SDE19}, It\^o's formula and \eqref{CC10}, it is easy to see that
\begin{align}\label{MM45}
\mE^{\mP_n}(M^n_{t_1} G_{t_0})=\mE^{\mP_n}(M^n_{t_0} G_{t_0}),
\end{align}
where	
$$
M^n_t:=u_n(t,\omega_t)-u_n(0,\omega_0)-\int^t_0f(\omega_r)\dif r,
$$
By taking limits for both sides of \eqref{MM45} and using the pointwise convergence \eqref{CC1}, we obtain \eqref{MM35}.

\medskip

{\bf (Uniqueness)} We divide the proof into three steps.

({\sc Step 1}) First of all, we show the uniqueness for linear SDE, i.e., $b(t,z,r,z')=b(t,z)$ does not depend on variables $r,z'$.
Let $\mP_1,\mP_2\in\widetilde\sM^{a,b}_{s,\nu}$ be two solutions of the generalized martingale problem.  Fix $T>s$ and $f\in C^\infty_c(\mR^{2d})$. 
Let $u^f$ be the unique weak solution of (see Theorem  \ref{Th42}),
$$
\p_t u^f+\tr(a\cdot \nabla^2_v u^f)+v\cdot\nabla_x u^f+b\cdot\nabla_v u^f=f,\ \ u^f(T)=0.
$$
By Definition \ref{MP81} and $u^f(T)=0$, we have
$$
\int_{\mR^{2d}}u^f(s,z)\nu(\dif z)=-\mE^{\mP_i}\int_s^Tf(\omega_r)\dif r,\quad i=1,2,
$$
which means that for each $T>s$,
\begin{align*}
\int_s^T\mE^{\mP_1} f(\omega_r)\dif r=\int_s^T\mE^{\mP_2} f(\omega_r)\dif r.
\end{align*}
Hence, for any $f\in C^\infty_c(\mR^{2d})$
\begin{align}\label{SD91}
\mE^{\mP_1} f(\omega_T)=\mE^{\mP_2} f(\omega_T),\ \ \forall T>s.
\end{align}
From this, by a standard way (see Theorem 4.4.2 in \cite{EK86}), we derive that
$$
\mP_1=\mP_2.
$$
Indeed, it suffices to prove the following claim by induction: 

{\bf (C$_n$)} for given $n\in\mN$, and for any $s\leq t_1<t_2<t_n<T$ and strictly positive and bounded measurable functions $g_1,\cdots, g_n$ on $\mR^{2d}$,
\begin{align}\label{GK1}
\mE^{\mP_1}(g_1(\omega_{t_1})\cdots g_n(\omega_{t_n}))=\mE^{\mP_2}(g_1(\omega_{t_1})\cdots g_n(\omega_{t_n})).
\end{align}
By \eqref{SD91}, {\bf (C$_1$)} holds.
Suppose now that {\bf (C$_n$)} holds for some $n\geq 2$. For simplicity we write
$$
\eta:=g_1(\omega_{t_1})\cdots g_n(\omega_{t_n})>0,
$$
and for $i=1,2$, we define new probability measures
$$
\dif\widetilde\mP_i:=\eta\dif\mP_i/(\mE^{\mP_i}\eta)\in\cP(\mD),\ \ \widetilde\nu_i:=\widetilde\mP_i\circ \omega^{-1}_{t_n}\in\cP(\mR^{2d}).
$$
Now we show
$$
\widetilde\mP_i\in\widetilde\cM^{a,b}_{t_n;\widetilde\nu_i},\ \ i=1,2.
$$
For any $f\in C^\infty_c(\mR^{2d})$, let 
$$
M_t:=u^f(t,\omega_t)-u^f(t_n,\omega_{t_n})-\int^t_{t_n}f(\omega_r)\dif r,\ \ t\in[t_n,T].
$$
We only need to prove that for any $T\geq t'>t\geq t_n$ and bounded $\cB_t$-measurable $\xi$,
$$
\mE^{\widetilde\mP_i}\left(M_{t'}\xi\right)=\mE^{\widetilde\mP_i}\left(M_t\xi\right)\Leftrightarrow \mE^{\mP_i}(M_{t'}\xi\eta)=\mE^{\mP_i}(M_t\xi\eta),
$$
which follows from $\mP_i\in\widetilde\cM^{a,b}_{s;\nu}$, $i=1,2$. Thus, by induction hypothesis and \eqref{SD91},
$$
\widetilde\nu_1=\widetilde\nu_2\Rightarrow\widetilde\mP_1\circ \omega^{-1}_{t_{n+1}}=\widetilde\mP_2\circ \omega^{-1}_{t_{n+1}},\ \ \forall T\geq t_{n+1}>t_n.
$$
which in turn implies that {\bf (C$_{n+1}$)} holds.

({\sc Step 2})  For general nonlinear SDE, we use Girsanov's trasformation method
(see \cite{RZ21, HZZZ21}). Without loss of generality, we assume $s=0$. Let $\mP_1,\mP_2\in\widetilde\sM^{a,b}_{0,\nu}$ be two solutions of the generalized martingale problem.
Let 
$$
a_n:=a*\varGamma_n,\ \ b_n:=b*\varGamma_n,\ \ \rho^n_0:=\rho_0*\varGamma_n,
$$
and for $\mP_i\circ\omega_t^{-1}(\dif z)=\rho^i_t(z)\dif z$,
\begin{align}\label{LK9}
\bar b^{(i)}_n(t,z):=\int_{\mR^{2d}}b_n(t,z, (\rho^i_t*\varGamma_n)(z), z')\rho^i_t(z')\dif z'.
\end{align}
Note that for any $T>0$,
$$
a_n, \bar b^{(i)}_n\in L^\infty_T(C^\infty_b(\mR^{2d})),\ \ \|a_n\|_\infty\leq\|a\|_\infty,\ \ \|\bar b^{(i)}_n\|_\infty\leq \|b\|_\infty.
$$
We consider the following approximation of linearized SDEs: for $i=1,2$,
\begin{align}\label{Asde0}
\dif X^{i,n}_t=V^{i,n}_t\dif t,\
\dif V^{i,n}_t=\bar b^{(i)}_n(t,Z^{i,n}_t)\dif t+\sqrt{2a_n}(t,X^{i,n}_t)\dif W_t,
\end{align}
where $Z^{i,n}_0$ is an $\sF_0$-measurable random variable and has smooth density $\rho^n_0$.
Since by Lemma \ref{Le511},
$$
\bar b^{(i)}_n(t,z)\stackrel{n\to\infty}{\to} b_{\rho^i}(t,z)=\int_{\mR^{2d}}b(t,z, \rho^i_t(z), z')\rho^i_t(z')\dif z'\mbox{ in Lebesgue measure $\dif t\dif z$,}
$$
as in the proof of the existence part, 
and due to the uniqueness of linear equations, for $i=1,2$, 
the law of $Z^{i,n}$ weakly converges to $\mP_i$ as $n\to\infty$. In particular, for any $\varphi\in C_b(\mR^{2d})$,
\begin{align}\label{AX1}
\bE\varphi(Z^{i,n}_t)\stackrel{n\to\infty}{\to} \mE^{\mP_i}\varphi(\omega_t),\ \ i=1,2.
\end{align}
Since $a_n, \bar b^{(i)}_n\in L^\infty_T(C^\infty_b(\mR^{2d}))$ and $\rho^n_0\in C^\infty_b$,
it is well known that the density $\rho^{i,n}_t$ of $Z^{i,n}_t$ is smooth and uniquely solves the following FPKE:
$$
\p_t \rho^{i,n}_t=\div_v(a_n\cdot\nabla_v\rho^{i,n}_t)-v\cdot\nabla_x \rho^{i,n}_t+\div_v(\bar b^{(i)}_n\rho^{i,n}_t),
$$
where we have used that $a_n$ does not depend on the velocity variable $v$. By Theorem \ref{Th44} with $b_1=0$, 
$b_2=\bar b^{(i)}_n$ and $g=0$, for any $T>0$, there is a constant $C>0$ such that for all $n\in\mN$,
\begin{align}\label{AX2}
\|\1_{[0,T]}\rho^{i,n}\|_{\mL^\infty}\lesssim_C \|\rho^{n}_0\|_{C^1_b}\lesssim_C \|\rho_0\|_{C^1_b},\ \ i=1,2.
\end{align}

({\sc Step 3})  To perform the Girsanov transformation, for $i=1,2$ and $n\in\mN$, we define
\begin{align}\label{LK91}
H^{(i)}_n(s,z):=(\sqrt{2a_n})^{-1}(s,x)\cdot\bar b^{(i)}_n(s,z),\ \ z=(x,v),
\end{align}
and
$$
A^{i,n}_t:=\exp\left\{-\int^t_0H^{(i)}_n(s,Z^{i,n}_s)\dif W_s-\frac12\int^t_0|H^{(i)}_n(s,Z^{i,n}_s)|^2\dif s\right\}.
$$
Fix $T>0$. Since  for some $C>0$ independent of $n$,
\begin{align}\label{AS9}
\|H^{(i)}_n\|_\infty\leq C,\ \ i=1,2,
\end{align}
by Girsanov's theorem, for $i=1,2$, under the new probability measure $Q^{i,n}:=A^{i,n}_T\bP$,
$$
\widetilde W^{i,n}_t:=\int^t_0H^{(i)}_n(s,Z^{i,n}_s)\dif s+W_t,\ t\in[0,T],
$$
is still a Brownian motion, and
$$
\dif X^{i,n}_t=V^{i,n}_t\dif t,\ \dif V^{i,n}_t=\sqrt{2a_n}(t,X^{i,n}_t)\dif \widetilde W^{i,n}_t.
$$
Since $\sqrt{a_n}$ is Lipschitz, the above SDE admits a unique weak solution. Thus
$$
Q^{1,n}\circ (Z^{1,n})^{-1}=Q^{2,n}\circ (Z^{2,n})^{-1},
$$
equivalently, for any nonnegative functional $G$ on $\mC_T:=C([0,T];\mR^{2d})$,
$$
\bE(G(Z^{1,n}_\cdot)A^{1,n}_T)=\bE(G(Z^{2,n}_\cdot)A^{2,n}_T).
$$
In particular, for any $\varphi\in C_b(\mR^{2d})$,
\begin{align*}
\bE\varphi(Z^{1,n}_T)=\bE(\varphi(Z^{2,n}_T)Y^n_T),
\end{align*}
where
$$
Y^n_T:=\exp\left\{\int^T_0(H^{(1)}_n-H^{(2)}_n)(s,Z^{2,n}_s)\dif W_s+\frac12\int^T_0(|H^{(1)}_n|^2-|H^{(2)}_n|^2)(s,Z^{2,n}_s)\dif s\right\}.
$$
Thus, by H\"older's inequality,
\begin{align}
|\bE\varphi(Z^{1,n}_T)-\bE\varphi(Z^{2,n}_T)|&=|\bE(\varphi(Z^{2,n}_T)(Y^n_T-1))|
\leq\|\varphi(Z^{2,n}_T)\|_{L^2(\Omega)}\|Y^n_T-1\|_{L^2(\Omega)}\no\\
&\lesssim\Big(\|\varphi\|_\infty\wedge\|\varphi\|_{\mL^2_z}\Big)\|Y^n_T-1\|_{L^2(\Omega)},\label{KB1}
\end{align}
where the implicit constant is independent of $n$, and we have used \eqref{AX2} to derive that
$$
\|\varphi(Z^{2,n}_T)\|_{L^2(\Omega)}=\left(\int_{\mR^{2d}}|\varphi(z)|^2\rho^{2,n}_T(z)\dif z\right)^{1/2}\lesssim\|\varphi\|_{\mL^2_z}.
$$
On the other hand, by It\^o's formula, we have
\begin{align*}
Y^n_T-1&=\int^T_0Y^n_s(H^{(1)}_n-H^{(2)}_n)(s,Z^{2,n}_s)\dif W_s\\
&+\frac12\int^T_0Y^n_s(|H^{(1)}_n|^2-|H^{(2)}_n|^2)(s,Z^{2,n}_s)\dif s\\
&+\frac12\int^T_0Y^n_s|H^{(1)}_n-H^{(2)}_n|^2(s,Z^{2,n}_s)\dif s.
\end{align*}
By BDG's inequality and \eqref{AS9}, we have
\begin{align*}
\bE|Y^n_T-1|^2\lesssim\int^T_0\bE|Y^n_s-1|^2\dif s+\int^T_0\bE\big(|H^{(1)}_n-H^{(2)}_n|^2(s,Z^{2,n}_s)\big)\dif s,
\end{align*}
which implies by Gronwall's inequality that
$$
\bE|Y^n_T-1|^2\lesssim\int^T_0\bE\big(|H^{(1)}_n-H^{(2)}_n|^2(s,Z^{2,n}_s)\big)\dif s.
$$
Note that by \eqref{LK9}, \eqref{LK91} and {\bf (H$_3$)},
\begin{align*}
|H^{(1)}_n-H^{(2)}_n|(t,z)\leq\|(\sqrt{2a_n})^{-1}\|_\infty|\bar b^{(1)}_n-\bar b^{(2)}_n|(t,z)
\lesssim|(\rho^1_t-\rho^2_t)*\varGamma_n|(z)+\|\rho^1_t-\rho^2_t\|_{\mL^1_z}.
\end{align*}
Therefore, by \eqref{AX2} again,
\begin{align}
\bE|Y^n_T-1|^2&\lesssim\int^T_0\bE|(\rho^1_s-\rho^2_s)*\varGamma_n|^2(Z^{2,n}_s)\dif s+\int^T_0\|\rho_s^1-\rho_s^2\|^2_{\mL^1_z}\dif s\no\\
&\lesssim\int^T_0\|(\rho^1_s-\rho^2_s)*\varGamma_n\|^2_{\mL^2_z}\dif s+\int^T_0\|\rho_s^1-\rho_s^2\|^2_{\mL^1_z}\dif s\no\\
&\lesssim\int^T_0\|\rho^1_s-\rho^2_s\|^2_{\mL^2_z}\dif s+\int^T_0\|\rho_s^1-\rho_s^2\|^2_{\mL^1_z}\dif s\no\\
&\lesssim\int^T_0\Big(\|\rho^1_s-\rho^2_s\|^2_{\mL^2_z}\vee\|\rho_s^1-\rho_s^2\|^2_{\mL^1_z}\Big)\dif s.\label{KB2}
\end{align}
Combining \eqref{KB1} and \eqref{KB2} and by \eqref{AX1}, we obtain that for all $\varphi\in C_b(\mR^{2d})$,
\begin{align*}
&|\mE^{\mP_1}\varphi(\omega_T)-\mE^{\mP_2}\varphi(\omega_T)|=\lim_{n\to\infty}|\bE\varphi(Z^{1,n}_T)-\bE\varphi(Z^{2,n}_T)|\\
&\qquad\lesssim\Big(\|\varphi\|_\infty\wedge\|\varphi\|_{\mL^2_z}\Big)
\left(\int^T_0\Big(\|\rho_s^1-\rho_s^2\|^2_{\mL^2_z}\vee \|\rho^1_s-\rho^2_s\|^2_{\mL^1_z}\Big)\dif s\right)^{1/2}.
\end{align*}
Since for Lebesgue almost all $T>0$, $\mP_i\circ\omega_T^{-1}(\dif z)=\rho_T^i(z)\dif z$, the above inequality means that for Lebesgue almost all $T>0$,
\begin{align}\label{KB4}
\|\rho_T^1-\rho_T^2\|_{\mL^2_z}\vee \|\rho^1_T-\rho^2_T\|_{\mL^1_z}\lesssim
\left(\int^T_0\Big(\|\rho_s^1-\rho_s^2\|^2_{\mL^2_z}\vee \|\rho^1_s-\rho^2_s\|^2_{\mL^1_z}\Big)\dif s\right)^{1/2},
\end{align}
which in turn implies that $\rho^1_T=\rho^2_T$ by Gronwall's inequality. Finally, as usual we use the uniqueness for linear equations in {\sc Step 1} to derive $\mP_1=\mP_2$. The proof is complete.
\end{proof}

\br\label{Re71}\rm
(i) From the proof of the existence part, it is easy to see that the unique generalized martingale solution $\mP\in\widetilde\sM^{a,b}_{0,\nu}$ also belongs to $\sM^{a,b}_\nu$, i.e., it is also a classical martingale solution.

(ii) If $b$ does not depend on the density variable $r$, then we can drop the assumption on the initial distribution $\nu$ since in this case the boundedness \eqref{AX2} is not needed in the proof and only the total variational norm $\|\rho^1_s-\rho^2_s\|^2_{\mL^1_z}$ is used in \eqref{KB4}.
\er

{\bf Acknowledgement}. The author is grateful to Zimo Hao for quite useful conversations.

\end{document}